\documentclass[11pt]{amsart}

\usepackage{amsmath,amssymb}

\def\COMMENT#1{}
\def\TASK#1{}

\def\noproof{{\unskip\nobreak\hfill\penalty50\hskip2em\hbox{}\nobreak\hfill%
        $\square$\parfillskip=0pt\finalhyphendemerits=0\par}\goodbreak}
\def\endproof{\noproof\bigskip}
\newdimen\margin   
\def\textno#1&#2\par{%
    \margin=\hsize
    \advance\margin by -4\parindent
           \setbox1=\hbox{\sl#1}%
    \ifdim\wd1 < \margin
       $$\box1\eqno#2$$%
    \else
       \bigbreak
       \hbox to \hsize{\indent$\vcenter{\advance\hsize by -3\parindent
       \sl\noindent#1}\hfil#2$}%
       \bigbreak
    \fi}

\newtheorem{firstthm}{Proposition}
\newtheorem{thm}[firstthm]{Theorem}

\newtheorem{lemma}[firstthm]{Lemma}
\newtheorem{cor}[firstthm]{Corollary}

\newtheorem{defin}[firstthm]{Definition}
\newtheorem{conj}[firstthm]{Conjecture}

\begin{document}
\title{Edge-disjoint Hamilton Cycles in Random Graphs}
\date{\today}
\author{Fiachra Knox, Daniela K\"uhn and Deryk Osthus}
\thanks {The authors were supported by the EPSRC, grant no.~EP/F008406/1. D. K\"uhn was also supported by the ERC, grant no.~258345.} 

\begin{abstract} We show that provided $\log^{50} n/n \leq p \leq 1 - n^{-1/4}\log^9 n$ we can with high probability find a collection of 
$\lfloor \delta(G)/2 \rfloor$ edge-disjoint Hamilton cycles in $G \sim G_{n, p}$, 
plus an additional edge-disjoint matching of size $\lfloor n/2 \rfloor$ if $\delta(G)$ is odd. This is clearly optimal and confirms, for the above range of $p$, 
a conjecture of Frieze and Krivelevich.
\end{abstract}

\maketitle
\section{Introduction} \label{introduction}

The question of how many edge-disjoint Hamilton cycles one can pack into a given graph $G$ has a long history, but
many of the main questions remain open. More precisely,
we say that a graph $G$ has \emph{property} $\mathcal{H}$ if $G$ contains 
$\lfloor \delta(G)/2 \rfloor$ edge-disjoint Hamilton cycles, together with 
an additional edge-disjoint optimal matching if $\delta(G)$ is odd.
(For a graph $G$ on $n$ vertices, call a matching $M$ in $G$ \emph{optimal} if $|M| = \lfloor n/2 \rfloor$.)
An old construction by Walecki (see e.g.~\cite{abs,lucas}) shows that the complete graph $K_n$  
has a Hamilton decomposition if $n$ is odd
(here a Hamilton decomposition of a graph $G$ is a set of edge-disjoint Hamilton cycles which together cover 
all of the edges of $G$).
It is well known that more generally $K_n$ has property $\mathcal{H}$ (see e.g.~\cite{Wallis}).
We discuss some further results at the end of the introduction.

The question was first investigated in a probabilistic setting by Bollob\'as and Frieze in the 1980's (see below). Later,
Frieze and Krivelevich~\cite{FK08} made the striking conjecture that whp $G_{n,p}$ has property $\mathcal{H}$ for any $p = p(n)$.
Here $G_{n,p}$ denotes a binomial random graph on $n$ vertices where every edge is present with probability $p$,
and we say that a property of a random graph on $n$ vertices holds \emph{whp} if the probability 
that it holds tends to 1 as $n$ tends to infinity.

\begin{conj}[Frieze and Krivelevich \cite{FK08}] \label{mainconj}
For any $p=p(n)$, whp a binomial random graph $G \sim G_{n,p}$ has property $\mathcal{H}$.
\end{conj}
Our main result confirms this conjecture as long as $p$ is not too small and not too large.

\begin{thm} \label{HamDecomp} Let $\frac{\log^{50} n}{n} \leq p \leq 1 - n^{-1/4}\log^9 n$. Then whp $G_{n, p}$ has
property $\mathcal{H}$.
\end{thm}
We deduce Theorem \ref{HamDecomp} from a purely deterministic result (Theorem~\ref{PRproof}) which states that every graph which satisfies certain
pseudorandomness conditions and which is close (but not too close) to being regular has property $\mathcal{H}$. 
Our proof shows that if both $\delta(G_{n, p})$ and $n$ are odd, then the optimal matching guaranteed by property $\mathcal{H}$ 
can be chosen to cover the vertex $x$ of minimum degree. So the Hamilton cycles and this matching together cover all the edges at $x$.

We now discuss some related results on packing edge-disjoint Hamilton cycles in random graphs.
Most of these actually consider the slightly weaker property $\mathcal{H}_\delta$ of containing $\lfloor \delta(G)/2 \rfloor$
edge-disjoint Hamilton cycles. 
Let $G_{n,m}$ denote a random graph chosen uniformly from the set of graphs on $n$ vertices with $m$ edges.
Bollob\'as and Frieze~\cite{BF85} showed that an analogue of Conjecture~\ref{mainconj} holds in the model $G_{n,m}$
when $2m/n=pn \le \log n+ O(\log \log n)$.
In this range of $m$, whp $\delta(G_{n,m})$ is bounded.
Their result generalizes a result of Bollob\'as~\cite{Bol84}, which gave a `hitting time' result for the property that
$G_{n,m}$ contains a Hamilton cycle.
Frieze and Krivelevich~\cite{FK08} extended the range of $p$ and showed that property $\mathcal{H}_\delta$ holds for all $p$ with 
$pn=(1+o(1)) \log n$. 
Recently, this was further extended by Ben-Shimon, Krivelevich and Sudakov~\cite{BKS} to all $p \le 1.02 \log n /n$.
In this range whp $G_{n,p}$ is far from being regular. 
In fact, as noted in~\cite{BKS} whp we have $\delta(G_{n,p}) \le  np/300$, and so only a small fraction of the edges of $G_{n,p}$ are contained in the Hamilton cycles guaranteed by property $\mathcal{H}_\delta$.
Using different ideas than in the present paper, very recently Krivelevich and Samotij~\cite{KrSa} were able to cover the range 
$\log n/n \le p \le n^{-1+\varepsilon}$. 
Finally, the `very dense' case was settled by K\"uhn and Osthus~\cite{KOappl}, who covered the range when $p \ge 2/3$.
So altogether, these results cover the entire range of probabilities $p$ and show that whp $\mathcal{H}_\delta$ holds for any $p$.
We emphasize that the `main' range of probabilities is covered by~\cite{KrSa} and the current paper.

As soon as $pn \gg \log n$, whp we have $\delta(G_{n,p}) =  (1+o(1))np$. 
So in this range, $G_{n, p}$ is close to being regular and so the above results yield an `approximate' Hamilton decomposition of 
$G_{n,p}$, i.e.~a set of edge-disjoint Hamilton cycles covering almost all edges of $G_{n,p}$.
Such an approximate result for constant $p$ was first proved by Frieze and Krivelevich~\cite{FK05}.
In~\cite{AHDoRG}, we extended this approximate result to all $pn \gg \log n$.
This was also proved independently by Krivelevich (personal communication).

As mentioned earlier, not many graphs are known which have property~$\mathcal{H}$.
Some examples are discussed in the surveys~\cite{abs,bondy} -- mainly these deal with 
classes of graphs which have a high degree of symmetry. A more recent result in the area is due to 
Alspach, Bryant and Dyer~\cite{abd}, who showed that Paley graphs have property~$\mathcal{H}$.
Kim and Wormald~\cite{KW01} proved that whp property $\mathcal{H}$ holds for random regular graphs of fixed degree. 
Nash-Williams~\cite{Nash-Williams71b} conjectured that any $d$-regular graph on $n$ vertices where $d \geq n/2$ has $\lfloor d/2 \rfloor$ 
edge-disjoint Hamilton cycles. (So such graphs would have property $\mathcal{H}$ if $d$ is even.) 
Approximate versions of this conjecture were proved in~\cite{CKO, KLO}.

A very general result in this direction was obtained very recently by K\"uhn and Osthus~\cite{monster}, who showed that
every regular `robustly expanding' digraph of linear degree has a decomposition into edge-disjoint Hamilton cycles.
The initial motivation for this result was that it implies that every (large) regular tournament has a Hamilton decomposition,
which proved a long-standing conjecture of Kelly.
However, as observed in~\cite{KOappl}, the main result of~\cite{monster} has a number of further applications. 
For instance, it can be used to show an analogue of Conjecture~\ref{mainconj} for random tournaments, which confirms a conjecture of Erd\H{o}s.
Similarly, the result in~\cite{KOappl} which confirms property $\mathcal{H}_\delta$ for $G_{n,p}$ in the `very dense' case is also derived
from the main result of~\cite{monster}.

Hypergraph versions of Conjecture~\ref{mainconj}
were also recently considered by Frieze and Krivelevich~\cite{FK10},
Frieze, Krivelevich and Loh~\cite{FKL10} and Bal and Frieze~\cite{BalFrieze}.

A related line of research was initiated by Glebov, Krivelevich and Szab\'o~\cite{GKS}, 
who investigated the `dual' problem of finding the smallest set of Hamilton cycles which together
cover all edges of $G_{n,p}$. As $G_{n,p}$ is whp not regular, this means that these Hamilton cycles will not be edge-disjoint
and the best result one can hope for is a set of $\lceil \Delta(G_{n,p})/2 \rceil$ such Hamilton cycles.
The main result in~\cite{GKS} shows that this bound is approximately true if $p \ge n^{-1+\varepsilon}$.
Subsequently, Hefetz, K\"uhn, Lapinskas and Osthus~\cite{HKLO} achieved an exact version of this result for a similar range
of $p$ as in Theorem~\ref{HamDecomp}.
\begin{thm}\label{thm:main-result}
Let $\frac{\log^{117}n}{n}\le p\le1-n^{-1/8}$. Then whp
the edges of $G_{n,p}$ can be covered by $\lceil \Delta(G)/2 \rceil$
Hamilton cycles.
\end{thm}
For $p$ as above, this proves a conjecture of~\cite{GKS}.
The crucial ingredient in the proof of Theorem~\ref{thm:main-result} is the main technical result of the current paper (Lemma~\ref{MiniHamDecomp}).
Roughly speaking, given a graph $H$ satisfying rather weak pseudorandomness conditions and a pseudorandom graph $G_1$, which is allowed to be surprisingly sparse
compared to $H$, Lemma~\ref{MiniHamDecomp} guarantees a set of edge-disjoint Hamilton cycles covering all edges of $H$.

Our proof of Theorem~\ref{HamDecomp} exploits the fact that in the given range of $p$, there is a small but significant gap between 
the two smallest degrees of $G_{n,p}$. This is not necessarily the case when $np/\log n \to \infty$ very slowly
(see e.g.~\cite{Bollobasbook}). This is one of the reasons for our restriction on the range of $p$ in Theorem~\ref{HamDecomp}.
Our method is based on a new technique of iterative improvements which is likely to have further applications.
In particular, this idea was an essential feature of the proof in~\cite{monster}.

\section{Sketch of proof and organisation of the paper} \label{sketch}

Recall that in~\cite{AHDoRG} we proved an approximate version of Conjecture~\ref{mainconj}, which 
finds a set of edge-disjoint Hamilton cycles covering all but an $\eta$ proportion of the edges of $G_{n,p}$.
Roughly speaking, the proof of this approximate result proceeds as follows: 
\begin{itemize}
\item[(1)] First we choose $0 < \varepsilon \ll \eta$ and remove a random subgraph $G_{thin}$ of density $\varepsilon p$ from $G_{n,p}$, and
call the remaining graph $G_{dense}$. $G_{thin}$ will be a pseudorandom graph, which is important in 
Step 4 below.
\item[(2)] Next we apply Tutte's theorem to find an $r$-regular subgraph $G_{reg}$ of $G_{dense}$, with $r=(1-\varepsilon)np$. 
\item[(3)] Using a counting argument, one can show that most $2$-factors of $G_{reg}$ have few cycles.
We remove such $2$-factors one by one to obtain a set of $(1-\varepsilon)r/2$ edge-disjoint $2$-factors $F_i$ which cover 
most edges of $G_{reg}$ (and thus of $G_{n,p}$).
\item[(4)] We now use the edges of $G_{thin}$ to transform each $F_i$ into a Hamilton cycle. 
More precisely, for each $i$ in turn, we swap some edges between $G_{thin}$ and $F_i$ in such a way that
the modified $2$-factor is in fact a Hamilton cycle. For the next $2$-factor, we use the `current' version of $G_{thin}$.
The fact that $F_i$ has few cycles means that we do not need to swap out many edges of $G_{thin}$ for this, and so $G_{thin}$ retains the pseudorandomness properties which allow us to perform this step.
\end{itemize}
Our approach to Theorem~\ref{HamDecomp} uses the above ideas, but some obvious obstacles arise.
For simplicity, assume that $n$  and $\delta(G_{n,p})$ are both even in what follows.
We can use Tutte's theorem to show that $G_{n,p}$ has an $r$-regular subgraph $G_{reg}$ with $r=\delta(G_{n,p})$
(this is harder than Step 2 but not too difficult for our range of $p$).
So we can decompose $G_{reg}$ into $r/2$ $2$-factors $F_i$ and would like to transform each $F_i$
into a Hamilton cycle. At first sight, this looks infeasible for 2 reasons:
\begin{itemize}
\item[(a)] The counting argument in Step 3 only works for reasonably dense subgraphs of $G_{n,p}$. So the $2$-factors produced by the later iterations in Step 3 will have too many cycles.
\item[(b)] We no longer have a graph $G_{thin}$ to use in Step 4.
\end{itemize}
For (a), it turns out that one can extend the counting argument so that it also works for sparser subgraphs of
$G_{n,p}$, provided these subgraphs are pseudorandom (see Lemmas~\ref{SizeOfP} and~\ref{Simple2Fac}).

A very useful observation which goes some way in solving problem (b) is the following:
Let $x_0$ be the vertex of minimum degree in $G_{n,p}$. Then 
whp there is a small but significant gap between $d(x_0)$ and $d(x)$ for any $x \neq x_0$ in $G_{n,p}$
(the gap has size close to $\sqrt{np}$).
Let $G_{left}$ be the subgraph of $G_{n,p}$ consisting of all the edges not in $G_{reg}$. Then $G_{left}$ has
density about $\sqrt{p/n}$, and we could try to use $G_{left}$ to merge $2$-factors
instead of the graph $G_{thin}$ in Step~4. One problem here is that $G_{left}$ is not actually pseudorandom as it 
is just the `leftover' from the Tutte application. 

We use the following idea to overcome this problem. Right at the start, we remove the edges of a random subgraph $G_5$ of density $p_5$ from $G_{n,p}$, where $np_5 \ll \sqrt{np}$. Let $x_0$ be the vertex of minimum degree of the remaining subgraph $G_1$ of $G_{n,p}$.
The choice of $p_5$ implies that $x_0$ will also be the vertex of minimum degree in $G_{n,p}$.
Now add all edges of $G_5$ which are incident to $x_0$ to $G_1$.
Then $\delta(G_1)=\delta(G_{n,p})$ and it turns out that we can still apply Tutte's theorem to obtain an $r$-regular subgraph $G_{reg}$ of $G_1$ (see Lemma~\ref{Tutte}). Again we obtain a decomposition of $G_{reg}$ into $2$-factors $F_i$.

The key advantage of this method is that $G_5 - {x_0}$ will be pseudorandom, and so one can use $G_5$ in the same way as $G_{thin}$
to merge the $F_i$ into Hamilton cycles. However, it turns out that $G_5$ is far too sparse to complete the process:
if e.g.~$p$ is a constant, then our bound on the total number of cycles in the $F_i$ is significantly larger than
$n^{3/2}$. On the other hand, $G_5$ has significantly fewer than $n^{3/2}$ edges. To reduce the number of cycles in $F_i$ by $1$, 
we need about $\log n$ edges from $G_5$. So even if we return the same number of edges from $F_i$ to $G_5$ each time, 
we cannot assume $G_5$ to be pseudorandom after dealing with even a small proportion of the $F_i$.

Our main idea is to use an iterative approach to overcome this. Choose $p_2$, $p_3$ and $p_4$ so that 
$p_5 \ll p_4 \ll p_3 \ll p_2 \ll p$. Then choose randomly edge-disjoint subgraphs $G_i$ of $G_{n, p}$, of density $p_i$, for $i = 2, 3, 4, 5$. (For simplicity, assume that $G_2$, $G_3$ and $G_4$ are regular of even degree in what follows.)
Let $G_1$ be the remaining subgraph of $G_{n,p}$ (with all of the edges of $G_5$ which are incident to $x_0$ added to $G_1$) and consider $\delta(G_1)/2$
edge-disjoint $2$-factors of $F_i$ of $G_1$.
Then as our first iteration we transform all the $F_i$ into Hamilton cycles, where $G_2$ plays the role of $G_{thin}$.
In the second iteration, we do the same for (the leftover of) $G_2$, with $G_3$ playing the role of $G_{thin}$, and so on until we reach $G_5$. (Note that there is no need to do anything with the leftover of $G_5$, as $G_5$ contains no edges incident to $x_0$.) In total, this gives us a set of $\lfloor (\delta(G_1) + \delta(G_2) + \delta(G_3) + \delta(G_4))/2\rfloor = \lfloor d(x_0)/2 \rfloor$ edge-disjoint Hamilton cycles which together contain all of the edges incident to $x_0$.

The new problem that we now face is that e.g. in the second iteration the graph $G_2$ is a `leftover' from the first iteration.
So its pseudorandom properties are not as strong as we need them to be (e.g.~in order to apply the 
strengthened version of the counting argument in Step 3) due to the existence of some `bad edges' which were 
moved into $G_2$ from the $F_i$ during the first iteration.

In order to deal with this, we perform some intermediate steps (about $\tau = \log n/\log \log n$ per iteration) using yet more random graphs
$G_{(2,j)}$ (which we must also remove from $G_{n,p}$ initially). So first $G_{(2, 1)}$ plays the role of $G_{thin}$ when transforming the $2$-factors of $G_1$ into Hamilton cycles. We then transform (the leftover of) $G_{(2, 1)}$ into Hamilton cycles using $G_{(2, 2)}$ and $G_{(2, 3)}$, and then the leftover of $G_{(2, 2)}$ and $G_{(2, 3)}$ using $G_{(2, 4)}$ and $G_{(2, 5)}$, etc. Roughly speaking, after $\tau$ iterations we will have
replaced $G_{(2, 1)}$ by a graph $G_{(2, 2\tau + 1)}$ of almost the same density as $G_{(2, 1)}$, but containing no bad edges (see Lemma~\ref{MiniHamDecomp}). So it is now possible to carry out the second iteration (as described in the previous paragraphs) with $G_{(2, 2\tau + 1)}$ in place of $G_2$ and $G_{(3, 1)}$ in place of $G_3$, and similarly the third and further iterations.

Thus Lemma~\ref{MiniHamDecomp} can be regarded as one of the key statements of the paper, and we believe it is of independent interest: 
as indicated above, it states that given a regular graph $H_0$ (which satisfies some fairly weak quasirandomness conditions)
and a graph $H$ which is the union of many quasirandom graphs $G_{(i,j)}$, we can find a set of edge-disjoint Hamilton cycles in $H_0 \cup H$ covering all edges of $H_0$, even if $H$
is allowed to be much sparser than $H_0$.

This paper is organized as follows:
In Section~\ref{PRgraphs}, we define the pseudorandomness properties we need in the proof and show that with very high probability they are satisfied by $G_{n, p}$. In Section \ref{TutteSection} we use Tutte's theorem find an analogue of $G_{reg}$. In Section \ref{2FactorSection} we use our extended counting argument to split this analogue of $G_{reg}$ into $2$-factors. In Section \ref{merging} we show how to merge the cycles in each of these $2$-factors into a Hamilton cycle. Finally in Section \ref{Completing} we combine the results from the rest of the paper to prove Theorem~\ref{HamDecomp}, using the iterative approach discussed above.

\section{Notation} \label{notation}

Throughout the paper we use the following notation: for a graph $G$ and sets $A, B$ of vertices of $G$, we write $e_G(A, B)$ for the number of edges of $G$ with one endpoint in $A$ and the other in $B$. Let $e_G(A) = e_G(A, A)$. For a graph $G$, let $e(G)$ denote the number of edges of $G$, and for a spanning subgraph $H$ of $G$, let $G \backslash H$ denote the graph obtained by removing the edges of $H$ from $G$. $N_G(A)$ is always taken to be the external neighbourhood of $A$, i.e., $N_G(A) = (\bigcup_{x \in A} N(x)) \backslash A$. $\log$ denotes the natural logarithm, and we write $\log^a n$ for $(\log n)^a$.

Since we are aiming to prove a result with high probability, we may and do assume throughout the paper that $n$ is always sufficiently large for our estimates to hold. Further we omit floor and ceiling symbols, and assume large quantities to be integers, whenever this does not have a significant effect on the argument.

\section{Pseudorandom graphs} \label{PRgraphs}

Our aim in this section is to establish several properties of random graphs which we will use later on, mainly regarding the degree sequence and expansion of small sets. For many of these, the fact that they hold whp in $G_{n, p}$ is well known; however, we need them to hold in $\log n/\log \log n$ random spanning subgraphs of $G_{n, p}$ simultaneously. To accomplish this we first show that they hold in $G_{n, p}$ with probability $1 - O(1/\log n)$, and then take a union bound.

The following useful definition is due to Thomason \cite{T87}. It involves tighter bounds than the similar and more common notion of $\varepsilon$-regularity. We will rely on these in the proof of Theorem~\ref{HamDecomp}.

\begin{defin} \label{jumbled} Let $p, \beta \geq 0$ with $p \leq 1$. A graph $G$ is $(p, \beta)$\emph{-jumbled} if $|e_G(S) - p{s \choose 2}| \leq \beta s$ for every $S \subseteq V(G)$ with $|S| = s$.
\end{defin}

We will often use the following immediate consequence of Definition \ref{jumbled}: Let $G$ be a $(p, \beta)$-jumbled graph and let $S, T \subseteq V(G)$ be disjoint. Then
\begin{equation} \label{DisjointJumbled}
|e_G(S, T) - p|S||T|| \leq 2 \beta (|S| + |T|).
\end{equation}
To see this, note that $e_G(S, T) = e_G(S \cup T) - e_G(S) - e_G(T)$; now applying Definition \ref{jumbled} and using the triangle inequality implies (\ref{DisjointJumbled}).

The following two definitions formalise the notion of `degree gap', which we need in our proof.

\begin{defin} \label{DegreeJump} Let $G$ be a graph on $n$ vertices with a vertex $x_0$ of minimum degree, and let $u \leq n$. Then $G$ is \emph{u-jumping} if every vertex of $G$ apart from $x_0$ has degree at least $\delta(G) + u$.
\end{defin}

\begin{defin} \label{StrongDegreeJump} Let $G$ be a graph on $n$ vertices. For a set $T \subseteq V(G)$, let $\overline{d}_G(T)$ be the average degree of the vertices of $T$ in $G$. Then $G$ is \emph{strongly 2-jumping} if $\overline{d}_G(T) \geq \delta(G) + \min\{|T|-1, \log^2 n\}$ for every $T \subseteq V(G)$.
\end{defin}
Note that if $G$ is strongly $2$-jumping then it is also $2$-jumping. In addition to these three properties we will use several other bounds concerning the degree sequence and edge distribution of a random graph. The following definition collects these properties together.

\begin{defin} \label{PseudoRandom} Call a graph $G$ on $n$ vertices $p$-\emph{pseudorandom} if all of the following hold:
\begin{itemize}
\item [\rm (a)] $G$ is $(p, 2\sqrt{np(1-p)})$-jumbled.
\item [\rm (b)] For any disjoint sets $S, T \subseteq V(G)$ with $|S| = s$ and $|T| = t$,
\begin{itemize}
\item [\rm (i)] if $\left(\frac{1}{s} + \frac{1}{t}\right) \frac{\log n}{p} \geq \frac{7}{2}$, then $e_G(S, T) \leq 2(s+t)\log n$,
\item [\rm (ii)] if $\left(\frac{1}{s} + \frac{1}{t}\right) \frac{\log n}{p} \leq \frac{7}{2}$, then $e_G(S, T) \leq 7stp$,
\item [\rm (iii)] if $\frac{\log n}{sp} \geq \frac{7}{4}$, then $e_G(S) \leq 2s\log n$, and
\item [\rm (iv)] if $\frac{\log n}{sp} \leq \frac{7}{4}$, then $e_G(S) \leq 7s^2p/2$.
\end{itemize}
\item [\rm (c)] $np - 2\sqrt{np \log n} \leq \delta(G) \leq np - 200\sqrt{np(1-p)}$ and $\Delta(G) \leq np + 2\sqrt{np \log n}$.
\item [\rm (d)] $G$ is strongly $2$-jumping.
\end{itemize}
\end{defin}

Definition \ref{PseudoRandom}(a) gives good bounds on the densities of large subgraphs but not on those of very small subgraphs, which is why we need (b).

The remainder of this section will be mainly devoted to showing that the random graphs we consider in the rest of the paper are in fact pseudorandom (see Lemma~\ref{ShowPR}).
For~(a), (b) and~(c) we will prove slightly stronger bounds, as these will be used in~\cite{HKLO}.
We will need the following large deviation bounds on the binomial distribution, proved in \cite{sj-tl-ar} (as Theorem 2.1, Corollary 2.3 and Corollary 2.4 respectively):

\begin{lemma} \label{cor23} Let $X \sim Bin(n, p)$. Then the following properties hold:
\begin{itemize}
\item [\rm (i)] If $h > 0$, then $\mathbb{P}\left[X \leq np - h\right] < e^{-h^2/2np}$.
\item [\rm (ii)] If $\varepsilon \leq 3/2$, then $\mathbb{P}[|X - np| \geq \varepsilon np] < 2e^{-\varepsilon^2 np/3}$.
\item [\rm (iii)] If $c > 1$ and $c' = \log c - 1 + 1/c$, then for any $a \geq cnp$, $\mathbb{P}[X \geq a] < e^{-c'a}$.
\end{itemize}
\end{lemma}

The following lemma implies that the property in Definition \ref{PseudoRandom}(a) holds with very high probability in $G_{n, p}$ for the desired range of $p$.

\begin{lemma} \label{ShowJumbled} Suppose that $np/\log^2 n \rightarrow \infty$ and $n(1-p)/\log^2 n \rightarrow \infty$, and let $G \sim G_{n, p}$.
Then the probability that $G$ is not $(p, \frac{3}{2}\sqrt{np(1-p)})$-jumbled is at most $2/n^2$.
\end{lemma}

\proof Without loss of generality we may assume that $0 < p \leq 1/2$ (by considering the complement of $G$ if $p > 1/2$). For each set $S \subseteq V(G)$ with $|S| = s$, we have that $e(S) \sim Bin({s \choose 2}, p)$. Let 
$$\varepsilon = \frac{3 \sqrt{n(1-p)}}{(s-1)\sqrt{p}} \text{ and } N = {s \choose 2}$$ 
and note that $\varepsilon pN = \frac{3}{2}\sqrt{np(1-p)}s$. Call $S$ \emph{bad} if $|e_G(S) - pN| \geq \varepsilon pN$. We consider the following cases: \medskip

\noindent \textbf{Case 1:} $s \geq 2 \sqrt{\frac{n(1-p)}{p}} + 1$. Then $\varepsilon \leq 3/2$, and hence by Lemma~\ref{cor23}(ii) we have that the probability of $S$ being bad is at most $2 e^{-\varepsilon^2 pN/3}$. But since $p \leq 1/2$, we have $\varepsilon^2 pN/3 = \frac{3 ns(1-p)}{2(s-1)} \geq \frac{3n}{4}$. So the probability that $S$ is bad is at most $e^{-3n/4}$. \medskip

\noindent \textbf{Case 2:} $s \leq 2 \sqrt{\frac{n(1-p)}{p}} + 1$. Since $\varepsilon \geq 3/2$ in this case, we have that $\varepsilon pN \geq pN$ and hence $\mathbb{P}[e_G(S) < pN - \varepsilon pN] = 0$. Let $c = \varepsilon + 1$, so $c \geq 5/2$. Let $c' = \log c - 1 + 1/c$ and note that since $c'$ is an increasing function of $c$ for $c > 1$, we have $c' \geq \log (5/2) - 1 + 2/5 > 1/4$. Now by Lemma~\ref{cor23}(iii) applied with $a = (1 + \varepsilon)pN$, the probability that $S$ is bad is at most%
    \COMMENT{We use that $np/\log^2 n \rightarrow \infty$ and $n(1-p)/\log^2 n \rightarrow \infty$ for the last inequality.}
$e^{-c'(1 + \varepsilon) pN} \leq e^{-c'\sqrt{np(1-p)}s} \leq e^{-4s \log n} = n^{-4s}$. \medskip

We now take a union bound on the probability that there exists a bad set $S$. Firstly for any $s \leq 2 \sqrt{\frac{n(1-p)}{p}} + 1$, the probability that there is some bad set $S$ with $|S| = s$ is at most ${n \choose s}n^{-4s} < n^{-3s} \leq n^{-3}$. Summing over all such $s$, we have an error bound of $n^{1 - 3} = 1/n^2$. Further the probability that there exists a bad set $S$ with $|S| = s \geq 2 \sqrt{\frac{n(1-p)}{p}} + 1$ is at most $2^n e^{-3n/4} \leq 1/n^2$. Now adding these two bounds completes the proof. \endproof

To check the properties in Definition \ref{PseudoRandom}(b) we can use the following two simple lemmas, which are slight
strengthenings of Lemmas 5 and 6 respectively in \cite{AHDoRG}.

\begin{lemma} \label{ABedges}
Let $G \sim G_{n,p}$. Then with probability at least $1 - 1/n^2$, the following properties hold for any disjoint $A, B \subseteq [n]$, with $|A| = a, |B| = b$:
\begin{itemize}
\item[{\rm (i)}] If $\left(\frac{1}{a} + \frac{1}{b}\right) \frac{\log n}{p} \geq \frac{7}{2}$, then $e_{G}(A, B) \leq \frac{3}{2}(a+b)\log n$.
\item[{\rm (ii)}] If $\left(\frac{1}{a} + \frac{1}{b}\right) \frac{\log n}{p} \leq \frac{7}{2}$, then $e_{G}(A, B) \leq 6abp$.
\end{itemize}
\end{lemma}
\proof
Let $X:=e(A,B)$. So $X\sim Bin(|A||B|,p)$. Write $t:=\frac{3}{2}(|A|+|B|)\log n$.
First suppose that $\left(\frac{1}{|A|}+\frac{1}{|B|}\right)\frac{\log n}{p}\ge\frac{7}{2}$.
Then $t\ge\frac{21}{4}|A||B|p$, and so by Lemma~\ref{cor23}(iii) applied with $c=\frac{21}{4}$ we have
\begin{equation} \label{STbound}
\mathbb{P}\left(X\ge t\right)<e^{-\left(\log\frac{21}{4}-\frac{17}{21}\right)t}<e^{-\frac{4t}{5}}=n^{-\frac{6}{5}(|A|+|B|)}.
\end{equation}
If $|A|+|B|\le20$ then $X\le|A||B|\le t$. Otherwise, (\ref{STbound}) implies that
\[
\mathbb{P}(X\ge t)<n^{-4}\frac{1}{n^{|A|}n^{|B|}},\]
and so the result follows by a union bound. 

Now suppose $\left(\frac{1}{|A|}+\frac{1}{|B|}\right)\frac{\log n}{p}\le\frac{7}{2}$,
so that $t\le\frac{21}{4}|A||B|p$. Then, from Lemma~\ref{cor23}(iii) applied with $c=6$ and $a=6|A||B|p$, we obtain
\[\mathbb{P}(X\ge6|A||B|p)\le e^{-\left(\log6-\frac{5}{6}\right)6|A||B|p}<e^{-\frac{21}{4}|A||B|p}\le e^{-t},\]
and so the result holds as before.
\endproof

The proof of the next result is almost the same as that of Lemma~\ref{ABedges}, and so we omit it.%
         \COMMENT{Let $X=e(A)$, so that $X\sim Bin\left(\binom{|A|}{2},p\right)$,
and write $t=\frac{3}{2}(|A|-1)\log n$. First suppose $\frac{\log n}{|A|p}\ge\frac{7}{4}$.
Then $t\ge\frac{21}{4}p\binom{|A|}{2}$, and so by Lemma~\ref{cor23}(iii)
$\mathbb{P}(X\ge t)<e^{-\frac{4t}{5}}=n^{-\frac{6}{5}(|A|-1)}$.
If $|A|\le 21$ then $X\le\binom{|A|}{2}\le t$. Otherwise,
we have $\mathbb{P}(X\ge t)< \frac{n^{-\frac{|A|}{5}+\frac{6}{5}}}{n^{|A|}}=n^{-3}\frac{1}{n^{|A|}}$
and so the result follows by a union bound.
Now suppose $\frac{\log n}{|A|p}\le\frac{7}{4}$, so that $t\le\frac{21}{4}p\binom{|A|}{2}$.
Then we obtain $\mathbb{P}\left(X\ge3|A|^{2}p\right)\le\mathbb{P}\left(X\ge6\binom{|A|}{2}p\right)<e^{-6\left(\log6-\frac{5}{6}\right)\binom{|A|}{2}}<e^{-\frac{21}{4}\binom{|A|}{2}}\le e^{-t}$
and so the result holds as before.}

\begin{lemma} \label{Aedges}
Let $G \sim G(n, p)$. Then with probability at least $1 - 1/n^2$, the following properties hold for every $A \subseteq V(G)$ with $|A| = a$:
\begin{itemize}
\item[{\rm (i)}] If $\frac{\log n}{ap} \geq \frac{7}{4}$, then $e_{G}(A) \leq \frac{3}{2} a \log n$.
\item[{\rm (ii)}] If $\frac{\log n}{ap} \leq \frac{7}{4}$, then $e_{G}(A) \leq 3a^2p$.
\end{itemize}
\end{lemma}

It remains to establish the bounds in Definition \ref{PseudoRandom}(c) on the minimum and maximum degree of $G_{n, p}$ (see Lemma~\ref{MinMaxDeg}), and the fact that $G_{n, p}$ is strongly $2$-jumping with probability $1 - O(1/\log n)$ (see Lemma~\ref{DegJumpProb2}). For this we need estimates on the binomial distribution which do not follow from standard Chernoff bounds (see Lemma~\ref{BinRatios}). These estimates use the following notation, where $X \sim Bin(n-1, p)$:
\begin{itemize}
\item $b(r) = \mathbb{P}\left[X = r \right] = {n-1 \choose r} p^r (1-p)^{n-r-1}$,
\item $B(m_1, m_2) = \mathbb{P}\left[m_1 \leq X \leq m_2 \right]$, and
\item $B(m) = \mathbb{P}\left[X \leq m \right]$.
\end{itemize}
$b'(r)$, $B'(m_1, m_2)$ and $B'(m)$ are defined similarly for $X \sim Bin(n-2, p)$.

\begin{lemma} \label{BinRatios} Suppose that $np(1-p) \rightarrow \infty$. Let $m_1 = np - 2\sqrt{np \log n}$, $m_2 = np - 200\sqrt{np(1-p)}$,
$m_3 = np - \frac{15}{8}\sqrt{np \log n}$ and $\lambda = 1 - 1/(8 \log^3 n)$. Then 
\begin{itemize}
\item[\rm (i)] $nB(m_2) \geq nb(m_2) \geq \sqrt{n}/\log n$, and
\item[\rm (ii)] $nB(m_1) \leq nB(m_3)\le 1/\sqrt{n}$.
\end{itemize}
Suppose in addition that $p \geq 48^2 \log^7 n/n$ and $1-p \geq 36 n^{-1/2} \log^{7/2} n$. Then
\begin{itemize}
\item[\rm (iii)] $\frac{b(r-1)}{b(r)} \geq \lambda$ for each $r \geq m_1 - 8 \log^3 n$, 
\item[\rm (iv)] $\frac{B(m_1 - 8 \log^3 n, r)}{b(r)} \geq 4 \log^3 n$ for each $r \geq m_1$, and
\item[\rm (v)] $\frac{b'(r)}{b(r)} \leq 1 + 1/\log n$ for each $r \geq m_1$.
\end{itemize}
\end{lemma}

\proof (i) Let $X \sim Bin(n-1, p)$, $\sigma = \sqrt{(n-1)p(1-p)}$ and $m'_2 = np - 201\sqrt{np(1-p)}$. Now the de Moivre-Laplace Theorem (see e.g., Theorem 1.6 of \cite{Bollobasbook}) states that if $\sigma \rightarrow \infty$ and $x_2 > x_1$ are constants, then
$$B((n-1)p + x_1\sigma, (n-1)p + x_2\sigma) = (1 + o(1))(\phi(x_2) - \phi(x_1)),$$
where $\phi(x)$ is the cumulative density function of the normal distribution. Clearly 
\begin{align*}
&B(np + x_1\sqrt{np(1-p)}, np + x_2\sqrt{np(1-p)}) \\
&= (1+o(1))B((n-1)p + x_1\sigma, (n-1)p + x_2\sigma).
\end{align*} 
To see this, note that the boundaries of the interval change by at most $2$ and that $b(r) \rightarrow 0$ for any $r$. Hence we have that 
$$B(m'_2, m_2) = (1 + o(1))c,$$ 
where $c = \phi(-200) - \phi(-201)$ is constant. Clearly $b(r) \leq b(m_2)$ for each $m'_2 \leq r \leq m_2$. So 
$$nb(m_2) \geq \frac{(1 + o(1))cn}{\sqrt{np(1-p)}+1} \geq \frac{\sqrt{n}}{\log n}.$$ 

(ii) Let $X \sim Bin(n-1, p)$. By Lemma~\ref{cor23}(i) we have that
\begin{align*}
\mathbb{P}\left[X \leq np - 15\sqrt{np \log n}/8\right] &< \mathbb{P}\left[X \leq (n-1)p - 7\sqrt{np \log n}/4\right] \\ 
&\leq e^{-(7\sqrt{np \log n}/4)^2/2np} \leq e^{-\frac{3}{2}\log n} = \frac{1}{n^{3/2}}.
\end{align*}

(iii) We have
$$\frac{b(r-1)}{b(r)} = \frac{{n-1 \choose r-1}p^{r-1} (1-p)^{n-r}}{{n-1 \choose r}p^{r} (1-p)^{n-r-1}} = \frac{(n-1)!r!(n-r-1)!(1-p)}{(n-1)!(r-1)!(n-r)!p} = \frac{r(1-p)}{(n-r)p},$$
and so
\begin{align*}
1 - \frac{b(r-1)}{b(r)} = \frac{np - r}{(n-r)p} &\leq \frac{2\sqrt{np \log n} + 8 \log^3 n}{(n - np + 2\sqrt{np \log n} + 8 \log^3 n)p} \\
&\leq \frac{3 \sqrt{np \log n}}{np - np^2} = \frac{3 \sqrt{\log n}}{(1-p)\sqrt{np}}.
\end{align*}
Now if $p \geq 1/2$ then $(1-p)\sqrt{np} \geq (36 n^{-1/2} \log^{7/2} n)(\sqrt{n/2})$ and the result follows. On the other hand if $p \leq 1/2$ then $(1-p)\sqrt{np} \geq (48 \log^{7/2} n)/2$ and again the result follows.

(iv) Note that
\begin{align*}
\frac{B(m_1 - 8 \log^3 n, r)}{b(r)} &\geq 1 + \frac{b(r-1)}{b(r)} + \frac{b(r-2)}{b(r)} + \ldots + \frac{b(r - 8\log^3 n)}{b(r)} \\
&\geq 1 + \lambda + \lambda^2 + \ldots + \lambda^{8 \log^3 n} = \frac{1 - \lambda^{8 \log^3 n + 1}}{1-\lambda} \\
&\geq \frac{1 - e^{-1}}{1/(8 \log^3 n)} \geq 4 \log^3 n,
\end{align*}
where in the final line we use that $\lambda \leq e^{-1/(8 \log^3 n)}$.

(v) \begin{align*}
\frac{b'(r)}{b(r)} &= \frac{{n-2 \choose r} p^{r} (1-p)^{n-r-2}}{{n-1 \choose r} p^r (1-p)^{n - r - 1}} = \frac{(n-2)!r! (n-r-1)!}{(n-1)!r!(n-r-2)!(1-p)} \\
&= \frac{n-r-1}{(n-1)(1-p)} = 1 + \frac{np - p - r}{(n-1)(1-p)} \leq 1 + \frac{2\sqrt{np\log n}}{(n-1)(1-p)} \\
&\leq 1 + \frac{3\sqrt{p \log n}}{\sqrt{n}(1-p)} \leq 1 + \frac{1}{\log n},
\end{align*}
as desired.
\endproof

\begin{lemma} \label{MinMaxDeg} Let $G \sim G_{n, p}$. Suppose that $np(1-p) \rightarrow \infty$. Then
\begin{itemize} 
\item [\rm (i)] $\mathbb{P}\left[\delta(G) \leq np - \frac{15}{8}\sqrt{np \log n}\right] \leq 1/\sqrt{n}$.
\end{itemize}
Suppose further that $p \geq 48^2 \log^7 n/n$ and $1-p \geq 36 n^{-1/2} \log^{7/2} n$. Let $m_1 = np - 2\sqrt{np \log n}$, $m_2 = np - 200\sqrt{np(1-p)}$ and let $m_1 \leq m \leq m_2$ be such that $nB(m) \geq \log n$ and $nB(m - 1) \leq \log n$. Then
\begin{itemize}
\item [\rm (ii)] $\mathbb{P}\left[\Delta(G) \geq np + \frac{15}{8}\sqrt{np \log n}\right] \leq 1/\log^2 n$,
\item [\rm (iii)] $\mathbb{P}\left[\delta(G) > m\right] \leq 4/\log n$, and
\item [\rm (iv)] $\mathbb{P}\left[\delta(G) > np - 200 \sqrt{np(1-p)}\right] \leq 4/\log n$.
\end{itemize}
\end{lemma}
Note that the results in Chapter 3 of \cite{Bollobasbook} give sharper bounds which hold with high probability. However as mentioned earlier, we need our error probability to be smaller. Since this does seem to affect the precise results we need to prove the bounds explicitly. Also note that by Lemma~\ref{BinRatios}(i) and (ii), $nB(m_1) \leq \log n \leq nB(m_2)$ (with room to spare), and so there exists $m_1 \leq m \leq m_2$ such that $nB(m) \geq \log n$, but $nB(m-1) \leq \log n$. Thus (iii) is not vacuous.

\proof (i) By taking a union bound we have that 
$$\mathbb{P}\left[\delta(G) \leq np - 15\sqrt{np \log n}/8\right] \leq nB(np - 15\sqrt{np \log n}/8) \leq 1/\sqrt{n},$$
where the last inequality follows from Lemma~\ref{BinRatios}(ii). 

(ii) Set $\varepsilon = \frac{7}{4}\sqrt{\log n/np}$, and note that the lower bound on $p$ implies that $\varepsilon \leq 3/2$.
Let $X \sim Bin(n-1, p)$ be the degree of a given vertex. 
By Lemma~\ref{cor23}(ii),
\begin{align*}
\mathbb{P}\left[X \geq np + 15\sqrt{np \log n}/8\right] &< \mathbb{P}\left[X \geq (n-1)p + \varepsilon (n-1)p\right] \\ 
&\leq 2 e^{-\varepsilon^2 (n-1) p/3} \leq 2 e^{-1.01 \log n} \leq \frac{1}{n \log^2 n},
\end{align*}
whereupon a union bound gives the result.

(iii) Let $Y$ be the number of vertices $v \in V(G)$ such that $m_1 \leq d(v) \leq m$. As $nB(m_1) \leq 1/\sqrt{n}$ by Lemma~\ref{BinRatios}(ii), we have 
$$\mathbb{E}(Y) = nB(m_1, m) \geq nB(m) - nB(m_1) \geq \log n - 1/\sqrt{n} \geq 2\log n/3$$ 
and 
$$\mathbb{E}_2(Y) = n(n-1)(pB'(m_1-1, m-1)^2 + (1-p)B'(m_1, m)^2) \leq n^2 B'(m_1, m)^2.$$ 
Hence
$$\frac{\sqrt{\mathbb{E}_2(Y)}}{\mathbb{E}(Y)} \leq \frac{\sum_{r = m_1}^{m} nb'(r)}{\sum_{r = m_1}^{m} nb(r)} = \frac{\sum_{r = m_1}^{m} nb(r)\frac{b'(r)}{b(r)}}{\sum_{r = m_1}^{m} nb(r)} \leq 1 + \frac{1}{\log n},$$
where the last inequality follows from Lemma~\ref{BinRatios}(v). Hence 
\begin{align*}
Var(Y) &= \mathbb{E}_2(Y) + \mathbb{E}(Y) - \mathbb{E}(Y)^2 \leq (1 + 1/\log n)^2 \mathbb{E}(Y)^2 + \mathbb{E}(Y) - \mathbb{E}(Y)^2 \\
&= (2/ \log n + 1/\log^2 n)\mathbb{E}(Y)^2 + \mathbb{E}(Y).
\end{align*}
So by Chebyshev's inequality, 
$$\mathbb{P}(Y = 0) \leq \frac{Var(Y)}{\mathbb{E}(Y)^2} \leq \frac{2}{\log n} + \frac{1}{\log^2 n} + \frac{1}{\mathbb{E}(Y)} \leq \frac{4}{\log n}.$$

(iv) This follows immediately from (iii) and the remark after the lemma statement. \endproof

\begin{lemma} \label{DegJumpProb2} Let $G \sim G_{n, p}$. Suppose that $p \geq 48^2 \log^7 n/n$ and $1-p \geq 36 n^{-1/2} \log^{7/2} n$. Then with probability at least $1 - 6/ \log n$, $G$ is strongly $2$-jumping. 
\end{lemma}
Similarly to Lemma~\ref{MinMaxDeg}, note that Theorem 3.15 in \cite{Bollobasbook} would imply Lemma~\ref{DegJumpProb2} if we only required the statement to hold with high probability.

\proof Let $m_1 = np - 2\sqrt{np \log n}$, $m_2 = np - 200\sqrt{np(1-p)}$ and $\lambda = 1 - 1/(8 \log^3 n)$. Let $m_1 \leq m \leq m_2$ be such that $nB(m) \geq \log n$, but $nB(m-1) \leq \log n$. By Lemma \ref{MinMaxDeg}(iii) there exists a vertex of degree at most $m$ with probability at least $1 - 4/\log n$. So it suffices to show that with probability at least $1 - 2/\log n$ there are no two vertices each of degree at most $m + 2\log^2 n$ whose degrees differ by at most $1$. 

Let $Z$ be the number of (unordered) pairs $v_1, v_2$ of vertices such that $m_1 \leq \min\{d(v_1), d(v_2)\} \leq m + 2\log^2 n$ and $|d(v_1) - d(v_2)| \leq 1$. We have
\begin{align*}
\mathbb{E}(Z) &\leq {n \choose 2} \sum_{r = m_1}^{m + 2\log^2 n} pb'(r-1)^2 + 2pb'(r-1)b'(r) + (1-p)b'(r)^2 + 2(1-p)b'(r)b'(r+1) \\
&\leq \frac{n^2}{2} \sum_{r = m_1}^{m + 2\log^2 n} 3b'(r)b'(r+1) \leq \frac{3n^2}{2} \left(1+\frac{1}{\log n}\right)^2 \cdot \sum_{r = m_1}^{m + 2\log^2 n} b(r)b(r+1) \\
&\leq 2n^2 b(m + 2\log^2 n + 1) B(m_1, m + 2\log^2 n),
\end{align*}
where the third inequality follows from Lemma~\ref{BinRatios}(v). Note that 
\begin{equation} \label{lambdaineq}
\lambda^{2\log^2 n + 2} \geq 1 - \frac{2\log^2 n + 2}{8 \log^3 n} \geq \frac{6}{7}.
\end{equation}
Now by Lemma~\ref{BinRatios}(iii), 
\begin{align*}
nB(m_1, m + 2\log^2 n + 1) &= \sum_{r = m_1}^{m + 2\log^2 n + 1} nb(r) \leq \sum_{r = m_1}^{m + 2\log^2 n + 1} \frac{nb(r- 2\log^2 n - 2)}{\lambda^{2\log^2 n + 2}} \\
&= \frac{nB(m_1 - 2\log^2 n - 2, m - 1)}{\lambda^{2\log^2 n + 2}} \\
&\leq \frac{nB(m-1)}{\lambda^{2\log^2 n + 2}} \leq \frac{\log n}{\lambda^{2\log^2 n + 2}} \stackrel{(\ref{lambdaineq})}{\leq} \frac{7 \log n}{6}.
\end{align*}
So Lemma~\ref{BinRatios}(iv) implies that
\begin{align*}
nb(m + 2\log^2 n +1) &\leq \frac{nB(m_1 - 8 \log^3 n, m + 2\log^2 n + 1)}{4 \log^3 n} \\
&= \frac{nB(m_1 - 8 \log^3 n, m_1 - 1) + nB(m_1, m + 2\log^2 n + 1)}{4 \log^3 n} \\
&\leq \frac{8 \log^3 n \cdot nb(m_1) + \frac{7}{6}\log n}{4 \log^3 n}
\end{align*}
and now since $b(m_1) \leq B(m_1)$, Lemma~\ref{BinRatios}(ii) implies that
$$nb(m + 2\log^2 n + 1) \leq \frac{2}{\sqrt{n}} + \frac{7}{24 \log^2 n} \leq \frac{1}{3 \log^2 n}.$$
Hence 
$$\mathbb{E}(Z) \leq 2 \cdot \frac{7 \log n}{6} \cdot \frac{1}{3 \log^2 n} \leq \frac{1}{\log n}.$$ 
So by Markov's inequality $\mathbb{P}[Z \geq 1] \leq 1/\log n$. Now by Lemma~\ref{MinMaxDeg}(i), the probability that there are any vertices at all of degree at most $m_1$ is at most $1/\sqrt{n}$. So with probability at least $1 - 1/\log n - 1/\sqrt{n} \geq 1 - 2/\log n$, there are no pairs of vertices of degree at most $m + 2\log^2 n$, whose degrees differ by at most $1$. \endproof

We now combine the above results to prove that our desired pseudorandomness conditions hold whp in $G_{n, p}$.

\begin{lemma} \label{ShowPR} Let $G \sim G_{n, p}$. Suppose that $p \geq 48^2 \log^7 n/n$ and $1-p \geq 36 n^{-1/2} \log^{7/2} n$. Then the probability that $G$ is \emph{not} $p$-pseudorandom is at most $11/\log n$. 
\end{lemma}

\proof By Lemmas \ref{ShowJumbled}, \ref{ABedges}, \ref{Aedges}, \ref{MinMaxDeg}(i), \ref{MinMaxDeg}(ii), \ref{MinMaxDeg}(iv) and \ref{DegJumpProb2} the probability that $G$ is not $p$-pseudorandom is at most
$$\frac{2}{n^2} + \frac{1}{n^2} + \frac{1}{n^2} + \frac{1}{\sqrt{n}} + \frac{1}{\log^2 n} + \frac{4}{\log n} + \frac{6}{\log n}< \frac{11}{\log n}.$$
\endproof

\begin{lemma} \label{SplitGraph} Let $G_0 \sim G_{n, p_0}$. Let $p_1, \ldots, p_5$ be positive reals such that $p_1 + p_2 + p_3 + p_4 + p_5 = p_0$. Let $m_2, m_3, m_4 \leq \log n/\log \log n$ be positive integers. For each $i = 2, 3, 4$, let $p_{(i, 1)}, \ldots, p_{(i, 2m_i+1)}$ be positive reals such that $p_{(i, 1)} + p_{(i, 2)} + \ldots + p_{(i, 2m_i+1)} = p_i$. Suppose that $48^2 \log^7 n/n \leq p \leq 1 - 36 n^{-1/2} \log^{7/2} n$ whenever $p = p_i$ for $0 \leq i \leq 5$, $p=p_{(i, j)}$ for $i = 2, 3, 4$ and $1 \leq j \leq 2m_i + 1$, or $p = p_i + p_{i + 1}$ for $i = 2, 3, 4$.

Define $G_1, G_5$ and $G_{(i, 1)}, \ldots, G_{(i, 2m_i+1)}$ for $i = 2, 3, 4$, as follows. For each edge $ab$ of $G_0$:
\begin{itemize}
\item With probability $p_1/p_0$ let $ab$ be an edge of $G_1$.
\item For $i = 2, 3, 4$ and $1 \leq j \leq 2m_i + 1$, with probability $p_{(i, j)}/p_0$ let $ab$ be an edge of $G_{(i, j)}$.
\item Otherwise (i.e. with probability $p_5/p_0$), let $ab$ be an edge of $G_5$.
\end{itemize}
For $i = 2, 3, 4$, let $G_i = \bigcup_{j=1}^{2m_i + 1} G_{(i, j)}$. Then whp, the following properties hold:
\begin{itemize}
\item [\rm (i)] $G_i$ is $p_i$-pseudorandom for each $0 \leq i \leq 5$,
\item [\rm (ii)] $G_{(i, j)}$ is $p_{(i, j)}$-pseudorandom for $i = 2, 3, 4$ and $1 \leq j \leq 2m_i + 1$, and
\item [\rm (iii)] $G_i \cup G_{i + 1}$ is $(p_i + p_{i+1})$-pseudorandom for $i = 2, 3, 4$.
\end{itemize}
\end{lemma}

\proof First note that $G_{(i, j)} \sim G_{n, p_{(i, j)}}$ for each $i, j$ and that $G_i \sim G_{n, p_i}$ for each $i$. Further, $G_i \cup G_{i+1} \sim G_{n, (p_i + p_{i + 1})}$ for $i = 2, 3, 4$.

Hence by Lemma \ref{ShowPR} for each graph, the probability that it is not $p$-pseudorandom for the relevant value of $p$ is at most $11/\log n$, and so by taking a union bound we have that the probability that at least one of these graphs is not pseudorandom is at most $67/\log \log n$, which tends to $0$ as $n \rightarrow \infty$. \endproof

At one point in the proof of Theorem \ref{HamDecomp} we will need $G_{n, p}$ to be $u$-jumping for some $u > 2$. For this we use the following result.

\begin{lemma} \label{DegJumpProb} Let $G \sim G_{n, p}$ with $np/\log n \rightarrow \infty$ and $n(1-p)/\log n \rightarrow \infty$. Then with high probability, $G$ is $8\sqrt{np(1-p)}/\log^{3/4} n$-jumping.
\end{lemma}

We will apply this lemma only once (in Section \ref{Completing}), so in this case a whp estimate is sufficient.

\proof The conditions of the lemma immediately imply that $np(1-p)/\log n \rightarrow \infty$. Let $m = \min\{(np(1-p)/\log n)^{1/5}, \log^{1/16} n\}$ and let $\alpha = 8\log^{-1/16} n$; note that $m \rightarrow \infty$ and $\alpha \rightarrow 0$. Applying Theorem 3.15 in \cite{Bollobasbook} to the complement of $G$, we have that $G$ is $u$-jumping where
$$u = \frac{\alpha}{m^2}\left(\frac{np(1-p)}{\log n}\right)^{1/2} \geq \frac{8\sqrt{np(1-p)}}{\log^{3/4} n}.$$ \endproof

\section{Constructing regular spanning subgraphs} \label{TutteSection}

The first step in our general strategy for finding a large collection of Hamilton cycles in a $p$-pseudorandom graph $G$ is to construct a regular spanning subgraph of $G$ of degree $\delta(G)$ if $\delta(G)$ is even, or of degree $\delta(G) - 1$ if $\delta(G)$ is odd. The aim of this section is to establish that this is always possible for our range of $p$.

We use Tutte's $r$-factor theorem: Given a graph $G$, an integer $r$, and disjoint subsets $S$ and $T$ of $V(G)$, let
$$R_r(S, T) = \sum_{v \in T} d(v) - e_G(S, T) + r(|S| - |T|).$$
Let $Q_r(S, T)$ be the number of components $C$ of $G -(S \cup T)$, such that $r|C| + e_G(C, T)$ is odd.

\begin{thm} [Tutte \cite{Tutte52}] \label{TutteThm} Let $G$ be a graph and $r$ be a positive integer. Then $G$ contains an $r$-factor if and only if $R_r(S, T) \geq Q_r(S, T)$ for every pair $S, T$ of disjoint subsets of $V(G)$.
\end{thm}

For the particular case of a $1$-factor (that is, a perfect matching), we use the following simpler result:

\begin{thm} [Tutte \cite{Tutte47}] \label{Tutte1Fac} Let $G$ be a graph. Then $G$ has a perfect matching if and only if for every $S \subseteq V(G)$, the number of components of $G - S$ which have an odd number of vertices is at most $|S|$.
\end{thm}

In order to make use of Theorem~\ref{TutteThm} we first need to bound $Q_r(S, T)$. For this we use the following lemma:

\begin{lemma} \label{Components} Let $G$ be an $r_G/n$-pseudorandom graph on $n$ vertices such that $r_G \geq \log n$. Let $r_H \geq \frac{549}{550}r_G$, and let $H$ be a spanning subgraph of $G$ such that $\delta(H) = r_H$. Let $W \subseteq V(G)$, and suppose that $\delta(H[W]) \geq r_H/3$. Then for any nonempty $B \subseteq W$, the number of components of $H[W \backslash B]$ is at most $|B|$. Moreover, $H[W]$ is connected.
\end{lemma}

\proof Suppose that a set $B$ violates the first assertion. \medskip

\noindent \textbf{Claim:} \emph{$H[W \backslash B]$ cannot have two disjoint isolated sets (i.e., unions of components of $H[W\backslash B]$) each of size at least $n/32$.} \medskip

\noindent To prove the claim, it suffices to show that for disjoint sets $S, T \subseteq W \backslash B$ such that $|S|, |T| \geq n/32$, we have $e_H(S, T) > 0$. To see that this holds, note that since $G$ is $(r_G/n, 2\sqrt{r_G})$-jumbled by Definition \ref{PseudoRandom}(a), (\ref{DisjointJumbled}) implies that
$$e_G(S, T) \geq \frac{r_G|S||T|}{n} - 4\sqrt{r_G}(|S| + |T|) \geq \frac{r_Gn}{1024} - 4\sqrt{r_G}n.$$
But since $r_H \geq \frac{549}{550}r_G$, it follows that 
$$|E(G)\backslash E(H)| = e(G) - e(H) \leq \frac{r_Gn}{2} + 2\sqrt{r_G}n - \frac{549r_Gn}{1100} = \frac{r_Gn}{1100} + 2\sqrt{r_G}n,$$
and hence $e_H(S, T) > 0$, which proves the claim. \medskip

Let $\mathcal{C}$ be the set of components of $H[W \backslash B]$ of size less than $n/32$, and let $A = \bigcup \mathcal{C}$. Suppose that $|A| \geq n/8$. Then we can form an isolated set $S$ such that $n/32 \leq |S| \leq n/16$ by taking successive unions of components in $\mathcal{C}$. Similarly we can form another isolated set $T$ with $n/32 \leq |T| \leq n/16$ using the remaining components in $\mathcal{C}$. But now $S$ and $T$ are disjoint, which contradicts the claim. Hence $|A| < n/8$. Moreover the claim implies that $H[W \backslash B]$ has exactly one component which is not in $\mathcal{C}$.

Now since by assumption, $H[W \backslash B]$ has more than $|B|$ components, we have $|A| \geq |\mathcal{C}| \geq |B|$. Now since $A$ is isolated in $H[W \backslash B]$, it follows that the $H[W]$-neighbourhood of $A$ lies entirely in $B$. Hence every edge of $H[W]$ which is incident to some vertex of $A$ lies in $E_H(A) \cup E_H(A, B)$. It follows that $\sum_{v \in A} N_{H[W]}(v) \leq 2e_H(A) + e_H(A, B)$, noting that edges in $e_H(A)$ will be counted twice on the left-hand side. So 
\begin{align*}
r_H|A|/3 &\leq |A|\delta(H[W]) \leq 2e_H(A) + e_H(A, B) \leq 2e_G(A) + e_G(A, B) \\
&\leq r_G|A|^2/n + 4\sqrt{r_G}|A| + r_G|A||B|/n + 4\sqrt{r_G}(|A| + |B|) \\
&\leq 2r_G|A|^2/n + 12\sqrt{r_G}|A| \leq \left(\frac{1}{4} + \frac{12}{\sqrt{r_G}}\right)r_G|A|,
\end{align*}
which is a contradiction unless $|A| = 0$. But if $|A| = 0$ then we have $B = \emptyset$, which proves the assertion. The moreover part follows from the special case where $|B| = 1$. \endproof

\begin{cor} \label{Qrst} Let $G$ be an $r_G/n$-pseudorandom graph with $r_G \geq \log^2 n$ and let $G'$ be a graph obtained from $G$ by first deleting an arbitrary matching $M$ and then adding an arbitrary set of additional edges. Let $r$ be an even integer. Then defining $Q_r(S, T)$ with respect to $G'$, any disjoint subsets $S, T$ of $V(G)$ satisfy $Q_r(S, T) \leq |S| + |T|$.
\end{cor}

\proof By Definition \ref{PseudoRandom}(c) we have that $\delta(G \backslash M) \geq r_G - 2\sqrt{r_G \log n} - 1 \geq \frac{549}{550}r_G$. Note that the number of components of $G' - B$ is always at most that of $(G\backslash M) - B$ for any $B \subseteq V(G)$. Now if $S \cup T \neq \emptyset$ then applying Lemma~\ref{Components} with $G = G$, $H = G \backslash M$, $W = V(G)$ and $B = S \cup T$ implies that $H - (S \cup T)$ has at most $|S \cup T|$ components. It follows that $Q_r(S, T) \leq |S \cup T| = |S| + |T|$. On the other hand, if $S = T = \emptyset$ then Lemma~\ref{Components} implies that $H$ (and thus $G'$) is connected. Let $C$ be the unique component of $G'$ and note that $r|C| + e_{G'}(C, T) = r|C|$, which is even. So in this case $Q_r(S, T) = 0 = |S| + |T|$. \endproof

We are now ready to prove that every pseudorandom graph $G$ has a regular spanning subgraph whose degree is equal to $2\lfloor \delta(G)/2 \rfloor$. In fact we will prove the following slightly stronger statement, which gives the same result even if $G$ is modified slightly from being pseudorandom.

\begin{lemma} \label{Tutte} Let $G$ be a $p$-pseudorandom graph on $n$ vertices such that $np(1-p)/\log^2 n \rightarrow \infty$. 
\begin{itemize}
\item[\rm (i)] For any vertex $x \in V(G)$, $G$ contains an optimal matching which covers $x$. 
\item[\rm (ii)] Let $u \leq 4\sqrt{np(1-p)}$, $u \neq 1$ and suppose in addition that $G$ is $2u$-jumping. Let $G'$ be formed from $G$ by adding $u$ edges at the vertex $x_0$ of minimum degree, and if $\delta(G) + u$ is odd, also removing an arbitrary matching $M$. Let $r$ be the greatest even integer which is at most $\delta(G) + u$. Then $G'$ has an $r$-regular spanning subgraph.
\end{itemize}
\end{lemma}

Note that if the matching we remove in (ii) covers $x_0$, then $r = d_{G'}(x_0)$. The case of Lemma \ref{Tutte}(ii) when $u = 0$ is the only place in the proof of Theorem \ref{HamDecomp} where we use the fact that a pseudorandom graph is strongly $2$-jumping (i.e., Definition \ref{PseudoRandom}(d)). One can probably replace this by a weaker condition. On the other hand, Defintion \ref{PseudoRandom}(a), (b) and (c) are probably not in themselves sufficient to prove Lemma \ref{Tutte}. 

\proof (i) If $n$ is even then Theorem~\ref{Tutte1Fac} and Lemma~\ref{Components} together imply that $G$ has a perfect matching, and the result follows immediately. If $n$ is odd remove any vertex $y \neq x$ from $G$. Now by Lemma~\ref{Components} the number of components of $G - (S \cup \{y\})$ is at most $|S| + 1$ for any $S \subseteq V(G)$. But if the number of components is exactly $|S| + 1$, then at least one component must have an even number of vertices (otherwise we have $n - |S| - 1 \equiv |S| + 1 \mod 2$ which cannot hold). Hence we can apply Theorem~\ref{Tutte1Fac} to $G - y$ and the result follows.

~(ii) For disjoint sets $S, T \subseteq V(G)$, let $s = |S|$ and $\rho = |T|/|S|$. Throughout the remainder of the proof let $R_r(S, T)$ and $Q_r(S, T)$ be defined with respect to the graph $G'$. By Corollary \ref{Qrst}, $Q_r(S, T) \leq |S| + |T| = s(\rho + 1)$ and so by Theorem~\ref{TutteThm} it suffices to show that $R_r(S, T) \geq s(\rho + 1)$ for all $S$ and $T$. If $T$ is nonempty, let $\overline{d}_G(T)$ be the average degree, in $G$, of the vertices of $T$, and define $\overline{d}_G(S)$, $\overline{d}_{G'}(T)$ and $\overline{d}_{G'}(S)$ similarly. Note that
$$R_r(S, T) = \overline{d}_{G'}(T) |T| - e_{G'}(S, T) + rs(1 - \rho).$$

Note also that $d_{G'}(x_0) \geq r$. We claim that $d_{G'}(x) \geq r+2$ for all $x \neq x_0$. Indeed, if $u = 0$ and $\delta(G)$ is even, we have
$$d_{G'}(x) = d_{G}(x) \geq \delta(G) + 2 = r+2$$ 
since $G$ is $2$-jumping by Definition \ref{PseudoRandom}(d). If $u = 0$ and $\delta(G)$ is odd, then similarly we have
$$d_{G'}(x) \geq d_G(x) - 1 \geq \delta(G) + 1 = r+2.$$
If $u \geq 2$ and $\delta(G) +u$ is even, we have
$$d_{G'}(x) = d_{G}(x) \geq \delta(G) + 2u \geq r+2,$$
and finally if  $u \geq 2$ and $\delta(G) +u$ is odd, we have
$$d_{G'}(x) \geq d_{G}(x) - 1 \geq \delta(G) + 2u - 1 \geq r+2.$$
This proves the claim. It follows immediately that unless $T = \{x_0\}$, we have $\overline{d}_{G'}(T) \geq r+1$. Now we consider the following cases: \medskip

\noindent \textbf{Case 1:} $|S| = 0$. Observe that $R_r(S, T) = \overline{d}_{G'}(T)|T| - r|T|$ and $Q_r(S, T) \leq |T|$ by Corollary \ref{Qrst}. If either $|T| = 0$ or $\overline{d}_{G'}(T) \geq r+1$ then trivially $R_r(S, T) \geq |T| \geq Q_r(S, T)$. So we are left with the case in which $T = \{x_0\}$ and $d_{G'}(x_0) = r$.

In this case it suffices to prove that $Q_r(S, T) = 0$. To see that this holds, note that applying Lemma~\ref{Components} with $H = G \backslash M$, $W = V(G)$ and $B = \{x_0\}$ implies that $G' - \{x_0\}$ is connected. But the unique component $C$ of $G' - \{x_0\}$ satisfies $r|C| + e_{G'}(C, \{x_0\}) = r(n-1) + r = rn$, and since $r$ is even $C$ cannot be odd. \medskip

\noindent \textbf{Case 2:} $|S| > 0$ and $\rho \leq 1/2$. Since $e_{G'}(S, T) \leq \overline{d}_{G'}(T) |T|$ and $|S| \geq 2|T|$, we have 
$$R_r(S, T) \geq r (|S| - |T|) \geq \frac{r}{3}(|S| + |T|) \geq |S| + |T| \geq Q_r(S, T),$$ 
where the last inequality follows by Corollary \ref{Qrst}. \medskip

\noindent \textbf{Case 3:} $|S| > 0$, $\rho \geq 1/2$ and $\left(\frac{1}{s} + \frac{1}{\rho s}\right)\frac{\log n}{p} \geq \frac{7}{2}$. In this case we have that $e_{G'}(S, T) \leq e_G(S, T) + (\rho + 1)s \leq 2 (\rho+1) s \log n + (\rho + 1)s$, where the last inequality follows from Definition \ref{PseudoRandom}(b). So
\begin{equation}
R_r(S, T) - Q_r(S, T) \geq \rho s(\overline{d}_{G'}(T) - 2\log n - r - 2) + s(r - 2 \log n - 2),
\label{TCase3}
\end{equation}
and it suffices to prove that the right-hand side of (\ref{TCase3}) is non-negative. Now observe that if $\rho \leq \frac{r - 2 \log n - 2}{2 \log n + 2}$, then this is immediate since $\overline{d}_{G'}(T) \geq r$. On the other hand, if $\rho \geq \frac{r - 2 \log n - 2}{2 \log n + 2}$, then $|T| \geq \frac{r - 2 \log n - 2}{2 \log n + 2} \geq \frac{r}{3 \log n} \geq 3 \log n$. So by Definition \ref{PseudoRandom}(d) we have $\overline{d}_{G'}(T) \geq r + 3 \log n - 1 \geq r + 2 \log n + 2$, and the result follows. \medskip

\noindent \textbf{Case 4:} $|S| > 0$, $\rho \geq 1/2$, $\left(\frac{1}{s} + \frac{1}{\rho s}\right)\frac{\log n}{p} \leq \frac{7}{2}$ and $\rho s \leq n/30$. In this case Definition \ref{PseudoRandom}(b) implies that $e_{G'}(S, T) \leq 7 \rho s^2 p + u$. So 
$$R_r(S, T) - Q_r(S, T) \geq \rho s(\overline{d}_{G'}(T) - r - 1) + s(r - 7\rho sp - u/s - 1),$$
and it suffices to prove that the right-hand side of this inequality is non-negative. Recall that $\overline{d}_{G'}(T) \geq r + 1$ as $|T| = \rho s \geq 2(\rho + 1) \log n/7p \geq 2$, and so the first bracket is non-negative. Also $7 \rho sp \leq np/4 < r/3$ by Definition \ref{PseudoRandom}(c), and $u/s < u < r/3$. Hence the second bracket is positive. \medskip

\noindent \textbf{Case 5:} $|S| > 0$, $\rho \geq \frac{1}{2}$ and $\rho s \geq n/30$. Since $G$ is $(p, 2\sqrt{np(1-p)})$-jumbled, (\ref{DisjointJumbled}) implies that
$$e_G(S, T) \leq \rho s^2 p + 4 \sqrt{np(1-p)} s(\rho + 1).$$
Note that $\overline{d}_{G'}(T)|T| - e_{G'}(S, T) \geq \overline{d}_{G}(T)|T| - e_{G}(S, T) - u$. Now we claim that $\overline{d}_{G}(T)|T| - u \geq (np - 130\sqrt{np(1-p)})|T|$. Indeed,
\begin{align*}
\overline{d}_{G}(T)|T| &= 2e_G(T) + e_G(T, V(G) \backslash T) \\
&\stackrel{(\ref{DisjointJumbled})}{\geq} 2p{|T| \choose 2} + p|T|(n-|T|) - 4\sqrt{np(1-p)}|T| - 4\sqrt{np(1-p)}n \\
&\geq np|T| - p|T| - (4n + 4|T|)\sqrt{np(1-p)} \geq (np - 130\sqrt{np(1-p)})|T| + u.
\end{align*}
Hence
\begin{align*}
R_r(S, T)/s &\geq \rho (np - 130\sqrt{np(1-p)}) - \rho sp - 4\sqrt{np(1-p)}(\rho + 1) + r(1-\rho) \\
&\geq \rho(np - r - 142\sqrt{np(1-p)}) + r - \rho sp.
\end{align*}
But 
\begin{align*}
\delta(G) - 1 \leq r \leq \delta(G) + u &\leq np - 200\sqrt{np(1-p)} + 4\sqrt{np(1-p)} \\
&= np - 196\sqrt{np(1-p)}.
\end{align*}
Hence
$$R_r(S, T)/s \geq 54\rho \sqrt{np(1-p)} + \delta(G) - 1 - \rho sp \geq 50\rho \sqrt{np(1-p)} + \delta(G) - \rho sp.$$
Further as $(\rho + 1)s \leq n$, we have 
$$\delta(G) - \rho sp \geq np - 2\sqrt{np \log n} - np\left(\frac{\rho}{\rho + 1}\right) = \frac{np}{\rho + 1} - 2\sqrt{np \log n}.$$
Now if $\rho \leq \frac{1}{6}\sqrt{\frac{np}{\log n}}$ then $\frac{np}{\rho + 1} \geq \frac{np}{3\rho} \geq 2 \sqrt{np \log n}$ and so $R_r(S, T)/s \geq 50\rho\sqrt{np(1-p)} \geq \rho + 1$.

On the other hand if $\rho \geq \frac{1}{6}\sqrt{\frac{np}{\log n}} \geq \sqrt{\frac{\log n}{1 - p}}$ then $2\rho\sqrt{np(1-p)} \geq 2\sqrt{np \log n}$ and so $R_r(S, T)/s \geq 48\rho\sqrt{np(1-p)} + np/(\rho + 1) \geq \rho + 1$. \endproof

We remark that as long as $np/\log^6 n \rightarrow \infty$ and $n(1-p)/\log^6 n \rightarrow \infty$ (i.e., under less restrictive assumptions than those of Theorem \ref{HamDecomp}), Lemma \ref{Tutte} implies that whp $G_{n, p}$ has a $2\lfloor \delta(G_{n, p})/2 \rfloor$-regular spanning subgraph. The fact that $G_{n, p}$ is whp $p$-pseudorandom for this range of $p$ follows by Lemmas \ref{ShowJumbled}--\ref{Aedges}, Lemma \ref{DegJumpProb} and Corollary 3.13 of \cite{Bollobasbook}.

\section{Splitting into 2-factors} \label{2FactorSection}

Our aim in this section is show that under certain conditions, an even-regular graph $H$ can be decomposed into $2$-factors so that the sum of the number of cycles in these $2$-factors is not too large. We would like such a result to hold for arbitrary even-regular graphs $H$, but this does not seem feasible for the densities that we consider. So we first show that it suffices for $H$ to be a spanning subgraph of a pseudorandom graph $G$ (see Corollary \ref{2FacSplit}), and then that it suffices for $H$ to be `partly' pseudorandom (see Corollary \ref{Absorb2Fac}). More precisely, we show that for any even-regular graph $H'$ of degree $r_{H'}$, the union of $H'$ with an even-regular graph $H$ which is close to being $p$-pseudorandom may be decomposed into $2$-factors with few cycles in total, provided only that $p$ is somewhat larger than $r_{H'}/n$.

In order to bound the number of cycles in each $2$-factor we will transfer the problem to a bipartite setting. This will allow us to use the following result (which was conjectured by van der Waerden), which gives a good bound on the total number of perfect matchings in a regular bipartite graph.

\begin{thm} [Egorychev \cite{VDW1}, Falikman \cite{VDW2}] \label{EFWThm}
Let $r \leq n$ be positive integers and let $B$ be an $r$-regular bipartite graph with vertex classes of size $n$. Then the number of perfect matchings of $B$ is at least $(\frac{r}{n})^n n!$.
\end{thm}

Given a bipartite graph $B$ on vertex classes $V = \{v_1, \ldots, v_n\}$ and $V' = \{v'_1, \ldots, v'_n\}$, and a matching $M$ in $B$, let $D(M)$ be the digraph on $[n]$ formed by including an edge from $a$ to $b$ if and only if there is an edge of $M$ between $v_a$ and $v'_b$. Note that if $P$ is a perfect matching then $D(P)$ is $1$-regular. For an integer $C$, let $\mathcal{P}_C(B)$ be the set of perfect matchings $P$ of $B$ such that $D(P)$ has at least $C$ cycles. Let $\mathcal{P}_{k, \ell}(B)$ be the set of perfect matchings $P$ such that $D(P)$ has at least $k$ cycles of length $\ell$. We now use a counting argument to show that $\mathcal{P}_C(B)$ is small whenever $C$ is large and $B$ satisfies a very weak pseudorandomness condition. The proof builds on ideas from \cite{FK05} and later developments in \cite{KO06} and \cite{AHDoRG}.

\begin{lemma} \label{SizeOfP} Let $B$ be an $r_B$-regular bipartite graph on vertex classes $V = \{v_1, \ldots, v_n\}$ and $V' = \{v'_1, \ldots, v'_n\}$. Let $\log^2 n \leq r_G \leq n$ and $\sqrt{r_G}\log n \leq r_B \leq r_G/2$, and let $C = 2n \sqrt{r_G \log n}/r_B$. Suppose that all $S \subseteq V$ and $T \subseteq V'$ such that $|S|, |T| \geq n/\sqrt{r_G}$ satisfy $e_B(S, T) \leq \frac{5r_G(|S| + |T|)^2}{2n}$. Let $X = (\frac{r_B}{n})^n n!$. Then $|\mathcal{P}_C(B)| < X$.
\end{lemma}
A similar result was proved in \cite{FK05} and later in \cite{AHDoRG}.\COMMENT{Didn't we use a different counting argument?} In the proof of Lemma~\ref{SizeOfP} we use the following result, which is part of the proof of Lemma 2 in \cite{FK05} (see inequality (5) there).

\begin{lemma} [Frieze and Krivelevich, \cite{FK05}] \label{Ykl}
Let $B$ be an $r$-regular bipartite graph on vertex classes $V$ and $V'$ of size $n$. Let $m < n$ and let $S \subseteq V$ and $T \subseteq V'$ each have size $n - m$. Let $B'$ be the bipartite subgraph of $B$ induced by $S$ and $T$, and let $Y_m$ be the number of perfect matchings of $B'$. Suppose that $e_B(S, T) \geq rn/2$. Then
$$Y_m \leq \left(\frac{(n-2m)r + e_B(V\backslash S, V'\backslash T)}{n - m}\right)^{n-m} e^{-n+m}(3n)^{10n/r}.$$
\end{lemma}
Note that Lemma~\ref{Ykl} is stated in \cite{FK05} with an estimate $m^2$ for $e_B(V\backslash S, V'\backslash T)$, rather than explicitly including the term $e_B(V\backslash S, V'\backslash T)$. Since the graphs we consider are subgraphs of random graphs we will have $e_B(V\backslash S, V'\backslash T) \ll m^2$. This will allow us to obtain nontrivial bounds on $|\mathcal{P}_C(B)|$ even for very sparse graphs.

\proof[Proof of Lemma~\ref{SizeOfP}] Let $k = n \log n/r_B$. For each $\ell$, let $X_{k, \ell} = |\mathcal{P}_{k, \ell}(B)|$, and let $Y_{k, \ell}$ be the maximum number of perfect matchings between two subsets of the vertex classes of $B$, each of size $n - k \ell$. Let $\ell_0 = r_B/\sqrt{r_G \log n}$. Since $\ell_0 \geq \sqrt{\log n}$ we may take $\ell_0$ to be an integer. Note that 
\begin{equation} \label{kl0Bound}
k \ell_0 = \frac{n\sqrt{\log n}}{\sqrt{r_G}} = o(n).
\end{equation}

We first derive a bound on $Y_{k, \ell}$, for $\ell \leq \ell_0$. Consider any $S \subseteq V$ and $T \subseteq V'$, each of size $n - k\ell$. We claim that $e_B(V \backslash S, V' \backslash T) \leq 10r_G k^2 \ell_0^2/n$. Indeed, let $S' \subseteq V$ and $T' \subseteq V'$ be sets of size $k \ell_0$ such that $V \backslash S \subseteq S'$ and $V' \backslash T \subseteq T'$. Then 
$$e_B(V \backslash S, V' \backslash T) \leq e_B(S', T') \leq \frac{5r_G(|S'| + |T'|)^2}{2n} = \frac{10r_G k^2 \ell_0^2}{n}.$$ Moreover, $|S| = |T| = n- o(n)$ by (\ref{kl0Bound}) and so $e_B(S, T) \geq r_B n/2$.\COMMENT{$e_B(S, T) = e_B(V, T) - e_B(V \backslash S, T) \geq |T|r_B - |V \backslash S|r_B \geq nr_B/2$.} Hence by Lemma~\ref{Ykl},
$$Y_{k, \ell} \leq \left(\frac{(n-2k\ell)r_B + 10r_G k^2 \ell_0^2/n}{n - k\ell}\right)^{n-k\ell} e^{-n+k\ell}(3n)^{10n/r_B}.$$
We can further estimate this as follows: Observe that
\begin{align*}
\frac{(n-2k\ell)r_B + 10r_G k^2 \ell_0^2/n}{n-k\ell} &= r_B\left(1 - \frac{k\ell}{n-k\ell} + \frac{10r_G k^2 \ell_0^2}{nr_B(n-k\ell)}\right) \\
&\leq r_B \exp{\left[-\frac{k\ell}{n-k\ell} + \frac{10r_G k^2 \ell_0^2 }{nr_B(n-k\ell)}\right]},
\end{align*}
and hence 
\begin{align} \label{Ykl2}
Y_{k, \ell} &\leq r_B^{n-k\ell} e^{-k \ell} \exp{\left[\frac{10r_G k^2 \ell_0^2 }{nr_B}\right]} e^{-n+k \ell}(3n)^{10n/r_B} \nonumber \\
&= r_B^{n-k\ell} e^{10k} e^{-n}(3n)^{10n/r_B},
\end{align}
where the last line follows from the fact that
$$\frac{r_G k\ell^2_0}{nr_B} \stackrel{(\ref{kl0Bound})}{=} \frac{r_G}{nr_B} \cdot \frac{n\sqrt{\log n}}{\sqrt{r_G}} \cdot \frac{r_B}{\sqrt{r_G \log n}} = 1.$$

Now we proceed to bound $X_{k, \ell}$. Note that if $M_1$ and $M_2$ are matchings such that $D(M_1)$ and $D(M_2)$ are cycles, then $M_1$ and $M_2$ are vertex-disjoint if and only if $D(M_1)$ and $D(M_2)$ are vertex-disjoint. We claim that there are at most ${n \choose k}r_B^{k(\ell-1)}\ell^{-k}$ ways of choosing $k$ vertex-disjoint matchings $M$ of size $\ell$ in $B$, such that $D(M)$ is a cycle of length $\ell$ for each $M$. Indeed, we can construct these matchings as follows: First choose an initial vertex in $V$ for each matching (${n \choose k}$ choices). For an initial vertex $v_{i_1}$ of $M$, select an unused neighbour $v'_{i_2} \in V'$ of $v_{i_1}$ (at most $r_B$ choices) and add $v_{i_1} v'_{i_2}$ to $M$. Then select an unused neighbour $v'_{i_3}$ of $v_{i_2}$ (recall that $v_{i_2}$ is the vertex of $V$ corresponding to $v'_{i_2}$), etc., until we have selected $v'_{i_{\ell}}$. Now if $v'_{i_1}$ is a neighbour of $v_{i_{\ell}}$ then add $v_{i_\ell} v'_{i_1}$ to $M$ and note that $D(M)$ is a cycle of length $\ell$. Repeat this process for each initial vertex, in each case using only vertices we have not used before. This process described above constructs every collection of $k$ vertex-disjoint matchings of size $\ell$ at least $\ell^k$ times, since we choose the initial vertex of each matching arbitrarily. The claim follows immediately.

Since $Y_{k, \ell}$ bounds the number of perfect matchings on the remaining vertices (i.e., those which are not contained in the above matchings), we have
$$X_{k, \ell} \leq {n \choose k} r_B^{k(\ell-1)}\ell^{-k}Y_{k, \ell}.$$
Estimating ${n \choose k} \leq \left(\frac{ne}{k}\right)^k$ and $n! \geq \left(\frac{n}{e}\right)^n$, we have
\begin{align*}
X_{k, \ell} X^{-1} &\stackrel{(\ref{Ykl2})}{\leq} \left(\frac{ne}{k}\right)^k r_B^{k(\ell-1)}\ell^{-k}r_B^{n-k\ell}e^{10k}e^{-n}(3n)^{10n/r_B} \left(\frac{n}{r_B}\right)^n \left(\frac{e}{n}\right)^n \\
&= \left(\frac{ne^{11}}{k}\right)^k r_B^{k\ell-k}\ell^{-k}r_B^{n-k\ell} (3n)^{10n/r_B} r_B^{-n} = \left(\frac{ne^{11}}{k\ell r_B}\right)^k(3n)^{10n/r_B} \\
&= \left(\frac{e^{11}}{\ell \log n}\right)^k (3n)^{10n/r_B}.
\end{align*}
Hence
$$X^{-1} \sum_{\ell=3}^{\ell_0} X_{k, \ell} \leq \left(\left(\frac{e^{11}}{\log n}\right)^{\log n/10} 3n\right)^{10n/r_B} \sum_{\ell=3}^{\ell_0} \ell^{-k}.$$
But $((e^{11}/\log n)^{\log n/10} 3n)^{10n/r_B} < 2^{-10n/r_B} < 1$ and $\ell^{-k} \leq 3^{- \log n} < 1/n$ for $\ell \geq 3$. Hence
$$X^{-1} \sum_{\ell=3}^{\ell_0} X_{k, \ell} < \sum_{\ell=3}^{\ell_0} 1/n < \ell_0/n < 1.$$

Note that since $k \ell_0 \leq C/2$ and $n/\ell_0 = C/2$, we have $C \geq k \ell_0 + n/\ell_0$. But any $1$-regular digraph which, for all $\ell \leq \ell_0$, contains fewer than $k$ cycles of length $\ell$ has fewer than $k \ell_0 + n/\ell_0$ cycles in total. Hence $|\mathcal{P}_C(B)| \leq \sum_{\ell=3}^{\ell_0} X_{k, \ell} < X$. \endproof

The following lemma and corollary accomplish our aim of decomposing an even-regular graph $H$ which is close to being pseudorandom into $2$-factors with few cycles in total.

\begin{lemma} \label{Simple2Fac} Let $r_G \geq \log^2 n$, and let $G$ be an $r_G/n$-pseudorandom graph on $n$ vertices. Let $r_H \geq 2\sqrt{r_G}\log n$ be even and let $H$ be an $r_H$-regular spanning subgraph of $G$. Then $H$ contains a $2$-factor with at most $C = 4n \sqrt{r_G \log n}/r_H$ cycles.
\end{lemma}

\proof Recall that Petersen's theorem \cite{P1891} states that every even-regular graph can be decomposed into $2$-factors. Therefore $H$ can be decomposed into a collection $F_1, \ldots, F_{r_B}$ of $2$-factors, where $r_B = r_H/2$. For each of these $2$-factors we orient the edges so that each cycle is an oriented cycle, thus forming a collection of $1$-regular digraphs on $V(H)$. Now taking the union of these digraphs yields an orientation $D$ of $H$ which is $r_B$-regular.

Label the vertices of $V(H)$ as $v_1, \ldots, v_n$. Form a bipartite graph $B$ on vertex classes $V = \{v_1, \ldots, v_n\}$ and $V' = \{v'_1, \ldots, v'_n\}$ by joining $v_i \in V$ to $v'_j \in V'$ if and only if there is an edge of $D$ from $v_i$ to $v_j$. Now for any sets $S \subseteq V$ and $T \subseteq V'$ with $|S|, |T| \geq n/\sqrt{r_G}$, let $S'$, $T'$ be the corresponding subsets of $V(H)$. Then
\begin{align*}
e_B(S, T) \leq e_H(S' \cup T') \leq e_G(S' \cup T') &\leq \frac{r_G(|S| + |T|)^2}{2n} + 2\sqrt{r_G}(|S| + |T|) \\
&\leq \frac{3r_G(|S| + |T|)^2}{2n},
\end{align*}
where the third inequality follows since $G$ is $(r_G/n, 2\sqrt{r_G})$-jumbled by Definition \ref{PseudoRandom}(a), and the fourth since $2 \leq \sqrt{r_G}(|S| + |T|)/n$. Now Lemma~\ref{SizeOfP} implies that $|\mathcal{P}_C(B)| <  \left(\frac{r_B}{n}\right)^n n!$. But by Theorem~\ref{EFWThm} the total number of perfect matchings of $B$ is at least $\left(\frac{r_B}{n}\right)^n n!$. So there exists a perfect matching $P$ of $B$ such that $D(P)$ has at most $C$ cycles. Now ignoring the orientation of $D(P)$ yields the desired $2$-factor of $H$. \endproof

For a $2$-regular graph $F$, let $c(F)$ be the number of cycles of $F$, and for a collection $\mathcal{F}$ of $2$-factors let $c(\mathcal{F}) = \sum_{F \in \mathcal{F}} c(F)$. Repeated application of Lemma \ref{Simple2Fac} gives us the following result.

\begin{cor} \label{2FacSplit} Let $G$ be an $r_G/n$-pseudorandom graph on $n$ vertices with $r_G \geq \log^2 n$. Let $H$ be an $r_H$-regular spanning subgraph of $G$, such that $r_H$ is even. Then $H$ can be decomposed into a collection of $2$-factors $\mathcal{F} = \{F_1, \ldots, F_{r_H/2}\}$, such that $c(\mathcal{F}) \leq 3n\sqrt{r_G \log^3 n}$.
\end{cor}

\proof If $r_H \leq 8\sqrt{r_G}\log n$ then the result follows from Petersen's theorem \cite{P1891}, noting that trivially $c(F_i) \leq n/3$ for each $i$. So suppose $r_H \geq 8\sqrt{r_G}\log n$. Let $i_0 = r_H/2 - \sqrt{r_G} \log n$ and note that since $i_0 \geq r_H/4 \geq \log^2 n$, we may assume that $i_0$ is an integer. By repeatedly applying Lemma~\ref{Simple2Fac} we have
\begin{align*}
c(\mathcal{F}) &\leq \sum_{i=1}^{i_0} \frac{2n\sqrt{r_G \log n}}{r_H/2 - i + 1} + \sum_{i = i_0 + 1}^{r_H/2} n/3 \leq 2n \sqrt{r_G \log n} \sum_{i=\sqrt{r_G} \log n}^{r_H/2} \frac{1}{i} + n \sqrt{r_G} \log n \\
&\leq 2n\sqrt{r_G \log n}(\log (r_H/2) + \sqrt{\log n}) \leq 3n\sqrt{r_G \log^3 n}.
\end{align*} \endproof

Now we consider an arbitrary even-regular graph $H'$ and aim to show that if we take an even-regular graph $H$ on $V(H')$ which is close to being pseudorandom and edge-disjoint from $H'$, then we can decompose $H \cup H'$ into $2$-factors with a similar bound on the total number of cycles as in Corollary \ref{2FacSplit}. Our strategy will be to first split $H'$ up into matchings. We then extend each of these matchings into two perfect matchings using edges of $H$ (see Lemma~\ref{HalfMatch}), and apply Lemma~\ref{SizeOfP} to transform each of these perfect matchings into a $2$-factor with few cycles (see Lemma~\ref{PMatchComplete}).

We start by considering the case of a single matching. If $n$ is odd, then call a graph $F$ on $n$ vertices a \emph{pseudomatching} if it has a unique vertex of degree $2$ and all other vertices are of degree $1$. A perfect pseudomatching in a graph $G$ is a pseudomatching covering every vertex of $G$.

\begin{lemma} \label{HalfMatch} Let $G$ be an $r_G/n$-pseudorandom graph on $n$ vertices with $r_G \geq \log^2 n$ and let $H$ be a spanning subgraph of $G$ such that $\delta(H) \geq \frac{549}{550}r_G + 5$. Let $M$ be a matching with $V(M) \subseteq V(G)$ and suppose $M$ and $H$ are edge-disjoint. Then 
\begin{itemize}
\item [\rm (i)] if $n$ is even, then there exists a partition of $M$ into submatchings $M_1$, $M_2$, such that each of $M_1, M_2$ can be completed to a perfect matching on $V(G)$ using disjoint sets of edges of $H$, and
\item [\rm (ii)] if $n$ is odd, then there exists a partition of $M$ into submatchings $M_1$, $M_2$, such that each of $M_1, M_2$ can be completed to a perfect pseudomatching on $V(G)$ using disjoint sets of edges of $H$.
\end{itemize}
\end{lemma}

It would be more convenient to extend $M$ directly into a perfect matching rather than splitting it first, but this does not seem to be possible.

\proof Let $r_H = \frac{549}{550}r_G$. We partition $M$ randomly by placing each edge into either $M_1$ or $M_2$, with equal probability. For a vertex $v \in V(G)$, let $Y_v$ be the set of edges $ww'$ of $M$ such that both $w$ and $w'$ lie in $N_H(v)$, and $Z_v$ the set such that exactly one of $w$ and $w'$ lies in $N_H(v)$. Note that $2|Y_v| + |Z_v| \leq d_H(v) \leq \Delta(G) \leq 2r_H$, and that $$|N_H(v) \backslash V(M_1)| \sim 2Bin(|Y_v|, 1/2) + Bin(|Z_v|, 1/2) + |N_H(v) \backslash V(M)|.$$
(Here we use that $M$ and $H$ are edge-disjoint and so the unique neighbour of $v$ in $M$ is not a neighbour of $v$ in $H$.)

Suppose $|Y_v| \geq 8 \log n$. Let $\varepsilon = 4 \sqrt{\log n/|Y_v|} \leq 3/2$ and note that $2\sqrt{r_H \log n} \geq \varepsilon |Y_v|/2$, since $r_H \geq |Y_v|$. By Lemma~\ref{cor23}(ii) we have that
$$\mathbb{P}\left[Bin(|Y_v|, 1/2) \leq |Y_v|/2 - 2\sqrt{r_H \log n}\right] < 2e^{-\varepsilon^2 |Y_v|/6} \leq 2e^{-16 \log n/6} < \frac{1}{n^2}.$$
On the other hand, if $|Y_v| \leq 8 \log n$ then $|Y_v|/2 \leq 2\sqrt{r_H \log n}$ and so $$\mathbb{P}\left[Bin(|Y_v|, 1/2) \leq |Y_v|/2 - 2\sqrt{r_H \log n}\right] = 0.$$
By the same argument we can show that\COMMENT{Note that we only get $|Z_v| \leq 2r_H$, not $|Z_v| \leq r_H$. This means we end up with a worse bound.}
$$\mathbb{P}\left[Bin(|Z_v|, 1/2) \leq |Z_v|/2 - 4\sqrt{r_H \log n}\right] < \frac{1}{n^2}.$$
So whp, every $v \in V(G)$ satisfies
\begin{align*}
|N_H(v) \backslash V(M_1)| &\geq |Y_v| + |Z_v|/2 - 8 \sqrt{r_H \log n} + |N_H(v) \backslash V(M)| \\
&= \frac{1}{2}|N_H(v) \cap V(M)| + |N_H(v) \backslash V(M)| - 8 \sqrt{r_H \log n} \geq r_H/3 + 7,
\end{align*}
and hence $\delta(H - V(M_1)) \geq r_H/3 + 7$. By a similar argument the same holds for $M_2$. Now choose a partition of $M$ into $M_1$ and $M_2$ such that the above bounds on $|N_H(v) \backslash V(M_1)|$ and $|N_H(v) \backslash V(M_2)|$ hold for all vertices $v \in V(G)$.

If $n$ is odd, then we first find a path $v_1 v_2 v_3$ of length $2$ in $H - V(M_1)$ and add it to $M_1$ to form a pseudomatching, so that $V(H) \backslash V(M_1)$ has an even number of vertices. We also remove $v_1 v_2$ and $v_2 v_3$ from $H$. (If $n$ is even then we simply omit this step.) Henceforth we proceed identically in cases ~(i) and ~(ii).

Let $W = V(H) \backslash V(M_1)$ and note that we still have $\delta(H) \geq r_H + 3$ and $\delta(H[W]) \geq r_H/3 + 4$. We now use Theorem~\ref{Tutte1Fac} to show that $H[W]$ contains a perfect matching. For this it suffices to prove that for any set $S \subseteq W$, the number of components of $H[W \backslash S]$ which have an odd number of vertices is at most $|S|$. If $S$ is nonempty then this follows from Lemma~\ref{Components}. On the other hand if $S = \emptyset$ then by Lemma~\ref{Components} $H[W \backslash S] = H[W]$ has exactly one component, namely $W$. But this component has an even number of vertices, and so by Theorem~\ref{Tutte1Fac} $H[W]$ contains a perfect matching $M'_1$. So $M_1 \cup M'_1$ is a perfect matching or pseudomatching on $V(G)$. We delete the edges in $M'_1$ from $H$.

We now repeat the process for $M_2$. First we extend $M_2$ to a pseudomatching (if $n$ is odd). Now we still have $\delta(H) \geq r_H$ and $\delta(H - V(M_2)) \geq r_H/3$, so the conditions of Lemma~\ref{Components} still hold and hence we can complete $M_2$ to a perfect matching or pseudomatching on $V(G)$. \endproof

\begin{lemma} \label{PMatchComplete} Let $G$ be an $r_{G}/n$-pseudorandom graph on $n$ vertices with $r_G \geq 2\log^2 n$, and let $H$ be a spanning subgraph of $G$ such that $\delta(H) \geq \frac{29}{30}r_G$. Let $P$ be a perfect matching (if $n$ is even) or a perfect pseudomatching (if $n$ is odd) on $V(G)$, edge-disjoint from $H$. Then $H$ contains a matching $P'$ such that $P \cup P'$ forms a $2$-regular graph with at most $\frac{3n \sqrt{r_G \log n}}{\delta(H)}$ cycles.
\end{lemma}

We will apply Lemma \ref{PMatchComplete} in the proof of Lemma \ref{PMatchSplit} with $P$ being one of the perfect matchings $M_i$ obtained from Lemma \ref{HalfMatch}.

\proof Let $r_H = \delta(H)$ and note that $r_H \geq \log^2 n$ by Definition \ref{PseudoRandom}(c). If $n$ is odd, then let $v_1$ be the vertex of $P$ of degree $2$ and $v_2, v_3$ its neighbours. For the remainder of the proof we treat $v_2 v_3$ as if it were an edge of $P$, and $v_1 v_2$, $v_1 v_3$ as if they were not edges of $P$. If $v_2 v_3$ is also an edge of $H$, then we delete it from $H$ (so that $P$ and $H$ remain edge-disjoint). The proof then proceeds identically whether $n$ is odd or even.

Label the edges of $P$ as $e_1, e_2, \ldots, e_{\lfloor n/2 \rfloor}$. For each edge $e_i$, label one of the vertices $a_i$ (chosen at random and independently from the other choices) and the other $b_i$. Let $B_1 = \{a_1, \ldots, a_{\lfloor n/2 \rfloor}\}$ and $B_2 = \{b_1, \ldots, b_{\lfloor n/2 \rfloor}\}$, and let $B$ be the bipartite graph on vertex classes $B_1, B_2$ with edges $E_H(B_1, B_2)$.

Now for each $v \in V(B)$, let $Y_v$ be the set of edges $ab \in P$ such that both $a$ and $b$ lie in $N_H(v)$, and let $Z_v$ be the set of edges $ab \in P$ such that exactly one of $a$ and $b$ lies in $N_H(v)$. Note that 
$$d_H(v) - 1 \leq 2|Y_v| + |Z_v| \leq d_H(v) \leq \Delta(H) \leq \Delta(G) \leq r_G + 2\sqrt{r_G \log n} \leq \frac{11r_H}{10},$$
where the second-last inequality follows from Definition \ref{PseudoRandom}(c). Now $d_B(v) \sim |Y_v| + Bin(|Z_v|, 1/2)$. (Here we use the condition that $P$ and $H$ are edge-disjoint and so the unique neighbour of $v$ in $P$ is not a neighbour of $v$ in $H$.)

Suppose $|Z_v| \geq 8 \log n$. Let $\varepsilon = 4 \sqrt{\log n/|Z_v|} \leq 3/2$ and note that $4\sqrt{r_H \log n} \geq \varepsilon |Z_v|/2$. By Lemma~\ref{cor23}(ii) we have that
$$\mathbb{P}\left[Bin(|Z_v|, 1/2) \leq |Z_v|/2 - 4\sqrt{r_H \log n}\right] < e^{-\varepsilon^2 |Z_v|/6} \leq 2e^{-16 \log n/6} < \frac{1}{n^2}.$$
On the other hand if $|Z_v| \leq 8 \log n$ then $|Z_v|/2 \leq 4\sqrt{r_H \log n}$ and so $$\mathbb{P}\left[Bin(|Z_v|, 1/2) \leq |Z_v|/2 - 4\sqrt{r_H \log n}\right] = 0.$$

So whp every $v \in V(B)$ satisfies $d_B(v) \geq |Y_v| + |Z_v|/2 - 4\sqrt{r_H \log n} \geq r_H/2 - 5\sqrt{r_H \log n}$. Similarly, whp every $v \in V(B)$ satisfies $d_B(v) \leq 11r_H/20 + 5\sqrt{r_H \log n}$. Now we choose $B$ such that $\delta(B) \geq r_H/2 - 5\sqrt{r_H \log n}$ and $\Delta(B) \leq 11r_H/20 + 5\sqrt{r_H \log n}$. Let $r_B = 2r_H/5$ and note that $r_B \geq r_G/3$. \medskip

\noindent \textbf{Claim:} \emph{$B$ contains a regular spanning subgraph $B'$ of degree $r_B$.} \medskip

\noindent To prove the claim, we use the Max-flow Min-cut theorem. Let each edge of $B$ have capacity $1$. Add a source $\sigma$ joined to each vertex of $B_1$ by an edge of capacity $r_B$ and a sink $\tau$ joined to each vertex of $B_2$ by an edge of capacity $r_B$. Let $n' = |B_1| = \lfloor n/2 \rfloor$. Now we show that the minimum cut must have capacity at least $n'r_B$ (indeed exactly $n'r_B$, since one could cut all of the edges incident to $\sigma$). It follows that the maximum flow from $\sigma$ to $\tau$ is $n'r_B$. Further there is some flow which achieves this maximum and such that the flow along each edge is an integer. Taking the edges of $B$ which have flow $1$ immediately yields the desired $r_B$-regular spanning subgraph.

So consider a cut $\mathcal{C}$. Let $T$ be the set of vertices $v \in B_1$ such that $\sigma v$ is not part of $\mathcal{C}$. Similarly let $S$ be the set of vertices $v \in B_2$ such that $v \tau$ is not part of $\mathcal{C}$. Now since every edge in $E_B(T, S)$ must be part of $\mathcal{C}$, the capacity of $\mathcal{C}$ is at least $e_B(T, S) + r_B(2n' - |S| - |T|)$. So noting that $e_B(T, S) = e_H(T, S)$, it suffices to prove that 
$$e_H(T, S) \geq n'r_B - r_B(2n' - |S| - |T|) = r_B(|S| + |T| - n').$$ 
Let $S' = B_2 \backslash S$. Then an equivalent statement is that 
$$e_H(T, B_2) - e_H(T, S') + r_B(|S'| - |T|) \geq 0.$$

We can assume without loss of generality that $|T| + |S'| \leq n' \leq n/2$; otherwise we can rearrange the inequality as $e_H(B_1, S) - e_H(T', S) + r_B(|T'| - |S|) \geq 0$ where $T' = B_1 \backslash T$ and the proof proceeds analogously. We now consider the following cases: \medskip

\noindent \textbf{Case 1:} $|T| \geq 2|S'|$. We can estimate $e_H(T, S') \leq |S'|\Delta(B)$ and $e_H(T, B_2) \geq |T|\delta(B)$. Now 
$$e_H(T, B_2) - e_H(T, S') + r_B(|S'| - |T|) \geq |T|(\delta(B) - r_B) - |S'|(\Delta(B) - r_B),$$
and the right-hand side is positive since $\delta(B) - r_B \geq r_H/10 - 5\sqrt{r_H \log n} > r_H/12$ and $\Delta(B) - r_B \leq 3r_H/20 + 5\sqrt{r_H \log n} < r_H/6$. \medskip

\noindent \textbf{Case 2:} $|S'| \geq |T|$. Then the inequality holds trivially. \medskip

\noindent \textbf{Case 3:} $|S'| \leq |T| \leq 2|S'|$. Note that $|T| \leq n/3$ (otherwise $|T| + |S'| > n/2$). So since $G$ is $(r_G/n, 2\sqrt{r_G})$-jumbled by Definition \ref{PseudoRandom}(c), (\ref{DisjointJumbled}) implies that
\begin{align*}
r_B |S'| - e_H(T, S') &\geq r_B |S'| - r_G|T||S'|/n - 4\sqrt{r_G}(|T| + |S'|) \\
&\geq |S'|(r_B - r_G|T|/n - 12\sqrt{r_G}) \geq |S'|(r_G/3 - r_G/3 - 12\sqrt{r_G}) \\
&= - 12|S'|\sqrt{r_G} \geq - 12|T|\sqrt{r_G}.
\end{align*} 
Hence $e_H(T, B_2) - e_H(T, S') + r_B(|S'| - |T|) \geq |T|\delta(B) - r_B|T| - 12|T|\sqrt{r_G} \geq 0$, which completes the proof of our claim. \medskip

Now Theorem~\ref{EFWThm} implies that the number of perfect matchings of $B'$ is at least $X = (\frac{r_B}{n'})^{n'} (n')!$. Note that for any $S \subseteq B_1$ and $T \subseteq B_2$ with $|S|, |T| \geq n'/\sqrt{r_G}$, 
\begin{align*}
e_{B'}(S, T) &\leq e_G(S, T) \stackrel{(\ref{DisjointJumbled})}{\leq} r_G|S||T|/n + 4\sqrt{r_G}(|S| + |T|) \\ 
&\leq r_G|S||T|/n' + 2r_G(|S| + |T|)^2/n' \\
&\leq 5r_G(|S| + |T|)^2/2n',
\end{align*}
where in the third inequality we use that $2 \leq \sqrt{r_G}(|S| + |T|)/n'$ and in the last inequality we use that $2|S||T| \leq (|S| + |T|)^2$. Set $C = 2n'\sqrt{r_G \log n'}/r_B$. Now applying Lemma~\ref{SizeOfP} with $V = B_1$, $V' = B_2$, and $v_i = a_i$ and $v'_i = b_i$ for each $1 \leq i \leq n'$ implies that $|\mathcal{P}_C(B')| < X$. So there exists a perfect matching $P'$ in $B'$ such that $D(P')$ has at most $C$ cycles. (Recall that $D(P')$ was defined in the paragraph before Lemma~\ref{SizeOfP}.) But we have a one-to-one correspondence between cycles $i_1 i_2 \ldots i_\ell i_1$ of $D(P')$ and cycles $b_{i_1} a_{i_1} b_{i_2} a_{i_2} \ldots b_{i_\ell} a_{i_\ell} b_{i_1}$ of $P \cup P'$. Hence $P \cup P'$ has at most $C$ cycles. Now note that $C \leq 3n \sqrt{r_G \log n}/r_H$. Finally, if $n$ is odd then on the cycle of $P \cup P'$ containing the edge $v_2 v_3$, we replace $v_2 v_3$ by the path $v_2 v_1 v_3$, so that $P \cup P'$ is $2$-regular. \endproof

We now combine Lemmas \ref{HalfMatch} and \ref{PMatchComplete} to show that under suitable conditions, the union of an arbitrary even-regular graph $H'$ and a graph $H$ which is close to being pseudorandom contains a collection of edge-disjoint $2$-factors which together cover the edges of $H'$. 

\begin{lemma} \label{PMatchSplit} Let $G$ be an $r_G/n$-pseudorandom graph on $n$ vertices with $r_G \geq 2\log^2 n$, and let $H$ be a spanning subgraph of $G$ such that $\delta(H) \geq (1 - 1/1100)r_G$. Let $H'$ be an arbitrary $r_{H'}$-regular graph on the same vertex set, and let $E_{bad} \subseteq E(H')$. Suppose that $H'$ is edge-disjoint from $H$ and that $r_{H'} + 1 + 10^6|E_{bad}|/n \leq r_G/5000$. Then there exists $m \leq r_G/5000$ and a collection $F_1, \ldots, F_{2m}$ of edge-disjoint $2$-factors in $H \cup H'$ whose union covers all of the edges of $H'$, such that each $F_i$ has at most $4n \sqrt{\log n}/\sqrt{r_G}$ cycles and each set $E(F_i) \cap E_{bad}$ is a matching of size at most $n/10^6$. \end{lemma}

We will apply Lemma \ref{PMatchSplit} with $E_{bad}$ being a set of `bad' edges which we will want to avoid when merging the cycles of each $F_i$ into a Hamilton cycle. The purpose of the restrictions on $F_i \cap E_{bad}$ is to spread the edges in $E_{bad}$ out among the $F_i$'s and thus make them easier to avoid.

\proof Let $r_H = (1 - 1/550)r_G + 5$. By Vizing's theorem, we can decompose $E(H')$ into edge-disjoint matchings $M'_1, M'_2, \ldots, M'_{r_{H'}+1}$. For each $i$, split $M'_i$ into $\lfloor 10^6|E(M'_i) \cap E_{bad}|/n + 1\rfloor$ matchings $M_j$ such that $|E(M_j) \cap E_{bad}| \leq n/10^6$. Let $M_1, \ldots, M_m$ be the resulting collection of matchings and note that $m \leq r_{H'} + 1 + 10^6|E_{bad}|/n \leq r_G/5000$. Now for each matching $M_j$, we will find a pair of edge-disjoint $2$-factors in $H \cup M_j$, each with at most $4n \sqrt{\log n}/\sqrt{r_G}$ cycles, which together cover the edges of $M_j$. Thus the total number of $2$-factors will be $2m$.

Firstly observe that for any $v \in V(H)$, no more than $4m - r_{H'} \leq r_G/1200$ edges of $H$ incident to $v$ will be used during this process to construct our $2$-factors (as each $2$-factor uses up at most $2$ edges of $H$ incident to $v$). Hence after deleting all the edges of $H$ lying in the $2$-factors found so far, we still have $\delta(H) \geq (1 - 1/1100 - 1/1200)r_G \geq r_H$. For each $M_j$, we use Lemma~\ref{HalfMatch} to decompose $M_j$ into two matchings and complete each one to a perfect matching or pseudomatching using edges of $H$. Now by Lemma~\ref{PMatchComplete} we can complete each of these perfect matchings or pseudomatchings to a $2$-factor using edges of $H$, such that each $2$-factor produced has at most $3n\sqrt{r_G \log n}/r_H \leq 4n \sqrt{\log n}/\sqrt{r_G}$ cycles. \endproof

Finally we combine Corollory \ref{2FacSplit} and Lemma \ref{PMatchSplit} to fully decompose $H \cup H'$ into $2$-factors with few cycles in total.

\begin{cor} \label{Absorb2Fac} Let $G$ be an $r_G/n$-pseudorandom graph with $r_G \geq 2\log^2 n$ and let $H$ be an even-regular spanning subgraph of $G$ of degree $r_H$, with $\delta(G) - 1 \leq r_H \leq \delta(G)$. Let $H'$ be an arbitrary even-regular graph of degree $r_{H'}$ on the same vertex set and let $E_{bad} \subseteq E(H')$. Suppose that $H'$ is edge-disjoint from $H$ and that $r_{H'} + 1 + 10^6|E_{bad}|/n \leq r_G/5000$. Let $m = (r_H + r_{H'})/2$. Then there exists a decomposition $\mathcal{F} = \{F_1, \ldots, F_m\}$ of $H \cup H'$ into $2$-factors such that $c(\mathcal{F}) \leq 4n\sqrt{r_G \log^3 n}$ and $E(F_i) \cap E_{bad}$ is a matching of size at most $n/10^6$ for each $i$.
\end{cor}

\proof Apply Lemma~\ref{PMatchSplit} to construct a collection $\mathcal{F'}$ of at most $r_G/2500$ edge-disjoint $2$-factors in $H \cup H'$, such that $\bigcup \mathcal{F}'$ covers $H'$,
$$c(\mathcal{F'}) \leq \frac{4n \sqrt{\log n}r_G}{2500\sqrt{r_G}} \leq n \sqrt{r_G \log n},$$
and the set $E(F) \cap E_{bad}$ is a matching of size at most $n/10^6$ for each $F \in \mathcal{F'}$. Let $H'' = (H \cup H') \backslash \bigcup \mathcal{F'}$ and note that $H''$ is an even-regular spanning subgraph of $G$. So by Corollary \ref{2FacSplit} with $H = H''$ we can decompose $H''$ into a collection $\mathcal{F''}$ of $2$-factors such that $c(\mathcal{F''}) \leq 3n\sqrt{r_G \log^3 n}$. Taking $\mathcal{F} = \mathcal{F'} \cup \mathcal{F''}$, we have $c(\mathcal{F}) \leq 4n\sqrt{r_G \log^3 n}$. \endproof

\section{Merging cycles} \label{merging}

So far we have the necessary tools to find large collections of disjoint $2$-factors in a pseudorandom graph $G$, by first finding a regular spanning subgraph and then decomposing this subgraph into a collection $\mathcal{F}$ of $2$-factors. The aim of this section is to show that we can transform these $2$-factors into Hamilton cycles, by `merging' the cycles of each $2$-factor together using edges which are taken from a pseudorandom graph $G'$. The crucial observation here is that $G'$ need not be as dense as $G$ originally was; in fact, under certain conditions $G'$ may be taken to be significantly sparser. To establish this we make use of the bounds on $c(\mathcal{F})$ proved in Section \ref{2FactorSection}.

Recall that $N(S)$ denotes the external neighbourhood of $S$, i.e. $N(S) = \bigcup_{s \in S} N(s) \backslash S$. All paths $P$ are considered to have a `direction' in the sense that the first and last endpoints are distinguished (so if $P = v_1 \ldots v_{\ell}$ and $P' = v_\ell \ldots v_1$ then we view $P$ and $P'$ as different paths). If $P = x \ldots y$ then we call $x$ the \emph{first endpoint} and $y$ the \emph{last endpoint} of $P$. We define the \emph{reverse} of a path $P = x \ldots y$ to be the path $P' = y \ldots x$ which has the same vertices and edges as $P$.

We will use the rotation-extension method. For this we need the following definitions.

\begin{defin}\label{rotation}
Let $P = v_1 v_2 \ldots v_\ell$ be a path. A \emph{rotation} of $P$ with \emph{pivot} $v_i$ is the operation of deleting the edge $v_{i-1} v_i$ from $P$ and adding the edge $v_i v_1$ to form a new path $v_{i-1} v_{i-2} \ldots v_1 v_i v_{i+1} \ldots v_{\ell}$ with endpoints $v_{i-1}$ and $v_\ell$. Call $v_{i-1} v_i$ the \emph{broken edge} of the rotation, and $v_i v_1$ the \emph{new edge}.
\end{defin}

\begin{defin}\label{reachable} Let $P$ be a path and $H$ a graph with $V(P) \subseteq V(H)$, such that $P$ and $H$ are edge-disjoint. Let $Q \subseteq V(H)$ and let $\tau \geq 1$ be an integer. Then a vertex $v$ of $P$ is $(H, Q, \tau)$-\emph{reachable} if there exists a sequence of at most $\tau$ rotations which make $v$ into the first endpoint of $P$, with all of the new edges being edges of $H$ and all of the pivots being elements of $Q$. 
\end{defin}

The main reason we include the set $Q$ in this definition is that there will be certain edges of $P$ that we do not want to break. We achieve this by ensuring that none of the endpoints of these edges lie in $Q$. Let
\begin{equation} \label{tau0}
\tau_0 = \frac{\log n}{\log \log n} + 3.
\end{equation}
In what follows, the sequences of rotations we consider will generally have length at most $\tau_0$.

\begin{defin} \label{extension}
Let $P = v_1 v_2 \ldots v_\ell$ be a path, and let $C = w_1 w_2 \ldots w_m w_1$ be a cycle which is vertex-disjoint from $P$. An \emph{extension of $P$ to incorporate $C$}, with \emph{join vertex} $w_i$ and \emph{broken edge} $w_{i-1} w_i$ is the operation of deleting the edge $w_{i-1} w_i$ from $C$ and adding the edge $v_1 w_i$ to form a new path $w_{i-1} w_{i-2} \ldots w_{i+1} w_{i} v_1 \ldots v_{\ell-1} v_{\ell}$. Call $v_1 w_i$ the \emph{new edge} of the extension.
\end{defin}

Given a $2$-regular graph $F$ and an edge-disjoint graph $H$ on the same vertex set, we say that a Hamilton cycle $C$ is formed by \emph{merging the cycles of $F$ using edges of $H$} if $E(C) \subseteq E(F) \cup E(H)$. Our strategy will be to first merge two of the cycles of $F$ together to form a long path $P$. We then show that if $Q \subseteq V(H)$ is `large', the set of $(H, Q, \tau_0)$-reachable vertices of $P$ must in general be large (see Corollary \ref{RotExp}). This will allow us to either extend the path by incorporating another cycle, or close $P$ to a cycle, thus reducing the number of cycles in $F$. Repeating this process will eventually produce a Hamilton cycle.

The following simple and well known lemma (see e.g., \cite[Lemma 23]{AHDoRG}, \cite[Proposition 3.1]{mk-bs-2003} for similar versions) will give us many reachable vertices provided that $H$ `expands', i.e., that $N_H(S) \cap Q$ is large (compared to $S$) for any $S \subseteq Q$ which is not itself too large. We include a proof here for completeness.

\begin{lemma}\label{expansion}
Let $H$ be a graph on $n$ vertices, and let $P = x \ldots y$ be a path on a subset of $V(H)$ which is edge-disjoint from $H$. Let $Q \subseteq V(P)$, and let $U_\tau$ be the set of $(H, Q, \tau)$-reachable vertices of $P$. Then $|U_{\tau+1}| \geq \frac{1}{2}|N_H(U_\tau) \cap Q| - |U_\tau|$.
\end{lemma}
\proof For a vertex $v \in P$, let $v^-$ and $v^+$ be the predecessor and successor, respectively, of $v$ along $P$. Let $T = \{v \in N_H(U_\tau) \cap Q \mid v^-, v^+ \notin U_\tau\}$. If $v \in T$, then since neither $v$ nor either of its neighbours on $P$ are in $U_\tau$, the neighbours of $v$ are preserved by every sequence of at most $\tau$ rotations of $P$ with pivots in $Q$; i.e., $v^+$ and $v^-$ are adjacent to $v$ along any path obtained from $P$ by at most $\tau$ rotations with pivots in $Q$. It follows that one of $v^-$ and $v^+$ must be in $U_{\tau+1}$. Indeed, starting from $P$, we can obtain by performing at most $\tau$ rotations a path with endpoints $z$ and $y$, such that $z \in U_\tau$ and $zv$ is an edge of $H$. Now by one further rotation with pivot $v$ and broken edge either $vv^+$ or $vv^-$, we obtain a path whose endpoints are either $v^+, y$ or $v^-, y$.

Now let $T^+ = \{v^+ \mid v \in T, v^+ \in U_{\tau+1}\}$ and $T^- = \{v^- \mid v \in T, v^- \in U_{\tau+1}\}$. 
It follows from the above that either $|T^+| \geq |T|/2$ or $|T^-| \geq |T|/2$, and both $T^+$ and $T^-$ are subsets of $U_{\tau+1}$. Hence $|U_{\tau+1}| \geq |T|/2 \geq (|N_H(U_\tau) \cap Q| - 2|U_\tau|)/2$, where the last inequality is immediate from the definition of $T$. \endproof

The fact that we have the required expansion property will follow from Lemma~\ref{AltExp}. Lemma~\ref{AltExp} relies on Corollary \ref{GnpExp}, which in turn relies on Lemma~\ref{DevBounds}. Lemma~\ref{DevBounds} is a simple consequence of pseudorandomness.

\begin{lemma} \label{DevBounds} Let $G$ be a $p$-pseudorandom graph on $n$ vertices, and let $S, T \subseteq V(G)$ be disjoint with $s = |S|$ and $t = |T|$. Let 
$$g(s, t) = \begin{cases}
2(s+t) \log n &\text{ if }\frac{\log n}{sp} \geq \frac{7}{2};\\
7s(s+t)p &\text{ otherwise},
\end{cases}$$ 
and 
$$h(s) = \begin{cases}
2s \log n &\text{ if }\frac{\log n}{sp} \geq \frac{7}{2};\\
7s^2 p &\text{ otherwise}.
\end{cases}$$ 
Then $e_G(S, T) \leq g(s, t)$ and $e_G(S) \leq h(s)$.
\end{lemma}

\proof Suppose first that $\frac{\log n}{sp} \geq \frac{7}{2}$. Then Definition \ref{PseudoRandom}(b)(i) implies that $e_{G}(S, T) \leq 2(s + t)\log n = g(s, t)$. If $\frac{\log n}{sp} \leq \frac{7}{2}$ and $\left(\frac{1}{s} + \frac{1}{t}\right)\frac{\log n}{p} \geq \frac{7}{2}$ then Definition \ref{PseudoRandom}(b)(i) implies that $e_G(S, T) \leq 2(s + t)\log n \leq 7s(s+t)p = g(s, t)$. Finally if $\left(\frac{1}{s} + \frac{1}{t}\right)\frac{\log n}{p} \leq \frac{7}{2}$ then Definition \ref{PseudoRandom}(b)(ii) implies that $e_G(S, T) \leq 7stp \leq 7s(s+t)p = g(s, t)$.

For the second part, if $\frac{7}{4} \leq \frac{\log n}{sp} \leq \frac{7}{2}$, then by Definition \ref{PseudoRandom}(b)(iii) we have $e_G(S) \leq 2s\log n \leq 7s^2 p = h(s)$, and otherwise the result follows immediately from Definition \ref{PseudoRandom}(b)(iii) and (iv). \endproof

\begin{cor} \label{GnpExp} Let $G$ be a $p$-pseudorandom graph on $n$ vertices. Let $S \subseteq V(G)$ with $|S| = s$, and for each vertex $x \in S$, let $A_x \subseteq N_G(x)$. Let $T = \bigcup_{x \in S} A_x \backslash S$, and $t = |T|$. Then the following properties hold:
\begin{itemize}
\item[{\rm (i)}] If $s \leq \frac{2\log n}{7p}$ and $\sum_{x \in S} |A_x| \geq 12s \log n$, then $t \geq \frac{\sum_{x \in S} |A_x|}{4 \log n}$.
\item[{\rm (ii)}] If $s \geq \frac{2\log n}{7p}$, then $t + 3s \geq \frac{\sum_{x \in S} |A_x|}{7sp}$.
\end{itemize}
\end{cor}

\proof ~(i) By Lemma~\ref{DevBounds} we have
$$\sum_{x \in S} |A_x| \leq e_G(S, T) + 2e_G(S) \leq g(s, t) + 2h(s) = 2(s+t) \log n + 4s \log n.$$
Hence $t \geq (\sum_{x \in S} |A_x| - 6 s\log n)/2\log n \geq \sum_{x \in S} |A_x|/4 \log n$.

~(ii) We have
$$\sum_{x \in S} |A_x| \leq e_G(S, T) + 2e_G(S) \leq g(s, t) + 2h(s) = 7sp(s+t) + 14s^2p$$
and the result follows immediately. \endproof

\begin{lemma} \label{AltExp} Let $0 < \varepsilon < 1$ be a constant. Let $G'$ be an $r_{G'}/n$-pseudorandom graph on $n$ vertices and let $G$ be an $r_G/n$-pseudorandom spanning subgraph of $G'$, with $r_G \leq r_{G'}$ and $\varepsilon r_{G} \geq 16 \log^2 n$. Let $H$ be a spanning subgraph of $G$ such that 
\begin{equation} \label{EGEHL38}
|E(G)\backslash E(H)| \leq 2n \sqrt{r_G \log n}. 
\end{equation}
Let $H'$ be a spanning subgraph of $G'$, such that
\begin{equation} \label{AltExpBound}
|E(H) \backslash E(H')|  \leq \frac{\varepsilon^2 nr_G^2}{10^4 r_{G'}} ~\emph{and}~ |E(H') \backslash E(H)| \leq \frac{\varepsilon^2 nr_G^2}{10^4 r_{G'}}.
\end{equation}
Let $Q', S \subseteq V(G)$ and suppose that $|N_{H'}(x) \cap Q'| \geq \varepsilon r_{G}$ for every vertex $x \in S$. Then at least one of the following holds:
\begin{itemize}
\item [\rm (i)] $\frac{1}{2}|N_{H'}(S) \cap Q'| - |S| \geq |S|\log n$,
\item [\rm (ii)] $|S| \leq \frac{\varepsilon nr_{G}}{50r_{G'}}$ and $\frac{1}{2}|N_{H'}(S) \cap Q'| - |S| \geq \frac{\varepsilon nr_{G}}{49r_{G'}}$, 
\item [\rm (iii)] $|S| \leq \varepsilon n/90$ and $\frac{1}{2}|N_{H'}(S) \cap Q'| - |S| \geq \varepsilon n/45$, 
\item [\rm (iv)]$|S| > |Q'|/6$,
\item [\rm (v)] $\frac{1}{2}|N_{H'}(S) \cap Q'| - |S| \geq |Q'|/6$.
\end{itemize}
\end{lemma}

Note that if $r_G$ is much smaller than $r_{G'}$ then (\ref{AltExpBound}) is more restrictive than simply requiring the symmetric difference of $E(H)$ and $E(H')$ to be $o(e(H))$. We need this more restrictive bound in Cases 3 and 4 below to obtain good expansion for small sets. (See also the remark after Corollary \ref{MergeCycles}.)

\proof Let $T = N_{H'}(S) \cap Q'$, $s = |S|$ and $t = |T|$. For each vertex $x \in V(G)$, let $B_x = N_{H'}(x) \cap Q'$ and note that
\begin{equation} \label{BxBound}
\sum_{x \in S} |B_x| \geq \varepsilon r_{G} s \geq 16 s \log^2 n.
\end{equation} 
We consider the following cases: \medskip

\noindent \textbf{Case 1:} $s \leq \frac{2 n \log n}{7r_{G'}}$. In this case we apply Corollary \ref{GnpExp}(i) with $G = G'$ and $A_x = B_x$ to obtain $t \geq \frac{\sum_{x \in S} |A_x|}{4 \log n} \geq \frac{16s \log^2 n}{4 \log n} = 4s\log n$. Hence (i) holds. \medskip

\noindent \textbf{Case 2:} $\frac{2 n \log n}{7r_{G'}} \leq s \leq \frac{\varepsilon nr_{G}}{50r_{G'}}$. Apply Corollary \ref{GnpExp}(ii) with $G = G'$ and $A_x = B_x$ to obtain 
$$t + 3s \geq \frac{n \sum_{x \in S} |A_x|}{7sr_{G'}} \geq \frac{\varepsilon nr_{G}}{7r_{G'}}.$$ 
Hence $t \geq \frac{4\varepsilon nr_{G}}{49r_{G'}}$ and so $\frac{1}{2}t - s \geq \frac{\varepsilon nr_{G}}{49r_{G'}}$, i.e., (ii) holds. \medskip

\noindent \textbf{Case 3:} $\frac{\varepsilon nr_{G}}{50r_{G'}} \leq s \leq \frac{2 n \log n}{7r_{G}}$. Let $A_x = B_x \cap N_H(x)$ and note that
$$\sum_{x \in S} |B_x| \stackrel{(\ref{BxBound})}{\geq} \varepsilon r_G s \geq \varepsilon r_{G} \frac{\varepsilon nr_{G}}{50r_{G'}} \stackrel{(\ref{AltExpBound})}{\geq} 8|E(H') \backslash E(H)|.$$
Hence $\sum_{x \in S} |A_x| \geq \sum_{x \in S} |B_x| - 2|E(H') \backslash E(H)| \geq \frac{3}{4}\sum_{x \in S} |B_x| \geq 12s \log^2 n$ by (\ref{BxBound}). So Corollary \ref{GnpExp}(i) with $G = G$ implies that $t \geq \frac{\sum_{x \in S}|A_x|}{4 \log n} \geq 3s \log n$, and hence (i) holds. \medskip

\noindent \textbf{Case 4:} $\max\{\frac{\varepsilon nr_{G}}{50r_{G'}}, \frac{2 n \log n}{7r_{G}}\} \leq s \leq \varepsilon n/90$. Let $A_x = B_x \cap N_H(x)$, and note similarly to Case 3 that $\sum_{x \in S} |A_x| \geq \frac{3}{4}\sum_{x \in S} |B_x|$. Hence Corollary \ref{GnpExp}(ii) with $G = G$ implies that 
$$t + 3s \geq \frac{n\sum_{x \in S} |A_x|}{7sr_G} \geq \frac{n\sum_{x \in S} |B_x|}{10sr_G} \stackrel{(\ref{BxBound})}{\geq} \frac{\varepsilon n}{10}.$$
So $t \geq \varepsilon n/15$, and $t/2 - s \geq \varepsilon n/45$. Hence (iii) holds. \medskip

\noindent \textbf{Case 5:} $s \geq \varepsilon n/90$. We may assume without loss of generality that $s \leq |Q'|/6$ (otherwise (iv) holds). In this case we must have $|Q'| \geq \varepsilon n/15$. \medskip

\noindent \textbf{Claim:} $|N_{H'}(S) \cap Q'| \geq 2|Q'|/3$. \medskip

\noindent Let $Q'' = Q' \backslash (N_{H'}(S) \cup S)$. To prove the claim, suppose for a contradiction that $|Q' \backslash N_{H'}(S)| \geq |Q'|/3$. Then 
$$|Q''| \geq |Q'|/3 - |S| \geq |Q'|/6 \geq s \geq \varepsilon n/90.$$ 
Now since $e_{H'}(S, Q'') = 0$, we have\COMMENT{Note that $n \sqrt{r_{G} \log n} \ll r_G n$, since $\sqrt{r_G} \geq \log n \gg \sqrt{\log n}$.}
\begin{align} \label{eSQBound}
e_{G}(S, Q'') \leq |E(G) \backslash E(H')| &\leq |E(G) \backslash E(H)| +|E(H) \backslash E(H')| \nonumber \\
&\stackrel{(\ref{EGEHL38}), (\ref{AltExpBound})}{\leq} 2n \sqrt{r_{G} \log n} + \frac{\varepsilon^2 nr_{G}^2}{10^4 r_{G'}} \leq \frac{\varepsilon^2 nr_G}{9000}.
\end{align}
But on the other hand Definition \ref{PseudoRandom}(a) and (\ref{DisjointJumbled}) imply that 
\begin{align*}
e_{G}(S, Q'') &\geq \frac{r_{G} s|Q''|}{n} - 4\sqrt{r_{G}}(s + |Q''|) \geq \frac{\varepsilon^2nr_G}{90^2} - 4n\sqrt{r_G} \geq \frac{\varepsilon^2 nr_G}{8500},
\end{align*}
contradicting (\ref{eSQBound}). This proves the claim. Now $\frac{1}{2}|N_{H'}(S) \cap Q'| - |S| \geq |Q'|/3 - |Q'|/6 = |Q'|/6$ and so (v) holds.\endproof

For a graph $H$ and a set $S \subseteq V(H)$, let $Int_{H}(S)$ be the set of vertices $x \in V(H)$ such that $x \in S$ and $N_H(x) \subseteq S$. If $S$ is a set of vertices such that $S \nsubseteq V(H)$, then we take $Int_H(S)$ to mean $Int_H(S \cap V(H))$. Further let $Cl_H(S) = S \cup N_H(S)$. Note that $V(H) \backslash Int_H(S) = Cl(V(H) \backslash S)$, and that if $P$ is a path then it is possible for an endpoint $x$ of $P$ to be in $Int_{P}(S)$.

Roughly speaking, the following lemma states that rotations (and extensions) of a path $P$ do not have much effect on $Int_P(S)$, provided that the pivots (or join vertices) themselves lie in $Int_P(S)$.

\begin{lemma} \label{PreserveNhood} Let $P = x \ldots y$ be a path, and let $Q$ be a set of vertices.
\begin{itemize}
\item [\rm (i)] Let $z \in Int_P(Q)$. Suppose we perform a rotation of $P$ with pivot $z$, and let $P'$ be the resulting path. Then $Int_P(Q) \backslash \{z\} \subseteq Int_{P'}(Q) \subseteq Int_P(Q)$. Further, if $x \in Q$ then $Int_{P'}(Q) = Int_{P}(Q)$.
\item [\rm (ii)] Let $C$ be a cycle vertex-disjoint from $P$. Let $zz^-$ be an edge of $C$ and suppose that $z \in Int_C(Q)$. Suppose we perform an extension of $P$ to incorporate $C$ with join vertex $z$ and broken edge $zz^-$, and let $P'$ be the resulting path. Then $(Int_P(Q) \cup Int_C(Q)) \backslash \{z\} \subseteq Int_{P'}(Q) \subseteq Int_P(Q) \cup Int_C(Q)$. Further, if $x \in Q$ then $Int_{P'}(Q) = Int_P(Q) \cup Int_C(Q)$.
\end{itemize}
\end{lemma}

\proof ~(i) Let $z^-$ be the predecessor of $z$ along $P$, $z^{--}$ the predecessor of $z^-$, and $x^+$ the successor of $x$. The only vertices whose neighbourhoods change as a result of the rotation are $x, z$ and $z^-$. However, if $x \in Int_P(Q)$ then $x^{+} \in Q$ and thus $x \in Int_{P'}(Q)$ and vice versa (here we use that $z \in Q$). Similarly, $z^{-} \in Int_P(Q)$ and $z^{-} \in Int_{P'}(Q)$ each hold if and only if $z^{--} \in Q$. Since $z \in Int_P(Q)$, the first part of (i) follows. If $x \in Q$ then we have $z \in Int_{P'}(Q)$ and hence $Int_{P'}(Q) = Int_{P}(Q)$.

~(ii) The proof of (ii) is similar. \endproof

The following corollary gives, under fairly weak conditions, a lower bound on the number of $(H', Q, \tau_0)$-reachable vertices of a long path $P$. This allows us to make any one of a large number of vertices of $P$ into the first endpoint of $P$ via a short sequence of rotations. Further, it allows us to `avoid' a specified set (namely $V(P) \backslash Q$) while doing so. We will use this second property for two main purposes: Firstly in order to make sure that certain edges of our $2$-factor $F$ which we want to keep are \emph{not} broken during the process of transforming $F$ into a Hamilton cycle, and secondly when we want to prevent one half of $P$ from being affected by the rotations at all. In each of these cases we need to construct sets satisfying the conditions on $Q$ and $Q'$ in Corollary \ref{RotExp}; Lemmas \ref{CreateW} and \ref{CreateQ} accomplish this.

\begin{cor} \label{RotExp} Let $0 < \varepsilon < 1$ be a constant. Let $G'$ be an $r_{G'}/n$-pseudorandom graph on $n$ vertices and let $G$ be an $r_{G}/n$-pseudorandom spanning subgraph of $G'$, with $r_{G} \leq r_{G'}$ and $\varepsilon r_{G} \geq 16 \log^2 n$. Let $H$ be a spanning subgraph of $G$ such that (\ref{EGEHL38}) holds, and let $H'$ be a spanning subgraph of $G'$ such that (\ref{AltExpBound}) holds.

Let $P = x \ldots y$ be a path in $G'$ such that $P$ and $H'$ are edge-disjoint. Let $Q \subseteq V(P)$ with $x \in Q$ and let $Q' \subseteq V(G')$ be such that $Q' \cap V(P) \subseteq Int_P(Q)$. Suppose that $|N_{H'}(v) \cap Q'| \geq \varepsilon r_{G}$ for every vertex $v \in Q$. 

Then either
\begin{itemize}
\item [\rm (i)] there exists a $(H', Q', \tau_0)$-reachable vertex $v$ of $P$ which has a neighbour in $H'$ lying in $Q' \backslash V(P)$, or
\item [\rm (ii)] the set of $(H', Q', \tau_0)$-reachable vertices of $P$ has size at least $|Q'|/6$.
\end{itemize}
\end{cor}

\proof Suppose that ~(i) does not hold, i.e., that $N_{H'}(v) \cap Q' \subseteq V(P)$ for every $(H', Q', \tau_0)$-reachable vertex $v$ of $P$. For each $\tau \leq \tau_0$, let $U_\tau$ be the set of $(H', Q', \tau)$-reachable vertices of $P$. Assume for contradiction that (ii) does not hold either, i.e., that $|U_{\tau}| < |Q'|/6$ for all $\tau \leq \tau_0$. By Lemma~\ref{expansion} with $Q = Q' \cap V(P)$ we have that $|U_{\tau+1}| \geq \frac{1}{2}|N_{H'}(U_\tau) \cap Q'| - |U_\tau|$ for each $\tau < \tau_0$, since $N_{H'}(U_\tau) \cap Q' \cap V(P) = N_{H'}(U_\tau) \cap Q'$.

Note that any rotation of $P$ whose pivot lies in $Q' \cap V(P) \subseteq Int_P(Q)$ produces a path $P'$ whose first endpoint lies in $Q$. Thus $Int_{P'}(Q) = Int_{P}(Q)$ by Lemma~\ref{PreserveNhood}(i) and the same holds for all paths obtained by further rotations with pivots in $Q' \cap V(P)$. So $U_{\tau} \subseteq Q$ and hence $|N_{H'}(v) \cap Q'| \geq \varepsilon r_G$ for each $v \in U_{\tau}$. Now Lemma~\ref{AltExp} with $S = U_\tau$ implies that $|U_{\tau + 1}| \geq \min\{\frac{\varepsilon nr_{G}}{49r_{G'}}, \frac{\varepsilon n}{45}, |U_{\tau}| \log n\}$ for any $\tau < \tau_0$. (Here we use that conclusions (iv) and (v) of Lemma~\ref{AltExp} cannot hold since $|U_\tau| < |Q'|/6$ for all $\tau \leq \tau_0$.) Hence there must exist some 
$$\tau_1 \leq \frac{\log\left(\min\left\{\frac{\varepsilon nr_{G}}{49r_{G'}}, \frac{\varepsilon n}{45}\right\}\right)}{\log \log n} + 1 \leq \tau_0 - 1,$$
such that $|U_{\tau_1}| \geq \min\{\frac{\varepsilon nr_{G}}{49r_{G'}}, \frac{\varepsilon n}{45}\}$.

Suppose first that $\frac{\varepsilon n}{45} \geq \frac{\varepsilon nr_{G}}{49r_{G'}}$. Then $|U_{\tau_1}| \geq \frac{\varepsilon nr_{G}}{49r_{G'}}$, and hence Lemma~\ref{AltExp} with $S = U_\tau$ implies that $|U_{\tau + 1}| \geq \min\{\frac{\varepsilon n}{45}, |U_{\tau}| \log n\}$ for each $\tau$ such that $\tau_1 < \tau < \tau_0$. Hence there exists $\tau_2 < \tau_0$ such that $|U_{\tau_2}| \geq \frac{\varepsilon n}{45}$. But now setting $S = U_{\tau_2}$, none of the conclusions of Lemma~\ref{AltExp} can hold. So we obtain the desired contradiction.

On the other hand, if $\frac{\varepsilon n}{45} \leq \frac{\varepsilon nr_{G}}{49r_{G'}}$ then we already have $|U_{\tau_1}| \geq \frac{\varepsilon n}{45}$. So we derive a contradiction in a similar way. \endproof

The following lemma will be used to obtain a set of vertices in which the expansion property we require still holds, and which also contains no endpoints of any `bad' edges that we want to avoid while rotating. ($W$ will be the set of these endpoints.)

\begin{lemma} \label{CreateW} Let $G'$ be an $r_{G'}/n$-pseudorandom graph on $n$ vertices and let $G$ be an $r_{G}/n$-pseudorandom spanning subgraph of $G'$, with $300 \log^3 n \leq r_{G} \leq r_{G'}$. Let $H$ be an even-regular spanning subgraph of $G$ with degree $r_H$ such that $\delta(G) - 1 \leq r_H \leq \delta(G)$. Let $H'$ be an $r_H$-regular spanning subgraph of $G'$, such that 
$$|E(H) \backslash E(H')| \leq \frac{nr_{G}^2}{2500r_{G'} \log^2 n},$$ 
and let $F$ be a $2$-factor of $G'$ which is edge-disjoint from $H'$. Let $W \subseteq V(G)$ with $|W| \leq n/400$.

Then there exist sets $V' \subseteq V(G)$ and $V'' \subseteq V' \backslash W$, such that
\begin{itemize}
\item $|Int_F(V'')| \geq n - 6|W|$, 
\item $|Int_F(V')| \geq n - |W|/\log^2 n$, and
\item $|N_{H'}(v) \cap Int_F(V'')| \geq r_H/2$ for every $v \in V'$.
\end{itemize}
\end{lemma}

\proof Note that by Definition \ref{PseudoRandom}(c),
$$\delta(G) \geq r_G - 2\sqrt{r_G \log n} \geq r_G\left(1 - \frac{2}{\sqrt{300} \log n}\right) \geq 290 \log^3 n,$$
and hence $r_H \geq 288 \log^3 n$.

If $W = \emptyset$ then we can take $V'' = V' = V(G)$, so we assume hereafter that $|W| \geq 1$. We define $V''$ as follows: Initially $V'' = V(G) \backslash W$. As long as there exists a vertex $v \in V''$ such that $|N_{H'}(v) \cap Int_F(V'')| \leq r_H/2$, we remove $v$ from $V''$. Let $W'$ be the set of removed vertices and note that $|N_{H'}(v) \cap Cl_F(W \cup W')| \geq r_H/2$ for every $v \in W'$. Suppose that this process continues until $|W'| = |W|$. Then for each $x \in W'$, let $B_x = N_{H'}(x) \cap Cl_F(W \cup W')$. Thus we have $|B_x| \geq r_H/2$ for each $x \in W'$. 

Note that 
\begin{equation} \label{sumBx}
\sum_{x \in W'} |B_x| \geq |W'|r_H/2 \geq 144 |W'|\log^3 n,
\end{equation}
and also that $|Cl_F(W \cup W')| \leq 3(|W| + |W'|) = 6|W|$. Now we separate into two cases: \medskip

\noindent \textbf{Case 1:} $|W'| = |W| \leq n r_H/150 r_{G'}$. Applying Corollary \ref{GnpExp} with $G = G'$, $S = W'$ and $A_x = B_x$, we obtain either 
$$|Cl_F(W \cup W')| \geq \left|\bigcup_{x \in W'} B_x\right| \geq \frac{|W|r_H}{8 \log n} \geq 7|W|$$
or
$$|Cl_F(W \cup W')| \geq \left|\bigcup_{x \in W'} B_x\right| \geq \frac{|W| r_H n}{14r_{G'} |W|} - 3|W| \geq 7|W|,$$
either of which yields an immediate contradiction. \medskip

\noindent \textbf{Case 2:} $|W'| = |W| \geq n r_H/150 r_{G'}$. Let $A_x = B_x \cap N_H(x)$. Note that 
$$\sum_{x \in W'} |B_x| \geq \frac{|W'|r_H}{2} \geq \frac{r_H}{2}\cdot \frac{n r_H}{150r_{G'}} \geq 4|E(H)\backslash E(H')|,$$ 
and that 
$$\sum_{x \in W'} |A_x| \geq \sum_{x \in W'} |B_x| - 2|E(H)\backslash E(H')| \geq \frac{1}{2}\sum_{x \in W'} |B_x| \stackrel{(\ref{sumBx})}{\geq} 72 |W| \log^3 n.$$ 
So we can apply Corollary \ref{GnpExp} with $G = G$ and $S = W'$ to obtain either 
$$|Cl_F(W \cup W')| \geq \left|\bigcup_{x \in W'} A_x\right| \geq \frac{|W|r_H}{16 \log n} \geq 7|W|$$
or 
$$|Cl_F(W \cup W')| \geq \left|\bigcup_{x \in W'} A_x\right| \geq \frac{|W| r_H n}{28 r_G |W|} - 3|W| \geq \frac{n}{30} - 3|W| \geq 7|W|,$$
either of which again yields an immediate contradiction. \medskip

So the process must terminate before $|W'| = |W|$. Fix $V''$ in its state at the point when the process terminates. Note that $|V(G) \backslash V''| \leq 2|W|$ and hence $|Cl_F(V(G) \backslash V'')| \leq 6|W|$, i.e., $|Int_F(V'')| \geq n - 6|W|$.

Let $W'' = \{v \in V(G) \mid |N_{H'}(v) \cap Cl_F(W \cup W')| \geq r_H/2 \}$ and $V' = V(G) \backslash W''$. Since $W'' \subseteq W \cup W'$ we have $|W''| \leq 2|W|$. We will show that $V''$ and $V'$ satisfy the assertions of the lemma. Since $W'' \subseteq W' \cup W = V(G) \backslash V''$, we have $V'' \subseteq V' \backslash W$ and so it remains to prove that $|Int_F(V')| \geq n - |W|/\log^2 n$. This holds provided that
\begin{equation} \label{SizeOfClW''}
|Cl(W'')| \leq |W|/\log^2 n.
\end{equation}
We now claim that $|Cl_F(W \cup W')| \geq 18|W''| \log^2 n$. Again we consider two cases: \medskip

\noindent \textbf{Case 1:} $|W''| \leq n r_H/300 r_{G'} \log^2 n$. Then applying Corollary \ref{GnpExp} with $G = G'$, $S = W''$ and $A_x = B_x$, we have either
$$|Cl_F(W \cup W')| \geq \frac{\sum_{x \in W''} |B_x|}{4 \log n} \geq \frac{|W''|r_H}{8 \log n} \geq 18|W''| \log^2 n$$
or
$$|Cl_F(W \cup W')| \geq \frac{|W''| n r_H}{14r_{G'} |W''|} - 3|W''| \geq |W''|(19 \log^2 n - 3) \geq 18|W''| \log^2 n,$$ 
as desired. \medskip

\noindent \textbf{Case 2:} $|W''| \geq n r_H/300 r_{G'} \log^2 n$. Note that
$$\sum_{x \in W''} |B_x| \geq \frac{r_H}{2}\cdot\frac{nr_H}{300r_{G'}\log^2 n} \geq 4|E(H)\backslash E(H')|,$$
and so setting $A_x  = N_H(x) \cap B_x$ again we have 
$$\sum_{x \in W''} |A_x| \geq \sum_{x \in W''} |B_x| - 2|E(H)\backslash E(H')| \geq \frac{1}{2}\sum_{x \in W''} |B_x| \geq 72 |W''| \log^2 n.$$
Now again applying Corollary \ref{GnpExp}, we have either
$$|Cl_F(W \cup W')| \geq \frac{|W''|r_H}{16 \log n} \geq 18|W''|\log^2 n$$
or 
$$|Cl_F(W \cup W')| \geq \frac{|W''| r_H n}{28 r_G |W''|} - 3|W''| \geq \frac{n}{30} - 3|W''| \geq 13|W| - 6|W| \geq 7|W|.$$
But the latter case cannot occur since $|Cl_F(W \cup W')| \leq 3|W \cup W'| \leq 6|W|$. This proves the claim. Hence $|W''| \leq \frac{|Cl_F(W \cup W')|}{18 \log^2 n} \leq |W|/3\log^2 n$ and so $|Cl_F(W'')| \leq |W|/\log^2 n$, which proves (\ref{SizeOfClW''}). \endproof

Roughly speaking, the following lemma states that if we have a graph $H'$ which is close to being pseudorandom and a long path $P$ subdivided into $\log n$ segments, then most of the vertices will have many neighbours in most of the segments. This will enable us to carry out rotations involving only the initial half of a long path in the proof of Lemma \ref{RotExt}. The proof of Lemma \ref{CreateQ} is similar to that of \cite[Lemma 22]{AHDoRG}.

\begin{lemma} \label{CreateQ} Let $G$ be an $r_{G}/n$-pseudorandom graph with $r_G \geq 10^5 \log^2 n$. Let $\varepsilon \leq 1/5$ be a positive constant and let $n'$ be an integer such that $n/10 \leq n' \leq n$.  Let $U \subseteq V(G)$ be such that $|V(G) \backslash U| \leq \varepsilon^2 n'/8$, and let $H'$ be a graph on $V(G)$ such that 
\begin{equation} \label{EGEH'}
|E(G) \backslash E(H')| \leq \frac{\varepsilon^3 r_{G} n}{32000}.
\end{equation} 

Let $P$ be a path which is edge-disjoint from $H'$, with $V(P) \subseteq V(G)$ and $|P| = n'$, divided into $\log n$ segments $J_1, \ldots, J_{\log n}$ whose lengths are as equal as possible. Then there exists a set $I \subseteq [\log n]$ such that $|I| \geq (1 - \varepsilon)\log n$ and for every $i \in I$, there exists $J'_i \subseteq J_i \cap U$ such that $|Int_P(J'_i)| \geq (1 - \varepsilon)n'/\log n$, which satisfies the following condition: For every $v \in J'_i$, and for all but at most $\varepsilon \log n$ indices $j \in I$, $|N_{H'}(v) \cap Int_P(J'_j)| \geq \frac{r_{G}}{25 \log n}$.
\end{lemma}
\proof We define $I$ and $\{J'_i \mid i \in I\}$ as follows: Initially $I = [\log n]$ and $J'_i = J_i \cap U$. For a vertex $v$ and for $j \in I$, call $v$ \textit{weakly connected} to $J'_j$ if $|N_{H'}(v) \cap Int_P(J'_j)| \leq \frac{r_{G}}{25 \log n}$. As long as there exists $i \in I$ and a vertex $v \in J'_i$, such that $v$ is weakly connected to $J'_j$ for more than $\varepsilon \log n$ values of $j \in I$, we remove $v$ from $J'_i$. 
Further, if at any stage there exists $i \in I$ such that $|Int_P(J'_i)| \leq (1-\varepsilon)n'/\log n$, then we remove $i$ from $I$. 

We claim that this process must terminate before $\varepsilon^2 n'/4 - |V(P) \backslash U|$ vertices are removed. Indeed, suppose we have removed $\varepsilon^2 n'/4 - |V(P) \backslash U|$ vertices and let $R$ be the set of removed vertices. Note that no vertices of $V(P) \backslash U$ are included in $R$, and that $|R| \geq \varepsilon^2 n'/8$ since $|V(P) \backslash U| \leq |V(G) \backslash U| \leq \varepsilon^2 n'/8$. Now $|R \cup (V(P) \backslash U)| = \varepsilon^2 n'/4$, and so $\sum_{i=1}^{\log n} |Int_P(J'_i)| \geq (1-3\varepsilon^2/4)n'$. Hence we have $|Int_P(J'_i)| \geq (1-\varepsilon)n'/\log n$ for at least $(1-3 \varepsilon/4)\log n$ values of $i$, i.e., at most $3 \varepsilon \log n/4$ indices have been removed from $I$. So for each vertex $v \in R$, there remain more than $\varepsilon \log n/4$ indices $i \in I$ for which $v$ is weakly connected to $J'_i$. For each $i \in I$, let $WC(i)$ be the set of vertices $v \in R$ which are weakly connected to $J'_i$. Let $I_0 = \{i \in I \mid |WC(i)| \geq \varepsilon^3 n'/64\}$. \medskip

\noindent \textbf{Claim:} \emph{For each $i \in I_0$, there are at least $\frac{|WC(i)|r_{G}}{50 \log n}$ edges in $E(G) \backslash E(H')$ between $Int_P(J'_i)$ and $WC(i)$.} \medskip

\noindent Let $S = Int_P(J'_i)$ and note that $|S| \geq (1 - \varepsilon)n'/\log n$, since $i \in I$. 
To prove the claim, note that $S \cap WC(i) = \emptyset$ since $WC(i) \subseteq R$, and hence (\ref{DisjointJumbled}) implies that 
\begin{align*}
e_{G}(S, WC(i)) &\geq \frac{|S||WC(i)|r_{G}}{n} - 4\sqrt{r_{G}}(|S| + |WC(i)|) \\
&\geq \frac{(1 - \varepsilon)n'|WC(i)|r_{G}}{n \log n} - 5\sqrt{r_{G}}|WC(i)| \\
&\geq \frac{|WC(i)|r_G}{50 \log n}\left(4 - \frac{250 \log n}{\sqrt{r_G}}\right) \geq \frac{3|WC(i)|r_{G}}{50 \log n},
\end{align*}
where the last inequality follows from the condition $r_G \geq 10^5 \log^2 n$. But the definition of $WC(i)$ implies that $e_{H'}(S, WC(i)) \leq \frac{|WC(i)|r_{G}}{25 \log n}$, and the claim follows immediately. \medskip

Observe that $\sum_{i \in I} |WC(i)| > |R|\varepsilon \log n/4 \geq \varepsilon^3 n'\log n/32$, since each $v \in R$ is weakly connected to $J'_i$ for more than $\varepsilon \log n/4$ values of $i \in I$. But now
$$\sum_{i \in I \backslash I_0} |WC(i)| \leq |I\backslash I_0|\frac{\varepsilon^3 n'}{64} \leq \frac{\varepsilon^3 n'\log n}{64},$$
and hence $\sum_{i \in I_0} |WC(i)| > \varepsilon^3 n'\log n/64$. Together with the claim, this implies that 
$$|E(G) \backslash E(H')| \geq \frac{r_{G}}{50 \log n} \sum_{i \in I_0} |WC(i)| > \frac{\varepsilon^3 r_{G} n' \log n }{3200 \log n} \geq \frac{\varepsilon^3 r_{G} n}{32000},$$ 
which contradicts (\ref{EGEH'}). This proves that the process must terminate before $\varepsilon^2 n'/4 - |V(P) \backslash U|$ vertices are removed. Now $|I| \geq (1 - \varepsilon) \log n$ at the point at which the process terminates, since as before at most $3 \varepsilon \log n/4$ indices have been removed, and so at this point $I$ and $\{J'_i \mid i \in I\}$ satisfy the assertions of the lemma. \endproof

In the next lemma we exhibit a method for merging cycles of a $2$-factor $F$ together. The basic idea is fairly simple: given $F$, we first use rotation-extension to transform some of the cycles in $F$ into a long path. When we can no longer extend the path, we apply Lemma \ref{CreateQ} to show that there is a large set $Q_1$ on the first half of the path whose vertices are reachable (by rotations which only use vertices in the first half of the path as pivots), and a similar set $Q_2$ in the last half of the path. Since $Q_1$ and $Q_2$ are large we can find an edge between them and close the path to a cycle, thus forming a new $2$-factor $F'$ with fewer cycles. However, the fact that we need to avoid bad edges forces us to be very careful.

\begin{lemma} \label{RotExt} Let $G'$ be an $r_{G'}/n$-pseudorandom graph on $n$ vertices, and let $G$ be an $r_{G}/n$-pseudorandom spanning subgraph of $G'$ with $300 \log^3 n \leq r_{G} \leq r_{G'}$. Let $F$ be a $2$-factor of $G'$ and let $H$ be an even-regular spanning subgraph of $G$ of degree $r_H$ such that $\delta(G) - 1 \leq r_H \leq \delta(G)$. Let $H'$ be an $r_H$-regular spanning subgraph of $G'$, such that 
\begin{equation} \label{RotExtEdgeBound}
|E(H) \backslash E(H')| = |E(H') \backslash E(H)| \leq \frac{nr_{G}^2}{2500 r_{G'} \log^2 n}.
\end{equation}
Let $E_{bad} \subseteq E(F)$ be a matching of size at most $n/220000$. Suppose that $F$ and $H'$ are edge-disjoint and that $2|E_{bad}|/\log^2 n \leq c(F) \leq \frac{4nr_G^2}{r_{G'} \log^3 n}$.

Then unless $F$ is a Hamilton cycle, we can obtain a new $2$-factor $F'$ with the following properties:
\begin{itemize}
\item $E(F') \subseteq E(F) \cup E(H')$,
\item $c(F') < c(F)$,
\item $|E(F') \cap E(H')| \leq \frac{5 \log n}{\log \log n}(c(F) - c(F'))$, and
\item $E_{bad} \subseteq E(F')$.
\end{itemize}
\end{lemma}

\proof Let $W$ be the set of endpoints of edges in $E_{bad}$. Note that $|W| \leq n/110000$ and so by Lemma~\ref{CreateW} there exist sets $V' \subseteq V(G)$ and $V'' \subseteq V' \backslash W$, with $|Int_F(V'')| \geq n - 6|W|$ and $|Int_F(V')| \geq n - |W|/\log^2 n$ and such that 
\begin{equation} \label{vxnhood}
|N_{H'}(v) \cap Int_F(V'')| \geq r_H/2.
\end{equation}
for every vertex $v \in V'$. Note that 
\begin{equation} \label{cF}
c(F) \geq 2|E_{bad}|/\log^2 n = |W|/\log^2 n \geq n - |Int_F(V')|.
\end{equation}
We will obtain $F'$ by deleting an edge $yy'$ in a cycle of $F$ to form a path $P$, performing a sequence of rotations and extensions to incorporate other cycles of $F$ into $P$, and finally closing $P$ to a cycle, thus reducing the number of cycles. The new edges for these rotations and extensions will be taken from a graph $H''$, which is defined as follows: Initially $H'' = H' \cup \{yy'\}$. Whenever a rotation or extension with new edge $e \in H''$ is performed, we remove $e$ from $H''$ and call $e$ a \emph{used} edge. For each edge of $H''$ used we will add back to $H''$ the edge of $F$ broken during the corresponding rotation or extension. 

We can view $H''$ as the `current' version of $H'$; note that we always have $|E(H'')| = |E(H')|$ or $|E(H'')| = |E(H')| + 1$. Similarly, let $F'$ be the current version of $F$ (so $F'$ is either a $2$-factor of $H' \cup F$, or consists of a path and various cycles). Note that since we do not change $F$ itself during the proof, the function $Int_F$ will not change either. 

As long as we use on average at most $5 \log n/\log \log n$ edges per cycle of $F$ merged, we have 
\begin{align} \label{EG'EH1}
|E(H'') \backslash E(H)| &\leq |E(H) \backslash E(H'')| + 1 \nonumber \\
&\leq |E(H) \backslash E(H')| + |E(H') \backslash E(H'')| + 1 \nonumber \\
&\stackrel{(\ref{RotExtEdgeBound})}{\leq} \frac{nr_{G}^2}{2500 r_{G'} \log^2 n} + \frac{5 c(F) \log n}{\log \log n} + 1 \leq \frac{2nr_G^2}{r_{G'} \log^2 n}.
\end{align}
This fulfils condition (\ref{AltExpBound}) of Lemma~\ref{AltExp} and Corollary \ref{RotExp} whenever they are applied with $H' = H''$. We also have that $e(G) \leq nr_{G}/2 + 2 n\sqrt{r_{G}}$ since $G$ is $(r_{G}/n, 2\sqrt{r_{G}})$-jumbled\COMMENT{We often use the fact that a $p$-pseudorandom graph $G$ is $(p, 2\sqrt{np})$-jumbled; maybe we should mention it for clarity?} by Definition \ref{PseudoRandom}(a), and 
$$e(H) = \frac{nr_H}{2} \geq \frac{n(\delta(G) - 1)}{2} \geq \frac{nr_G}{2} - n\sqrt{r_{G} \log n} - \frac{n}{2}$$ 
by Definition \ref{PseudoRandom}(c). So since $E(H) \subseteq E(G)$, 
\begin{equation} \label{EG'EH2}
|E(G) \backslash E(H)| \leq 2 n\sqrt{r_{G}} + n\sqrt{r_{G} \log n} + \frac{n}{2} \leq 2n \sqrt{r_{G} \log n}.
\end{equation}
This fulfils condition (\ref{EGEHL38}) of Lemma~\ref{AltExp} (and Corollary \ref{RotExp}). Further, still assuming that on average at most $5 \log n/\log \log n$ edges are used per cycle merged, we have
\begin{align}
|E(G) \backslash E(H'')| &\leq |E(G) \backslash E(H)| + |E(H) \backslash E(H'')| \nonumber \\ 
&\stackrel{(\ref{EG'EH1}),(\ref{EG'EH2})}{\leq} 2n\sqrt{r_{G} \log n} + \frac{2nr^2_{G}}{r_{G'} \log^2 n} \leq \left(\frac{1}{15}\right)^3\frac{r_{G}n}{32000}. \label{EG'EH}
\end{align} \medskip

\noindent \textbf{Claim:} \emph{$H'[Int_F(V')]$ is connected.} \medskip

\noindent To prove the claim, suppose for contradiction that $H'[Int_F(V')]$ has two components $S$ and $T$. Since (\ref{vxnhood}) implies that $|N_{H'}(v) \cap Int_F(V')| \geq r_H/2$ for any $v \in Int_F(V')$, we can apply Lemma~\ref{AltExp} with $\varepsilon = 1/3$, $Q' = Int_F(V')$ and $S = S$. But $N_{H'}(S) \cap Q' = \emptyset$, and so of the possible conclusions of Lemma~\ref{AltExp} only (iv) can hold. This implies that $|S| \geq |Q'|/6 \geq n/7$, and similarly $|T| \geq |Q'|/6 \geq n/7$. Now since $G$ is $(r_{G}/n, 2\sqrt{r_{G}})$-jumbled, (\ref{DisjointJumbled}) implies that 
$$e_{G}(S, T) \geq r_G|S||T|/n - 4\sqrt{r_{G}}(|S| + |T|) \geq nr_{G}/50.$$ 
But by (\ref{EG'EH}) this implies that $e_{H'}(S, T) > 0$. Hence $H'[Int_F(V')]$ has only one component, which proves the claim. \medskip

Our first aim is to find a path $P_1 = x'' \ldots y''$ whose vertices span at least two cycles of $F$ and so that $x''$ and $y''$ lie in $Int_F(V'')$. Let $C_1$ be a cycle of $F$ containing a vertex of $Int_F(V')$. It is easy to see that $Int_F(V') \backslash V(C_1)$ is nonempty. Indeed, if $V' = V(G)$ then it is sufficient to recall that $F$ is not a Hamilton cycle. Otherwise $|Int_F(V')| \leq n - 3$ and so by (\ref{cF}) the number of vertices which are not part of $C_1$ is at least $3(c(F) - 1) \geq 3(n - |Int_F(V')| - 1) > n - |Int_F(V')|$. Hence one of these vertices must be in $Int_F(V')$. (This is the only place in the proof where we use the lower bound on $c(F)$.)

Since $H'[Int_F(V')]$ is connected, there exists an edge $xy$ of $H'$ with $x \in Int_F(V') \cap V(C_1)$ and $y \in Int_F(V') \backslash V(C_1)$. Let $C_2$ be the cycle of $F$ in which $y$ lies. Since $E_{bad}$ is a matching, there will be at least one edge $yy' \notin E_{bad}$ in $C_2$ incident to $y$. Delete $yy'$ from $C_2$ to form a path $C'_2 = y \ldots y'$. Note that $Int_F(V'') \cap V(C'_2) \subseteq Int_{C'_2}(V'')$.\COMMENT{Since breaking the edge $yy'$ can only make it `easier' for $y$ and $y'$ to lie in $Int_{C'_2}(V'')$.}

Again since $E_{bad}$ is a matching we can find $xx' \notin E_{bad}$ in $C_1$ incident to $x$. Perform an extension of $C'_2$ to incorporate $C_1$, with join vertex $x$ and broken edge $xx'$. Let $P_0 = x' \ldots y'$ denote the resulting path and note that $(Int_F(V'') \cap V(P_0)) \backslash \{x, y\} \subseteq Int_{P_0}(V'')$.\COMMENT{Again, for vertices other than $x, y$ it only becomes easier to be in the interior.}

Now by (\ref{vxnhood}), $x'$ has a neighbour $z_1$ in $H''$, which lies in $Int_F(V'') \backslash \{x, y\}$. If $z_1 \in V(P_0)$ then we perform a rotation with pivot $z_1$; note that the first endpoint $x''$ of the resulting path lies in $V''$. Otherwise, let $C$ be the cycle containing $z_1$ and let $z_1 x'' \notin E_{bad}$ be an edge of $C$. Perform an extension of $P_0$ to incorporate $C$ and again note that the new first endpoint $x''$ will lie in $V''$. Call the resulting path $P'_0$ and let $P''_0 = y' \ldots x''$ be the reverse of $P'_0$. (Recall that the reverse of a path was defined in the second paragraph of Section \ref{merging}.) Now similarly we can perform a rotation (or extension) of $P''_0$ with some pivot $z_2$ (or join vertex $z_2$), so that the new first endpoint $y''$ also lies in $V''$. Call the resulting path $P'''_0$ and let $P_1 = x'' \ldots y''$ be the reverse of $P'''_0$. 

Applying Lemma~\ref{PreserveNhood} twice with $P = P_0$ and $P''_0$ respectively, $Q = V''$ and $z = z_1$ and $z_2$ respectively implies that 
\begin{equation} \label{PresInt1}
(Int_F(V'') \cap V(P_1)) \backslash \{z_1, z_2, x, y\} \subseteq Int_{P_1}(V'').
\end{equation}
All rotations performed during the remainder of the proof will have pivots in $(Int_F(V'') \cap V(P_1)) \backslash \{z_1, z_2, x, y\}$ (where $P_1$ is the current path). Thus by (\ref{PresInt1}) the pivots lie in $Int_{P_1}(V'')$. Also the join vertices of all extensions performed during the remainder of the proof will lie in $Int_F(V'')$. Hence no edges in $E_{bad}$ will be broken. Moreover, the endpoints of $P_1$ will always lie in $V''$ and hence Lemma~\ref{PreserveNhood} implies that (\ref{PresInt1}) will always hold.

Note that for each $v \in V'$, we have 
$$|N_{H''}(v) \cap Int_F(V'')| \geq |N_{H'}(v) \cap Int_F(V'')| - 2$$ 
since $H' \subseteq H'' \cup F'$ and $\Delta(F') = 2$. So (\ref{vxnhood}) implies that
$$|N_{H''}(v) \cap Int_F(V'')| \geq r_G/3 + 4$$
with room to spare. We now incorporate additional cycles into $P_1$ by an iterative procedure. Each iteration proceeds as follows: \medskip

\noindent Apply Corollary \ref{RotExp} with $\varepsilon = 1/3$, $H' = H''$, $P = P_1$, $Q = V'' \cap V(P_1)$ and $Q' = Int_F(V'') \backslash \{z_1, z_2, x, y\}$. We obtain one of two possible conclusions: \medskip

\noindent \textbf{Case 1:} \emph{There exists an $(H'', Q', \tau_0)$-reachable vertex $v$ of $P_1$ which has a neighbour $z \in Int_F(V'') \backslash V(P_1)$.} Let $C$ be the cycle of $F$ on which $z$ lies. Now we perform the necessary rotations to make $v$ the first endpoint of $P_1$, and then an extension with join vertex $z$ to incorporate $C$ into $P_1$. We then redefine $P_1$ to be the resulting path and begin the next iteration. \medskip

\noindent \textbf{Case 2:} \emph{The number of $(H'', Q', \tau_0)$-reachable vertices of $P_1$ is at least $|Q'|/6 \geq 3n/20$.} Note that this immediately implies that $|P_1| \geq |Q'/6| \geq n/10$. Let $P_1$ be divided into $\log n$ segments $J_i$ whose lengths are as equal as possible. Noting that $|V(G) \backslash Q'| \leq 6|W| + 4 \leq n/18000$, we may apply Lemma~\ref{CreateQ} with $\varepsilon = 1/15$, $H' = H''$, $U = Q'$ and $n' = |P_1|$. Note that (\ref{EGEH'}) is satisfied due to (\ref{EG'EH}). Thus we obtain a set $I \subseteq [\log n]$ of size at least $14 \log n/15$ and sets $J'_i \subseteq J_i \cap Q'$ for each $i \in I$, such that $|Int_{P_1}(J'_i)| \geq 14|P_1|/15 \log n$ and the following holds: For every $v \in J'_i$, and for all but at most $\log n/15$ indices $j \in I$, $|N_{H''}(v) \cap Int_{P_1}(J'_j)| \geq \frac{r_{G}}{25 \log n}$. 

Note that $|\bigcup_{i \in I} Int_{P_1}(J'_i)| \geq (1 - 2/15)|P_1|$. Hence for some $i \in I$, $J'_i$ contains a $(H'', Q', \tau_0)$-reachable vertex of $P_1$. Perform the necessary rotations to make one of these vertices into the first endpoint, breaking at most $\tau_0$ edges in the process. Redefine $P_1$ to be the resulting path, and let $P'_1$ be the reverse of $P_1$. Now apply Corollary \ref{RotExp} with $\varepsilon = 1/3$, $H' = H''$, $P = P'_1$, $Q = V'' \cap V(P'_1)$ and $Q' = Int_F(V'') \backslash \{z_1, z_2, x, y\}$. Again we obtain one of two possible conclusions: \medskip

\noindent \textbf{Case 2a:} \emph{There exists an $(H'', Q', \tau_0)$-reachable vertex $v$ of $P'_1$ which has a neighbour $z \in Int_F(V'') \backslash V(P'_1)$.} In this case we extend $P'_1$ as in Case 1, redefine $P_1$ to be the resulting path and begin the next iteration. Note that in future instances of Case 2, the sets $I$, $J_i$ and $J'_i$ will be redefined for each new path $P'_1$. \medskip

\noindent \textbf{Case 2b:} \emph{The number of $(H'', Q', \tau_0)$-reachable vertices of $P'_1$ is at least $|Q'|/6 \geq 3n/20$.} Note that any instance of Case 2b is immediately preceded by an instance of Case 2, and so the segments $J_i$ have been defined (possibly redefined) so that they partition $V(P'_1)$. Call a segment $J_i$ \emph{broken} if one of its edges was used as a broken edge in one of the rotations in Case 2 (or, for later uses, in Case 2b), and let $I' \subseteq I$ be the set of all indices $i \in I$ such that $J_i$ is an unbroken segment. Note that since each rotation breaks at most one segment, $|\bigcup_{i \in I'} Int_{P'_1}(J'_i)| \geq (1 - 2/15 - \tau_0/\log n)|P'_1| > (1 - 3/20)|P'_1|$. Hence for some $i \in I'$, $J'_i$ contains a $(H'', Q', \tau_0)$-reachable vertex of $P'_1$. Perform the necessary rotations to make one of these vertices into the first endpoint, breaking at most $\tau_0$ edges in the process. Call the resulting path $P_2 = x_2 \ldots y_2$. \medskip

Eventually, the iterative procedure has to conclude by entering Case~2b.
Let $I'' \subseteq I$ be the set of all indices $i \in I$ such that $J_i$ was broken in neither Case 2 nor Case 2b. Let $I_1$, $I_2$ be subsets of $I''$ such that $|I_1| = |I_2| = |I''|/2$ and for every $i_1 \in I_1$, $i_2 \in I_2$ the segment $J_{i_1}$ precedes $J'_{i_2}$ where the ordering is taken along $P_2$. Now
$$|I_1|, |I_2| \geq \frac{1}{2}((1-1/15)\log n - 2\tau_0) \geq \frac{1}{2}(1-1/10)\log n.$$ 
Let $Q_1 = \{x_2\} \cup \bigcup_{i \in I_1} J'_i$, $Q_2 = \{y_2\} \cup \bigcup_{i \in I_2} J'_i$ and note that
\begin{equation} \label{SizeOfQ12}
|Int_{P_2}(Q_1)|, |Int_{P_2}(Q_2)| \geq \frac{1}{2}\left(1 - \frac{1}{10}\right)\left(1 - \frac{1}{15}\right)|P_1| \geq \frac{2|P_1|}{5} \geq \frac{n}{25}.
\end{equation}

Further, for each vertex $v$ of $Q_1$, and for at least $|I_1| - \log n/15 \geq \log n/3$ indices $i \in I_1$, we have that $|N_{H''}(v) \cap Int_{P_2}(J'_i)| \geq \frac{r_G}{25 \log n}$ (noting that $Int_{P_2}(J'_i) = Int_{P_1}(J'_i)$ since $J'_i$ is unbroken). Hence
$$|N_{H''}(v) \cap Int_{P_2}(Q_1)| \geq r_G/75$$
for all $v \in Q_1$, and the corresponding statement holds for vertices of $Q_2$.

Now applying Corollary \ref{RotExp} with $P = P_2$, $Q = Q_1$, $Q' = Int_{P_2}(Q_1)$, $H' = H''$ and $\varepsilon = 1/75$ implies that there exists a set $A$ of $|Int_{P_2}(Q_1)|/6$ $(H'', Int_{P_2}(Q_1), \tau_0)$-reachable vertices of $P_2$ (since $Q' \subseteq V(P_2)$, it is impossible for the first conclusion of Corollary \ref{RotExp} to hold). Let $P'_2$ be the reverse of $P_2$ and apply Corollary \ref{RotExp} again with $P = P'_2$, $Q = Q_2$ and $Q' = Int_{P'_2}(Q_2)$ to obtain a set $B$ of $|Int_{P'_2}(Q_2)|/6$ $(H'', Int_{P'_2}(Q_2), \tau_0)$-reachable vertices of $P'_2$. Now $|A|, |B| \geq n/150$ by (\ref{SizeOfQ12}), and so using (\ref{DisjointJumbled}) we have
$$e_{G}(A, B) \geq \frac{|A||B|r_{G}}{n} - 4\sqrt{r_{G}}(|A| + |B|) \geq \frac{nr_{G}}{22500} - 4n\sqrt{r_{G}}\geq \frac{nr_{G}}{30000}.$$
Hence
$$e_{H''}(A, B) \geq \frac{nr_{G}}{30000} - |E(G)\backslash E(H'')| \stackrel{(\ref{EG'EH})}{>} 0.$$

Now let $x_3 y_3 \in E_{H''}(A, B)$. Noting that $x_3$ is $(H'', Int_{P_2}(Q_1), \tau_0)$-reachable, perform a set of at most $\tau_0$ rotations of $P_2$ with pivots in $Int_{P_2}(Q_1)$, to form a new path $P''_2$ whose first endpoint is $x_3$. Note that $Int_{P''_2}(Q_2) = Int_{P_2}(Q_2)$. Thus setting $P'''_2$ to be the reverse of $P''_2$, we have that $y_3$ is still $(H'', Int_{P'''_2}(Q_2), \tau_0)$-reachable. Perform a set of at most $\tau_0$ rotations to make $y_3$ the first endpoint of $P'''_2$. Note that all the pivots of these rotations lie in $Q_1 \cup Q_2 \subseteq \bigcup_{i \in I} J'_i \subseteq Q'$, and so as discussed after (\ref{PresInt1}) no edge in $E_{bad}$ is broken. Now use $x_3 y_3$ to close $P'''_2$ to a cycle. 

It remains to estimate the number of edges broken during the process. In the initial formation of $P_1$ we break only $4$ edges (two while forming $P_0$ and one for each of the two subsequent rotations or extensions). During each instance of Case 1 or Case 2a we break at most $\tau_0 + 1$ edges, during Case 2 we break at most $\tau_0$ edges, during Case 2b we break at most $\tau_0$ edges, and during the remainder of the proof we break at most $2\tau_0$ edges. Each time we merge a new cycle into $P_1$, we need either to go through Case 1 once, or through both Cases 2 and 2a once. Case 2b occurs only once (at the end, after Case 2).

Let $k_1$ be the number of times we repeat Case 1 and $k_2$ be the number of times we repeat Cases 2 and 2a. Then the number of cycles of $F$ decreases by $k_1 + k_2 + 1$. Further the total number of edges broken is at most 
$$4 + k_1(\tau_0 + 1) + k_2(2\tau_0 + 1) + 4\tau_0 \leq 4(k_1 + k_2 + 1)\tau_0 + 4$$
So on average the number of edges broken per cycle merged is at most $4\tau_0 + 4 \leq 5 \log n/\log \log n$. \endproof

We can now apply Lemma~\ref{RotExt} repeatedly to transform $2$-factors into Hamilton cycles. Given a $2$-regular graph $F$ and an edge-disjoint graph $H$ on the same vertex set, we say that a Hamilton cycle $C$ is formed by \emph{merging} the cycles of $F$ using edges of $H$ if $E(C) \subseteq E(F) \cup E(H)$. $E(H) \cap E(C)$ is the set of \emph{used edges}, and $E(F) \backslash E(C)$ is the set of \emph{broken edges}.

\begin{cor} \label{MergeCycles} Let $G'$ be an $r_{G'}/n$-pseudorandom graph on $n$ vertices, and let $G$ be an $r_{G}/n$-pseudorandom spanning subgraph of $G'$ with $300 \log^3 n \leq r_G \leq r_{G'}$. Let $H$ be an even-regular spanning subgraph of $G$ with degree $r_H$, such that $\delta(G)-1 \leq r_H \leq \delta(G)$. Let $\mathcal{F}$ be a collection of edge-disjoint $2$-factors $F_1, F_2, \ldots, F_m$ of $G'$, such that each $F_i$ is edge-disjoint from $H$.

Let $E_{bad} \subseteq E(\bigcup \mathcal{F})$ be such that $E_{bad} \cap E(F_i)$ is a matching and $|E_{bad} \cap E(F_i)| \leq n/10^6$ for each $F_i$. Suppose that $c(\mathcal{F}) r_{G'} \log^3 n \leq 4n r_{G}^2$. Then we can merge the cycles of each $F_i$ into a Hamilton cycle $C_i$ using the edges of $H$. Further, we can ensure that the number of edges of $E_{bad}$ broken during the process is at most $|E_{bad}|/\log n$, and that all of the $C_i$'s are pairwise edge-disjoint. 
\end{cor}

When $r_G$ is much smaller than $r_{G'}$ the bound on $c(\mathcal{F})$ is more restrictive than simply requiring that $c(\mathcal{F})$ is small compared to $e(H)$ as one might at first expect. The assumption is necessary due to (\ref{RotExtEdgeBound}) which in turn arises from (\ref{AltExpBound}). In the proof of Lemma \ref{MiniHamDecomp} (and thus of Theorem \ref{PRproof}) this assumption will be the limiting factor in determining how small we can make $r_G$ compared to $r_{G'}$, and thus how many iterations we need to use.

\proof We merge cycles by repeatedly applying Lemma~\ref{RotExt}. During this process we will remove certain edges from $H$ (namely those which lie in the new $2$-factor obtained by Lemma~\ref{RotExt}) and add certain edges to $H$ (namely those edges which are removed from the old $2$-factor in Lemma~\ref{RotExt} to obtain the new one). Let $H'$ denote the `current' version of $H$ (so $H$ always denotes the original version).

We use Lemma~\ref{RotExt} repeatedly to reduce $c(F_i)$ until $c(F_i) = 1$, i.e., $F_i$ is a Hamilton cycle for each $i$. We make use of the fact that on average at most $5 \log n/\log \log n$ edges of $H'$ are used by Lemma~\ref{RotExt} for each cycle that needs to be merged. Hence we have 
$$|E(H)\backslash E(H')| = |E(H')\backslash E(H)| \leq \frac{5 c(\mathcal{F}) \log n}{\log \log n} \leq \frac{r_{G}^2 n}{2500 r_{G'} \log^2 n}$$
throughout the process. So (\ref{RotExtEdgeBound}) is satisfied.

More precisely, we proceed as follows: Given $1 \leq i \leq m$, suppose that the number of cycles of $F_i$ is at least $2|E_{bad} \cap E(F_i)|/\log^2 n$. Then we can apply Lemma~\ref{RotExt} with $E_{bad} = E_{bad} \cap E(F_i)$ to merge one or more cycles of $F_i$ together, without breaking any edges of $E_{bad}$ and using on average at most $5\log n/\log \log n$ edges of $H'$ per cycle merged. (Note that the required upper bound on $c(F_i)$ holds since $c(F_i) \leq c(\mathcal{F})$.)

On the other hand if the number of cycles is fewer than $2|E_{bad} \cap E(F_i)|/\log^2 n$, then we can apply Lemma~\ref{RotExt} with $E_{bad} = \emptyset$ repeatedly to merge all of the remaining cycles of $F_i$ together, thus forming a Hamilton cycle. Again we use at most $5\log n/\log \log n$ edges of $H'$ per cycle and thus the total number of edges of $E_{bad}$ which are broken during this process is at most 
$$\frac{2|E_{bad} \cap E(F_i)|}{\log^2 n} \cdot \frac{5 \log n}{\log \log n} \leq \frac{|E_{bad} \cap E(F_i)|}{\log n}.$$
The bound on the total number of edges of $E_{bad}$ broken follows immediately. \endproof

For a given application of Corollary \ref{MergeCycles} during the remainder of the paper, define the \emph{leftover graph} to be the final state of $H'$. Thus $H'$ is obtained from $H$ by deleting the edges of the Hamilton cycles produced by Corollary \ref{MergeCycles}, and adding those edges of $\bigcup \mathcal{F}$ which do not lie in any of these Hamilton cycles.

We now prove variants of Lemma~\ref{RotExt} and Corollary \ref{MergeCycles}. These will be used in cases where instead of needing to avoid a set of bad edges when rotating, we only want to avoid a vertex $x_0$ (when we apply the lemmas $x_0$ will be the vertex of minimum degree in $G_{n, p}$). For clarity we prove these as separate lemmas as the conditions are different, but the proofs proceed along similar lines in each case.

\begin{lemma} \label{FinalMerge1} Let $G'$ be an $r_{G'}/n$-pseudorandom graph on $n$ vertices, and let $G$ be an $r_{G}/n$-pseudorandom spanning subgraph of $G'$ with $10^5 \log^2 n \leq r_{G} \leq r_{G'}$.\COMMENT{This suffices here, whereas in Lemma \ref{MergeCycles} we needed a better bound in order to apply Lemma \ref{CreateW}.} Let $x_0$ be a vertex of $G$ and let $H$ be formed from $G$ by removing all of the edges of $G$ incident to $x_0$.  Let $H'$ be a  spanning subgraph of $G'$, such that (\ref{RotExtEdgeBound}) holds and $d_H(v) = d_{H'}(v)$ for each $v \in V(G)$.

Let $F$ be a $2$-factor of $G'$, edge-disjoint from $H'$, such that $c(F) \leq \frac{4nr_G^2}{r_{G'} \log^3 n}$. Then unless $F$ is a Hamilton cycle, we can obtain a new $2$-factor $F'$ such that the following hold:
\begin{itemize}
\item $E(F') \subseteq E(F) \cup E(H')$,
\item $c(F') < c(F)$, and
\item $|E(F') \cap E(H')| \leq \frac{5 \log n}{\log \log n}(c(F) - c(F'))$.
\end{itemize}
\end{lemma}

\proof The proof is similar to that of Lemma~\ref{RotExt}. We will again obtain $F'$ from $F$ by performing a sequence of rotations and extensions, with all the new edges for these rotations and extensions taken from a graph $H''$ which is defined as in the proof of Lemma~\ref{RotExt}. Again, we say that an edge $e$ of $H''$ is \emph{used} if it is the new edge of some rotation or extension performed during the proof. Define the current version $F'$ of $F$ as in the proof of Lemma~\ref{RotExt}. Observe that $e(H) \geq e(G) - n$, and so assuming we use on average at most $\frac{5 \log n}{\log \log n}$ edges for each cycle of $F$ merged, the bounds (\ref{EG'EH1}), (\ref{EG'EH2}) and (\ref{EG'EH}) still hold. Moreover, 
\begin{align} \label{deltaHx0}
\delta(H'' - \{x_0\}) &\geq \delta(H' - \{x_0\}) = \delta(H - \{x_0\}) \nonumber \\
&\geq \delta(G) - 1 \geq r_G - 2\sqrt{r_G \log n} - 1 \geq \frac{r_G}{2}.
\end{align}

\noindent \textbf{Claim:} \emph{$H' - \{x_0\}$ is connected.} \medskip

\noindent To prove the claim, suppose for a contradiction that $H' - \{x_0\}$ has two components, $S$ and $T$. By (\ref{deltaHx0}) we can apply Lemma~\ref{AltExp} with $\varepsilon = 1/3$, $Q' = V(G) \backslash\{x_0\}$ and $S = S$. But now $N_{H' - \{x_0\}}(S) = \emptyset$, and so of the possible conclusions of Lemma~\ref{AltExp} only (iv) can hold. This implies that $|S| \geq (n-1)/6 > n/7$. Similarly $|T| \geq (n-1)/6 > n/7$. Now since $G$ is $(r_{G}/n, 2\sqrt{r_{G}})$-jumbled we have that 
$$e_{G}(S, T) \geq r_G|S||T|/n - 4\sqrt{r_{G}}(|S| + |T|) \geq nr_{G}/50.$$ 
But by (\ref{EG'EH}) this implies that $e_{H'}(S, T) > 0$, which contradicts our assumption that $S$ and $T$ were components. Hence $H' - \{x_0\}$ has only one component, which proves the claim. \medskip

Let $C_1$ be a cycle of $F$. Since $H' - \{x_0\}$ is connected, there exists an edge $xy$ of $H'$ joining two distinct cycles $C_1$ and $C_2$ of $F$, with $x$ on $C_1$ and $y$ on $C_2$ and such that $x, y \neq x_0$. Let $yy'$ be an edge of $C_2$ incident to $y$ such that $y' \neq x_0$. Delete $yy'$ from $C_2$ to form a path $C'_2$. Let $xx'$ be an edge of $C_1$ such that $x' \neq x_0$, and perform an extension of $C'_2$ with join vertex $x$ and broken edge $xx'$ to incorporate $C_1$. Let $P_1 = x' \ldots y'$ denote the resulting path. Note that $x_0$ cannot be an endpoint of $P_1$.

Let $V'' = V(G) \backslash \{x_0\}$ and $Q' = Int_F(V'')$. Note that $Int_{P_1}(V'') = Q' \cap V(P_1)$. All rotations performed during the remainder of the proof will have pivots in $Int_F(V'') \cap V(P_1)$, where $P_1$ is the current path. Similarly, all extensions performed during the remainder of the proof will have join vertices in $Q'$. Thus by applying Lemma~\ref{PreserveNhood} with $Q = V''$, we always have that
$$Int_{P_1}(V'') = Int_F(V'') \cap V(P_1) = Q' \cap V(P_1).$$
So in particular $x_0$ will never be an endpoint of $P_1$.

By (\ref{deltaHx0}) we can apply Corollary \ref{RotExp} with $\varepsilon = 1/3$, $H' = H''$, $P = P_1$, $Q = V'' \cap V(P_1)$ and $Q' = Q'$ and use the same case analysis as in Lemma~\ref{RotExt}. The remainder of the argument is also identical to that in the proof of Lemma~\ref{RotExt}. (Note that $x_0 \notin J'_i$ for any $i$ since $J'_i \subseteq J_i \cap Q'$, and hence $x_0 \notin Q_1$ and $x_0 \notin Q_2$.) \endproof

\begin{cor} \label{FinalMerge2} Let $G'$ be an $r_{G'}/n$-pseudorandom graph on $n$ vertices, and let $G$ be an $r_{G}/n$-pseudorandom spanning subgraph of $G'$ with $10^5 \log^2 n \leq r_{G} \leq r_{G'}$. Let $x_0$ be a vertex of $G$ and let $H$ be formed from $G$ by removing all of the edges of $G$ incident to $x_0$. Let $\mathcal{F}$ be a collection of edge-disjoint $2$-factors $F_1, F_2, \ldots, F_m$ of $G'$, such that each $F_i$ is edge-disjoint from $H$ and $c(\mathcal{F}) \leq \frac{4n r_{G}^2}{r_{G'} \log^3 n}$. Then we can merge the cycles of each $F_i$ into a Hamilton cycle $C_i$ using the edges of $H$, such that the $C_i$'s are pairwise edge-disjoint. (Recall that merging was defined in the paragraph before Corollary \ref{MergeCycles}.)
\end{cor}

\proof We merge cycles by repeatedly applying Lemma~\ref{FinalMerge1}. During this process we will remove certain edges from $H$ (namely those which lie in the new $2$-factor obtained by Lemma~\ref{FinalMerge1}) and add certain edges to $H$ (namely those edges which are removed from the old $2$-factor in Lemma~\ref{FinalMerge1} to obtain the new one). Let $H'$ denote the `current' version of $H$ (so $H$ always denotes the original version).

We use Lemma~\ref{FinalMerge1} repeatedly to reduce $c(F_i)$ until $c(F_i) = 1$, i.e., $F_i$ is a Hamilton cycle for each $i$. We make use of the fact that on average at most $5 \log n/\log \log n$ edges of $H'$ are used by Lemma~\ref{FinalMerge1} for each cycle that needs to be merged. So
$$|E(H)\backslash E(H')| = |E(H')\backslash E(H)| \leq \frac{5 c(\mathcal{F}) \log n}{\log \log n} \leq \frac{r_{G}^2 n}{2500 r_{G'} \log^2 n},$$
and hence (\ref{RotExtEdgeBound}) is satisfied throughout the process. \endproof

\section{Completing the proof} \label{Completing}

In this section we combine our results to prove Theorem~\ref{HamDecomp}. Roughly speaking, the following lemma states that given a graph $H_0$ which is close to being pseudorandom and given about $\log n$ pseudorandom graphs $H_1, \ldots, H_{2m+1}$, we can find a set of edge-disjoint Hamilton cycles in the union of these graphs which cover all of the edges of $H_0$. While the $H_i$ cannot be too sparse, they need not be as dense as $H_0$. The point is that the remaining `uncovered' graph is much sparser than $H_0$ and is also not too far from being pseudorandom. In the proof of Theorem~\ref{PRproof} we apply this lemma three times in succession to obtain an uncovered graph which is very sparse.

\begin{lemma} \label{MiniHamDecomp} Let $p_0 \geq \log^{14} n/n$, and let $p_1 \geq ((np_0)^3 \log^{10} n)^{1/4}/n$. Let $m = \frac{\log(n^2 p_1)}{\log \log n}$ and let $p_2, \ldots, p_{2m+1}$ be positive reals such that $p_i = p_1$ for odd $i$ and $p_i = 10^{10} p_1$ for even $i$.\COMMENT{The fact that $m \leq \frac{2 \log n}{\log \log n}$ implies that $\sum_{i = 1}^{2m+1} p_i \leq p_0$. So the lemma is not vacuous.} Let $G_0$ be a $p_0$-pseudorandom graph on $n$ vertices. Suppose that $G_i$ is a $p_i$-pseudorandom spanning subgraph of $G_0$, and let $H_i$ be an even-regular spanning subgraph of $G_i$ such that $\delta(G_i) - 1 \leq \delta(H_i) \leq \delta(G_i)$, for each $1 \leq i \leq 2m + 1$. Suppose that the graphs $G_i$ are pairwise edge-disjoint for $1 \leq i \leq 2m+1$ and let $H_0$ be an even-regular spanning subgraph of $G_0$ which is edge-disjoint from $\bigcup_{i = 1}^{2m+1} H_i$. Then there exists a collection $\mathcal{HC}$ of edge-disjoint Hamilton cycles such that $H_0 \subseteq \bigcup \mathcal{HC} \subseteq \bigcup_{i = 0}^{2m+1} H_i$.
\end{lemma}

Formally the assumption $p_0 \geq \log^{14} n/n$ can be omitted. It is included for clarity since if $p_0$ is significantly smaller then $p_0 \leq \sum_{i=1}^{2m+1} p_i$ and so the lemma becomes vacuous.

To prove Lemma~\ref{MiniHamDecomp} we first decompose $H_0$ into $2$-factors which on average have few cycles. We then use edges of $H_1$ to transform these $2$-factors into Hamilton cycles. Because edges are exchanged between $H_1$ and the $2$-factors there will still be some edges of $H_0$ left uncovered (the `bad' edges). We decompose $H'_1 \cup H_2$ (where $H'_1$ is the leftover of $H_1$ and $H_0$) into $2$-factors with few cycles and then use $H_3$ to transform them into Hamilton cycles (we cannot decompose $H'_1$ on its own since it is no longer close to being pseudorandom). Again some edges of $H_0$ will be left uncovered, but we can guarantee that the number of such edges will be reduced (by a factor of about $\log n$). After about $\log n/\log \log n$ iterations we arrive at a leftover graph which contains no edges of $H_0$, i.e., all of the edges of $H_0$ are covered.

\proof[Proof of Lemma \ref{MiniHamDecomp}] Note that $p_1 \geq \log^{13} n/n$. Corollary \ref{2FacSplit} with $G = G_0$ and $H = H_0$ implies that $H_0$ can be decomposed into a collection $\mathcal{F}_1$ of $2$-factors such that $c(\mathcal{F}_1) \leq 3n\sqrt{np_0 \log^3 n}$. Since
\begin{equation} \label{cFcalc}
c(\mathcal{F}_1) \cdot np_0 \log^3 n \leq 4n (np_0 \log^3 n)^{3/2} \leq 4n(np_1)^2,
\end{equation}
we may apply Corollary \ref{MergeCycles} with $G' = G_0$, $G = G_1$, $H = H_1$, $\mathcal{F} = \mathcal{F}_1$ and $E_{bad} = \emptyset$. This allows us to merge the cycles of each $2$-factor of $\mathcal{F}_1$ to form a collection of edge-disjoint Hamilton cycles. Let $H'_1 \subseteq H_0 \cup H_1$ be the leftover graph from Corollary \ref{MergeCycles} (as defined after the proof of Corollary \ref{MergeCycles}), and let $E_{bad} = E(H'_1) \cap E(H_0)$. Note that $|E_{bad}| \leq |E(H'_1)| = |E(H_1)| \leq n^2 p_1/2$.

Now $H'_1$ is even-regular with degree at most $np_1$, and $np_1 + 1 + 10^6np_1/2 \leq np_2/5000$. Hence we may apply Corollary \ref{Absorb2Fac} with $G = G_2$, $H = H_2$ and $H' = H'_1$ to obtain a decomposition $\mathcal{F}_2$ of $H'_1 \cup H_2$ into $2$-factors, such that $c(\mathcal{F}_2) \leq 4n\sqrt{np_2 \log^3 n}$ and $E(F) \cap E_{bad}$ is a matching of size at most $n/10^6$ for each $F \in \mathcal{F}_2$. 

By (\ref{cFcalc}) (noting that $p_1 = p_3$) we may apply Corollary \ref{MergeCycles} with $G' = G_0$, $G = G_3$, $H = H_3$, $\mathcal{F} = \mathcal{F}_2$ and $E_{bad} = E_{bad}$. This yields another collection of edge-disjoint Hamilton cycles and a leftover graph $H'_3$. Furthermore, if we redefine $E_{bad} = E(H'_3) \cap E(H_0)$, then $|E_{bad}|$ is reduced by a factor of $\log n$.

We now repeat this process a further $m-1$ times, using up the graphs $H_4, H_5, \ldots,$ $ H_{2m+1}$. Now we have $|E(H'_{2m+1}) \cap E(H_0)| \leq \frac{n^2p_1}{2(\log n)^m} < 1$, i.e., $E(H'_{2m+1}) \cap E(H_0) = \emptyset$. Let $\mathcal{HC}$ be the union of all the collections of Hamilton cycles produced by this process and note that $H_0 \subseteq \bigcup \mathcal{HC} \subseteq \bigcup_{i=0}^{2m + 1} H_i$. 
\endproof

Before proving Theorem~\ref{HamDecomp}, we state a `pseudorandom' version of the theorem. The conditions in this version are significantly more complicated than those of Theorem~\ref{HamDecomp}; however, it has the advantage of being entirely deterministic and we believe it to be of independent interest.

\begin{thm} \label{PRproof} Let $\log^{50} n/n \leq p_0 \leq 1 - n^{-1/4}\log^9 n$. Let 
\begin{align} \label{p1p2p3}
&p_2 = \frac{(np_0)^{3/4} \log^{7/2} n}{n},~p_3 = \frac{(np_2)^{3/4} \log^{7/2} n}{n}, \nonumber \\
&p_4 = \frac{(np_3)^{3/4} \log^{7/2} n}{n},~p_5 = \frac{\sqrt{np_0(1-p_0)}}{n \log^{3/4} n}
\end{align}
and $p_1 = p_0 - p_2 - p_3 - p_4 - p_5$. For $i = 2, 3, 4$ let 
\begin{equation*} \label{mi}
m_i = \frac{2\log(n^2 p_i)}{\log \log n}.
\end{equation*}
For $1 \leq j \leq 2m_i + 1$ let $p_{(i, j)} = \frac{p_i}{(10^{10} + 1)m_i + 1}$ if $j$ is odd and $p_{(i, j)} = \frac{10^{10}p_i}{(10^{10} + 1)m_i + 1}$ if $j$ is even. 

Let $G_0$ be a $p_0$-pseudorandom graph on $n$ vertices. Suppose that $G_0$ has a decomposition into graphs $G_1, G_2, G_3, G_4, G_5$, and that $G_i$ has a decomposition into graphs $G_{(i, 1)}, G_{(i, 2)}, \ldots, G_{(i, 2m_i + 1)}$ for $i = 2, 3, 4$, such that the following conditions hold:
\begin{itemize}
\item [\rm (i)] $G_i$ is $p_i$-pseudorandom for each $1 \leq i \leq 5$,
\item [\rm (ii)] $G_{(i, j)}$ is $p_{(i, j)}$-pseudorandom for $i = 2, 3, 4$ and $1 \leq j \leq 2m_i+1$,
\item [\rm (iii)] $G_i \cup G_{i + 1}$ is $(p_i + p_{i+1})$-pseudorandom for $i = 2, 3, 4$, and
\item [\rm (iv)] $G_0$ is $4u$-jumping where 
\begin{equation} \label{definu}
u = 2np_5 = \frac{2 \sqrt{np_0(1-p_0)}}{\log^{3/4} n}.
\end{equation}
\end{itemize}
Then $G_0$ has property $\mathcal{H}$.
\end{thm}

For the moment we will assume the truth of Theorem~\ref{PRproof} and use it to prove Theorem~\ref{HamDecomp}.

\noindent \proof[Proof of Theorem~\ref{HamDecomp}] Let $p_0 = p$ and let $G_0 \sim G_{n, p}$. Define as in the statement of Theorem~\ref{PRproof} the real numbers $p_i$ for $1 \leq i \leq 5$, integers $m_i$ for $i = 2, 3, 4$, and reals $p_{(i, j)}$ for $i = 2, 3, 4$ and $1 \leq j \leq 2m_i + 1$.

Form graphs $G_1, G_5$ and $G_{(i, j)}$ for $i = 2, 3, 4$ and $1 \leq j \leq 2m_i + 1$ as follows: For each edge $e$ of $G_0$, place $e$ in $G_1$ with probability $p_1/p_0$, in $G_5$ with probability $p_5/p_0$, and in $G_{(i, j)}$ with probability $p_{(i, j)}/p_0$ for each $i = 2, 3, 4$ and $1 \leq j \leq 2m_i + 1$. Let $G_i = \bigcup_{j = 1}^{2m_i + 1} G_{(i, j)}$ for $i = 2, 3, 4$.

Note that $p_i = o(p_0)$ for $i = 2, 3, 4, 5$ and that
\begin{equation} \label{piconditions}
np_5 \geq \log^{24} n,~np_i \geq \log^{14} n \text{ and } np_{(i, j)} \geq \log^{13} n
\end{equation}
for $i = 2, 3, 4$, where the second inequality holds since $x \geq \log^{14} n$ implies that $x^{3/4} \log^{7/2} n \geq \log^{14} n$. Thus the bounds on each $p_i$ and $p_{(i, j)}$ in Lemma \ref{SplitGraph} hold, and hence Lemma \ref{SplitGraph} implies that $G_0$ is $p_0$-pseudorandom and that conditions (i), (ii) and (iii) of Theorem~\ref{PRproof} hold whp. Moreover, condition (iv) holds whp by Lemma~\ref{DegJumpProb}. Hence Theorem~\ref{PRproof} implies that $G_0$ has property $\mathcal{H}$. \endproof

It remains to prove Theorem~\ref{PRproof}.

\proof[Proof of Theorem~\ref{PRproof}] Recall from the proof of Theorem \ref{HamDecomp} that $p_i = o(p_0)$ for $i = 2, 3, 4, 5$. Thus
\begin{equation} \label{uapprox}
p_1 = (1 - o(1))p_0.
\end{equation}
Note also that (\ref{piconditions}) holds and 
\begin{equation} \label{np3bound}
np_4 = (np_0)^{27/64} \log^{\frac{7}{2}\left(1 + \frac{3}{4} + \frac{9}{16}\right)} \leq (np_0)^{27/64} \log^{49/6} n.
\end{equation}

Let $x_0$ be the vertex of $G_0$ of minimum degree. If $\delta(G_0)$ is odd, then at this point we use Lemma~\ref{Tutte}(i) to remove an optimal matching $M_{Opt}$ which covers $x_0$ from $G_1$, and let $G'_1$ be the remainder. If $\delta(G_0)$ is even then let $G'_1 = G_1$ and $M_{Opt} = \emptyset$.

Form $H_5$ from $G_5$ by removing all of the edges incident to $x_0$, and add the removed edges to $G'_1$. For each $G_{(i, j)}$, apply Lemma~\ref{Tutte}(ii) with $u = 0$ to form a regular spanning subgraph $H_{(i, j)}$ whose degree is either $\delta(G_{(i, j)})$ (if $\delta(G_{(i, j)})$ is even) or $\delta(G_{(i, j)}) - 1$ (otherwise). If there are edges of $G_{(i, j)} \backslash H_{(i, j)}$ which are incident to $x_0$, move these edges into $G'_1$. Let $H_i = \bigcup_{j = 1}^{2m_i + 1} H_{(i, j)}$ for $i = 2, 3, 4$. Now all edges of $G_0$ which are incident to $x_0$ lie in $G'_1 \cup H_2 \cup H_3 \cup H_4 \cup M_{Opt}$, i.e.,
\begin{equation} \label{x0edges}
d_{G_0}(x_0) = d_{G'_1}(x_0) + d_{H_2}(x_0) + d_{H_3}(x_0) + d_{H_4}(x_0) + d_{M_{Opt}}(x_0).
\end{equation}

Let $\sum_{i, j}$ denote the summation $\sum_{i = 2}^4 \sum_{j=1}^{2m_i + 1}$. Note that (\ref{x0edges}) implies that
\begin{equation} \label{sumij}
d_{G_1}(x_0) - 1 \leq d_{G'_1}(x_0) \leq \delta(G_0) - 2 \sum_{i, j} \lfloor \delta(G_{(i, j)})/2 \rfloor \leq \delta(G_0) - \sum_{i, j} (\delta(G_{(i, j)}) - 1).
\end{equation}
The next claim shows that the number of edges incident to $x_0$ which we added to $G_1$ to form $G'_1$ is at most $u$. \medskip

\noindent \textbf{Claim 1:}
\begin{equation} \label{DegreeDiff}
\Delta(G_5) + \sum_{i, j}\left(\Delta(G_{(i, j)}) - \delta(G_{(i, j)}) + 1\right) \leq u.
\end{equation}
Indeed, by Definition \ref{PseudoRandom}(c) we have that 
$$\Delta(G_5) \leq np_5 + 2\sqrt{np_5\log n} \leq \frac{4np_5}{3} \stackrel{(\ref{definu})}{=} \frac{2u}{3}.$$ 
Further, note that 
\begin{equation} \label{p0bound}
(5 \log^{4} n)^8 \leq np_0(1-p_0)^4.
\end{equation}
Also, Definition \ref{PseudoRandom}(c) implies that 
$$\Delta(G_{(2, j)}) - \delta(G_{(2, j)}) \leq 4\sqrt{np_{(2, j)} \log n} \leq 4\sqrt{np_2 \log n}$$ 
for each $1 \leq j \leq 2m_2+1$. Hence using $m_2 = o(\log n)$, we have
\begin{align*}
\sum_{j=1}^{2m_2 + 1} \left(\Delta(G_{(2, j)}) - \delta(G_{(2, j)}) + 1\right) &\leq 4(2m_2 + 1)\left(\sqrt{np_2 \log n} + 1\right) \\ 
&\leq \log n \sqrt{np_2 \log n} = \log^{13/4} n (np_0)^{3/8} \\
&\stackrel{(\ref{p0bound})}{\leq} \frac{\sqrt{np_0(1-p_0)}}{5\log^{3/4} n} = \frac{np_5}{5} \stackrel{(\ref{definu})}{=} \frac{u}{10}.
\end{align*}
Similarly 
$$\sum_{j=1}^{2m_i + 1} \left(\Delta(G_{(i, j)}) - \delta(G_{(i, j)}) + 1\right) \leq \log n \sqrt{np_i \log n} \leq \frac{u}{10}$$ 
for $i = 3, 4$. Thus the left-hand side of (\ref{DegreeDiff}) is at most $2u/3 + 3(u/10) \leq u$, which proves the claim. \medskip

\noindent \textbf{Claim 2:} \emph{Each $x \neq x_0$ satisfies $d_{G_1}(x) \geq d_{G_1}(x_0) + 2u$.} \medskip

\noindent In other words, $x_0$ is the vertex of minimum degree in $G_1$ and $G_1$ is $2u$-jumping. To prove Claim 2, recall our assumption that $G_0$ is $4u$-jumping. Hence we have
\begin{align*}
d_{G_1}(x) &\geq \delta(G_0) + 4u - \Delta(G_5) - \sum_{i, j} \Delta(G_{(i, j)}) \\
&\stackrel{(\ref{DegreeDiff})}{\geq} \delta(G_0) - \sum_{i, j} (\delta(G_{(i, j)}) - 1) + 3u \stackrel{(\ref{sumij})}{\geq} d_{G_1}(x_0) - 1 + 3u \geq d_{G_1}(x_0) + 2u,
\end{align*}
which proves the claim. \medskip

Since each $H_{(i, j)}$ is even-regular and since $M_{Opt}$ covers $x_0$ if and only if $d_{G_0}(x_0)$ is odd, ({\ref{x0edges}) implies that $d_{G'_1}(x_0)$ is even. Also recall that the number of edges added to $G_1$ at the vertex $x_0$ (after removing $M_{Opt}$) to form $G'_1$ is at most the left-hand side of (\ref{DegreeDiff}), and hence is at most $u$. So $x_0$ is the vertex of minimum degree in $G'_1$. Moreover, since $1 - p_0 < 1 - p_1$,
$$u \stackrel{(\ref{definu}), (\ref{uapprox})}{\leq} \frac{2\sqrt{n\cdot 2p_1 (1 - p_1)}}{\log^{3/4} n} \leq 4\sqrt{np_1(1-p_1)}.$$ 
Hence we may apply Lemma~\ref{Tutte} with $G = G_1$ and $G' = G'_1$ to form a regular spanning subgraph $H_1$ of $G'_1$ with degree $\delta(G'_1) = d_{G'_1}(x_0)$. Note that $H_1$ contains every edge of $G'_1$ incident to $x_0$.

By Lemma~\ref{MiniHamDecomp} where $p_0 = p_0$, $p_j = p_{(2, j)}$ for $1 \leq j \leq 2m_2 +1$, $G_0 = G_0$, $H_0 = H_1$, $G_j = G_{(2, j)}$, and $H_j = H_{(2, j)}$ for $1 \leq j \leq 2m_2+1$, there exists a collection $\mathcal{HC}_1$ of edge-disjoint Hamilton cycles in $H_1 \cup H_2$ which cover the edges of $H_1$. Let $H'_2$ be the graph formed by the edges of $H_2$ which are not contained in one of these Hamilton cycles.

Applying Lemma~\ref{MiniHamDecomp} again with $p_0 = p_2 + p_3$, $p_j = p_{(3, j)}$ for $1 \leq j \leq 2m_3+1$, $G_0 = G_2 \cup G_3$, $H_0 = H'_2$, $G_j = G_{(3, j)}$, and $H_j = H_{(3, j)}$, we obtain a collection $\mathcal{HC}_2$ of edge-disjoint Hamilton cycles in $H'_2 \cup H_3$ which cover $H'_2$. Let $H'_3$ be the graph formed by the edges of $H_3$ which are not covered by one of the Hamilton cycles.

Applying Lemma~\ref{MiniHamDecomp} again with $p_0 = p_3 + p_4$, $p_j = p_{(4, j)}$ for $1 \leq j \leq 2m_4+1$, $G_0 = G_3 \cup G_4$, $H_0 = H'_3$, $G_j = G_{(4, j)}$, and $H_j = H_{(4, j)}$, we obtain a collection $\mathcal{HC}_3$ of edge-disjoint Hamilton cycles in $H'_3 \cup H_4$ which cover $H'_3$. Let $H'_4$ be the graph formed by the edges of $H_4$ which are not covered by one of the Hamilton cycles.

Note that $H'_4$ is a subgraph of $G_4$. Hence by Corollary \ref{2FacSplit} we can find a decomposition $\mathcal{F}$ of $H'_4$ into $2$-factors such that $c(\mathcal{F}) \leq 3n\sqrt{np_4\log^3 n}$. Now we claim that
\begin{equation} \label{EndMergeClaim}
36(np_4)^3 \log^9 n \leq (np_5)^4.
\end{equation}

To prove (\ref{EndMergeClaim}), note that $(np_5)^4 = (np_0(1-p_0))^2/\log^3 n$. So by (\ref{np3bound}) it suffices to prove that $36(np_0)^{81/64} \log^{67/2} n \leq (np_0(1-p_0))^2/\log^3 n$, or equivalently that $(np_0)^{47/64}(1-p_0)^2 \geq 36\log^{73/2} n$. But if $p_0 \leq 1/2$ then we have 
$$(np_0)^{47/64}(1-p_0)^2 \geq (np_0)^{47/64}/4 \geq (\log^{50} n)^{47/64}/4 \geq 36\log^{73/2} n,$$
and if $p_0 \geq 1/2$ then 
$$(np_0)^{47/64}(1-p_0)^2 \geq n^{47/64}(n^{-1/4})^2/2 \geq 36\log^{73/2} n,$$
with room to spare, which proves (\ref{EndMergeClaim}).\COMMENT{We don't actually need the last iteration but it improves the bounds.}

It follows immediately from (\ref{EndMergeClaim}) that if $p_5 \leq p_4$, then
$$c(\mathcal{F}) \leq 3n\sqrt{np_4\log^3 n} \leq \frac{n(np_5)^2}{2np_4\log^3 n} \leq \frac{n(np_5)^2}{n(p_4+p_5) \log^3 n}.$$
On the other hand if $p_5 \geq p_4$, then 
$$c(\mathcal{F}) \leq 3n\sqrt{np_4\log^3 n} \leq 3n\sqrt{np_5\log^3 n} \leq \frac{n(np_5)}{2\log^3 n} \leq \frac{n(np_5)^2}{n(p_4+p_5) \log^3 n}.$$
Hence in either case we can apply Corollary \ref{FinalMerge2} with $G' = G_4 \cup G_5$, $G = G_5$, $H = H_5$ and $\mathcal{F} = \mathcal{F}$ to obtain a collection $\mathcal{HC}_4$ of edge-disjoint Hamilton cycles in $H'_4 \cup H_5$ which cover $H'_4$. Now let $\mathcal{HC} = \mathcal{HC}_1 \cup \mathcal{HC}_2 \cup \mathcal{HC}_3 \cup \mathcal{HC}_4$. Observe that $\mathcal{HC}$ covers every edge of $H_1$, $H_2$, $H_3$ and $H_4$. But recall that every edge of $G_0$ incident to $x_0$ is contained in either $H_1$, $H_2$, $H_3$, $H_4$ or $M_{Opt}$. Hence $\mathcal{HC}$ contains exactly $\lfloor d_{G_0}(x_0)/2\rfloor = \lfloor\delta(G_0)/2\rfloor$ Hamilton cycles. \endproof

Note that the only place where we use the full strength of the condition on $p_0$ is in the proof of (\ref{p0bound}) and (\ref{EndMergeClaim}). Also note that if we omit one of the iterations (i.e., if instead of defining $p_4$, $G_4$ and $H_4$ we simply use $H_5$ to finish the decomposition of $H'_3$) then the proof of Theorem \ref{HamDecomp} still works as long as $\log^{125} n/n \leq p \leq 1 - n^{-1/7}$ (say). On the other hand, we could have improved the lower bound on $p$ in Theorem \ref{HamDecomp} somewhat by adding extra iterations. However, even a large number of iterations will only reduce the lower bound to approximately $\log^{30} n/n$. Some further small improvements could be made by using tighter calculations in some places.

\medskip

{\footnotesize \obeylines \parindent=0pt

\begin{tabular}{lll}

Fiachra Knox                        &\ &  Daniela K\"{u}hn \& Deryk Osthus \\
School of Mathematical Sciences			&\ &  School of Mathematics \\
Queen Mary, University of London  	&\ &  University of Birmingham \\
Mile End Road                       &\ &  Edgbaston \\
London															&\ &  Birmingham\\
E1 4NS															&\ &  B15 2TT \\
UK																	&\ &  UK\\
\end{tabular}

\begin{flushleft}
{\tt{E-mail addresses}:
f.knox@qmul.ac.uk, \{kuehn,osthus\}@maths.bham.ac.uk}
\end{flushleft} }

\end{document}